\documentclass[12pt]{article}
\usepackage{amsmath,amssymb}
\usepackage{times}
\usepackage{array}

\usepackage[pdftex]{graphicx}
\usepackage{epstopdf}

\newcommand{\abs}[1]{\lvert#1\rvert}
\newcommand{\norm}[1]{\lVert#1\rVert}

\def\div{\text{div}}

\input epsf

\def\Om{\Omega}
\def\qed{\hfill$\Box$}
\def\bR{\mathbb R}

\def\di{\hbox{div}}

\def\De{\Delta}

\def\cH{{\cal H}}

\def\cL{{\cal L}}

\newcommand{\sgn}{{\rm \,sgn\,}}
\renewcommand{\d}{\partial}
\newcommand{\R}{\mathbb R}

\numberwithin{equation}{section}

\newtheorem{proposition}{Proposition}[section]
\newtheorem{theorem}{Theorem}[section]
\newtheorem{lemma}{Lemma}[section]
\newtheorem{corollary}{Corollary}[section]
\begin{document}


\title{\bf Two cases of squares evolving by anisotropic diffusion}

\author{Piotr B. Mucha$^1$, Monika Muszkieta$^2$, Piotr Rybka$^1$}

\maketitle
\date{}

{\small 
\begin{center}
 1. Institute of Applied Mathematics and Mechanics, University of Warsaw\\
PL 02-097 Warszawa, Poland \\

\smallskip

2. Institute of Mathematics and Computer Science, Wroc{\l}aw University of
Technology, PL 50-370 Wroc{\l}aw, Poland\\

E-mails: p.mucha@mimuw.edu.pl, monika.muszkieta@pwr.wroc.pl, p.rybka@mimuw.edu.pl

\end{center}
}

\begin{abstract}
We are interested in an anisotropic singular diffusion equation in the plane and in its regularization. We establish
existence, uniqueness and basic regularity of solutions to both equations. We
construct explicit solutions showing the
creation of facets, i.e. flat regions of solutions. By using  the formula for
solutions, we  rigorously prove that both equations 
create ruled surfaces out of convex initial conditions as well as do not admit
point (local) extrema. 
We present numerical experiments suggesting that the two flows seem not differ
much. Possible applications to image reconstruction is pointed out, too.
\end{abstract}

\section{Introduction}
We study two examples of singular diffusion equations. One of them is 
an anisotropic total variation (TV) flow, the other one is the same equation
with the additive
isotropic linear 
diffusion,
\begin{equation} \label{eq_anisoTV} 
\frac{\partial u}{\partial t}=\beta\,\div 
\left(\frac{u_{x_1}}{|{u_{x_1}}|},\frac{u_{x_2}}{|{u_{x_2}}|}\right),
\end{equation}
\begin{equation} \label{eq_diffanisoTV} 
\frac{\partial u}{\partial t}= \gamma \Delta u +
\beta\,\div \left(\frac{u_{x_1}}{|{u_{x_1}}|},\frac{u_{x_1}}{|{u_{x_1}}|}\right).
\end{equation}
These problems are considered  
on a domain in $\bR^2$. Our goal is to study  features of solutions
like facets, i.e. flat parts of solutions with normals corresponding to the
singular directions. 
Our study was inspired by the phase transition
theory appearing in  crystal growth problems and
image restoration, where presence of walls and edges plays a significant
role, \cite{fukui}, \cite{ROF},
\cite{Spohn}.

Let us describe ideas behind this note.  The key element of the
systems we study is the
anisotropy. In both cases this determines the features of solutions. We will see
that numerical experiments appear to give 
almost the same despite fact that the second equation is not
degenerate.
The most spectacular phenomenon which is observed for this type of problems are
flat parts of solutions, connected with very strong diffusion, where $\nabla
u=0$. Such effects have been well studied for the isotropic
 total variation flow. We note that the interest in the TV flow arose from its
application to image analysis and reconstruction, \cite{osher}, \cite{beceen},
\cite{mazon-die}. Namely, any regular level sets of solutions to this flow
evolve by
the mean curvature. This property is used to smooth out contours and in
deblurring. We stress that
numerical algorithms exploit properties of this flow even implicitly. 

The case of anisotropic diffusion is not so well studied. Despite the
available papers like \cite{moll}, \cite{mazon-ann}, 
the mathematical theory is still far from the excellence. This changes however,
because of the interest in algorithms detecting or retaining special image
features
like edges and corners. A conspicuous example is the paper \cite{choksi}
on 2D bar codes. We observe a growing body of literature devoted to this
subject,
\cite{Preusser}, \cite{Birkholz}, \cite{gao}, \cite{chinczycy}, \cite{kuroda}.
We see the need to study
evolution equations which are likely  to preserve pronounced features of
solution
or its graphs besides facets, e.g. edges or corners. It turns out that
the equations
we study here may serve this purpose. The rigorous goals will be stated below,
but we also present numerical simulations in section \ref{ronu}, which
illustrate the qualitative features of solutions.

We set the following main goals of this paper:

\smallskip

\noindent
$\clubsuit$ \ to study facets, the flat parts of solutions, defined by $\nabla
u=0$;

\smallskip

\noindent
$\clubsuit$  \  to exhibit ruled surfaces, arising when one of the
components of the gradient of solution $u$ vanishes;

\smallskip

\noindent
$\clubsuit$  \ to construct special solutions, given by the explicit formulas
which
shows characteristic features of solutions;

\smallskip

\noindent
$\clubsuit$  \ to present numerical experiments and to show
evolution of interesting model shapes, which
were the motivation for looking for analytical results;

\smallskip

\noindent
$\clubsuit$ \ to propose a possible application to image
processing.

\smallskip

It is surprising that facets and ruled surfaces are the attributes of solutions
to both systems. The lack of degeneration 
in the second model results only in smoothing out effect appearing near
`corners' and some dispersion. Hence 
in practice, one could find the second equation as more suitable for
practical applications. Here, we present a series of solutions to both systems,
represented by  the gray scale. It shows some interesting differences, which nonetheless are very subtle.
  The initial data are represented by the last picture on Fig. \ref{fig1}.

\begin{figure}[h!]  
   \begin{center}
    \setlength{\fboxsep}{0pt}
    \setlength{\fboxrule}{0.5pt}
\fbox{\includegraphics[scale=0.28]
{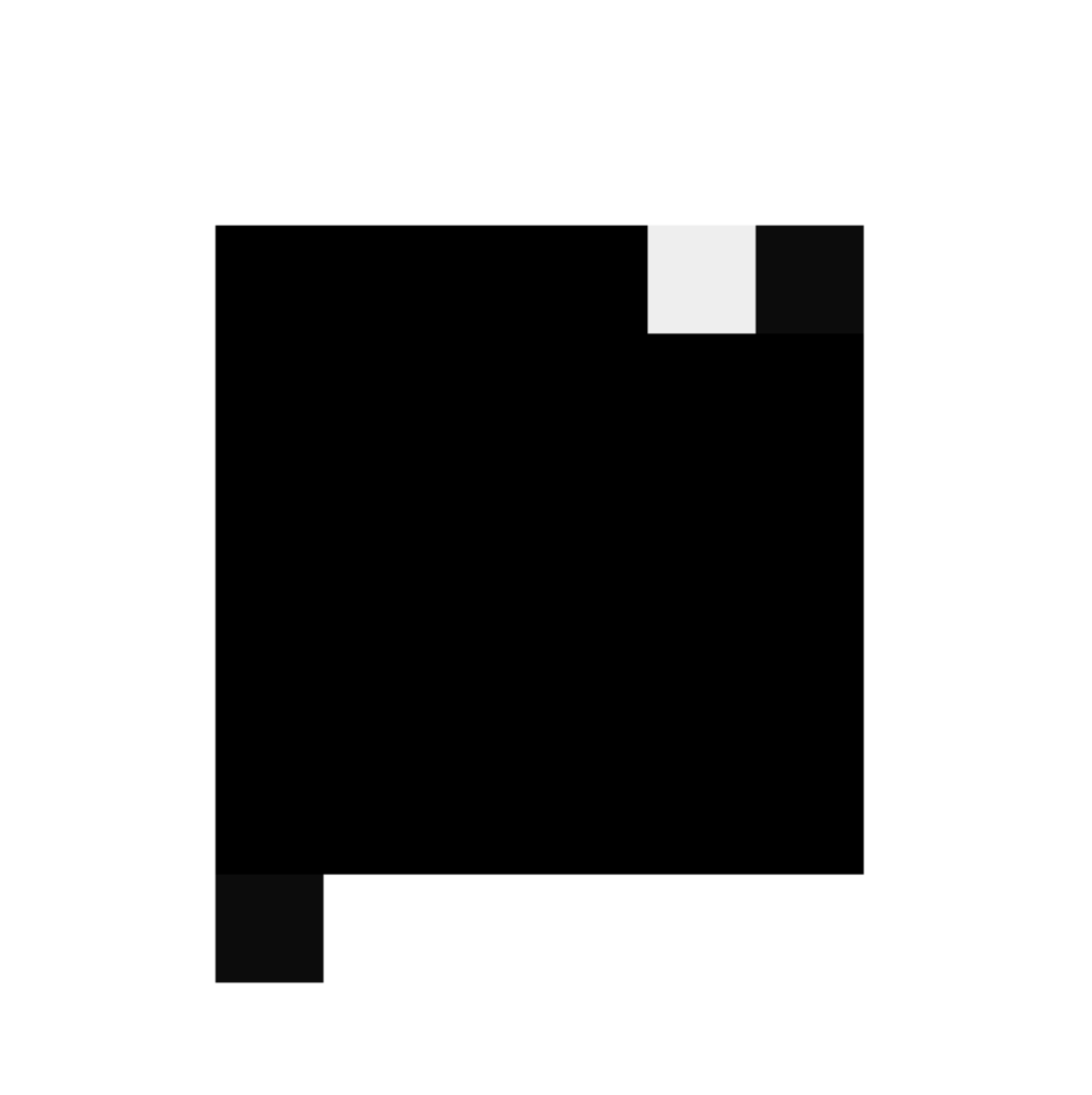}}\hspace{0.2cm}
\fbox{\includegraphics[scale=0.28]
{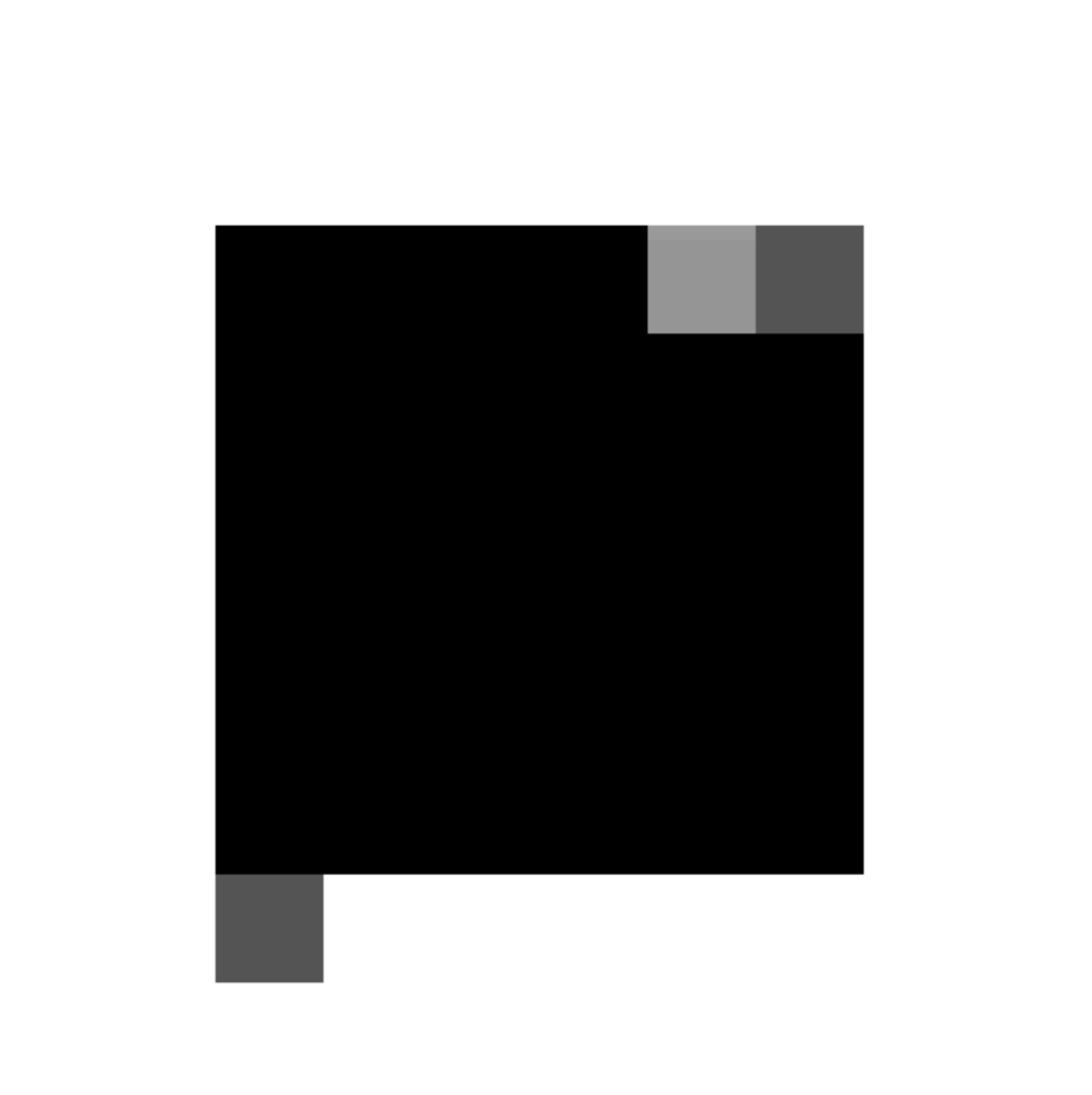}}\hspace{0.2cm}
\fbox{\includegraphics[scale=0.28]
{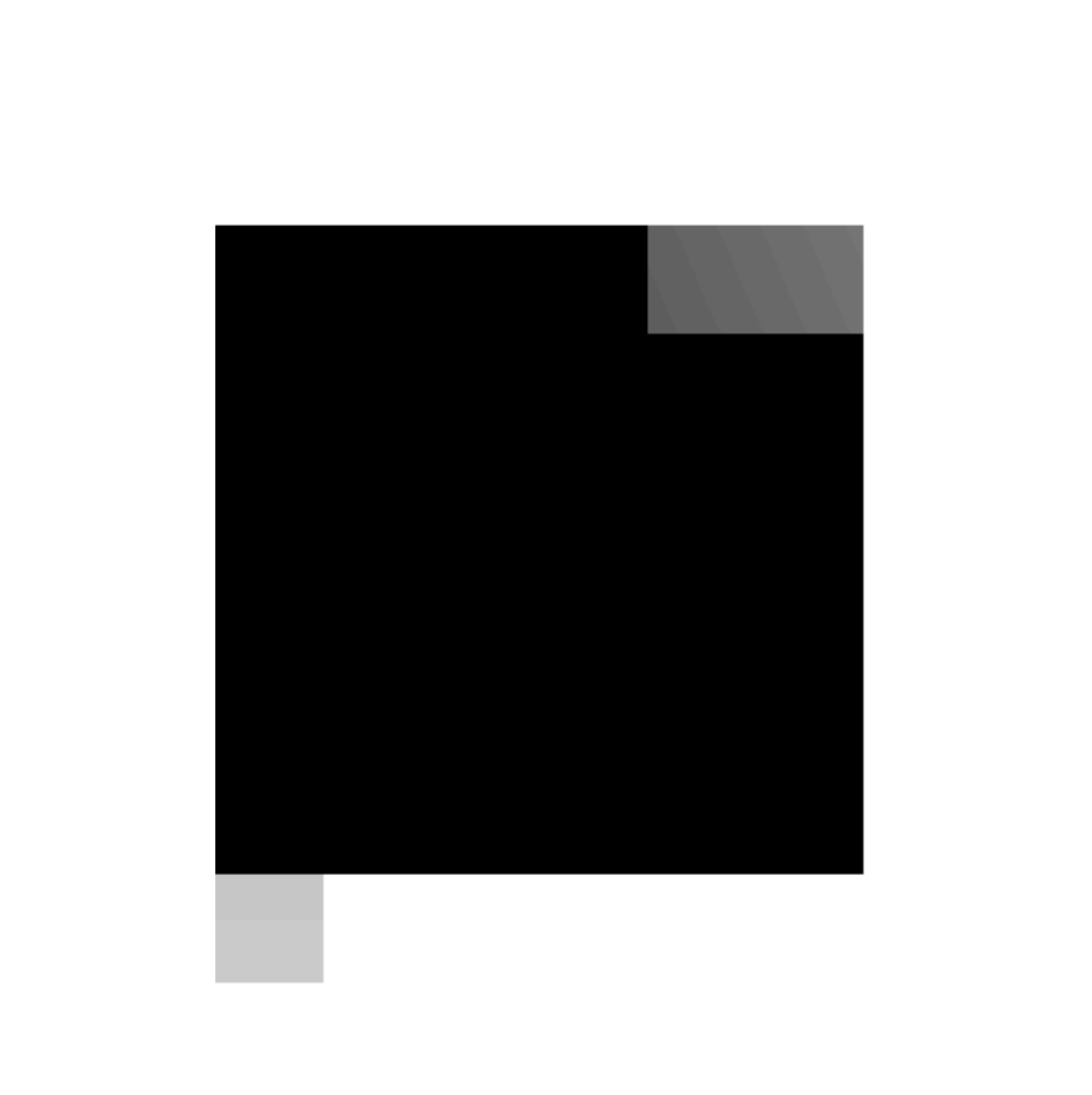}}\hspace{0.2cm}
\fbox{\includegraphics[scale=0.28]
{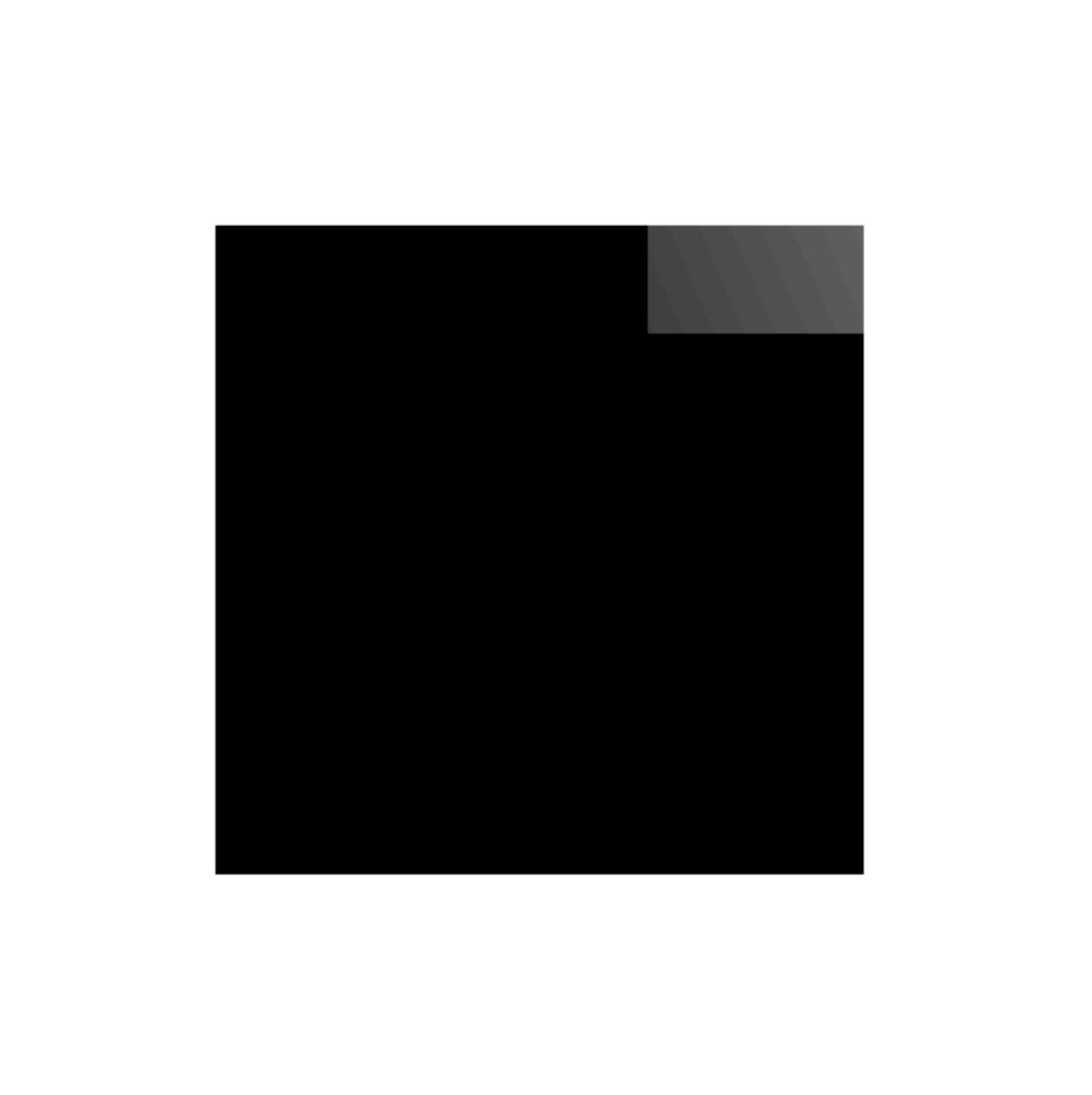}}\hspace{0.2cm}
\\[0.2cm]
\fbox{\includegraphics[scale=0.28]
{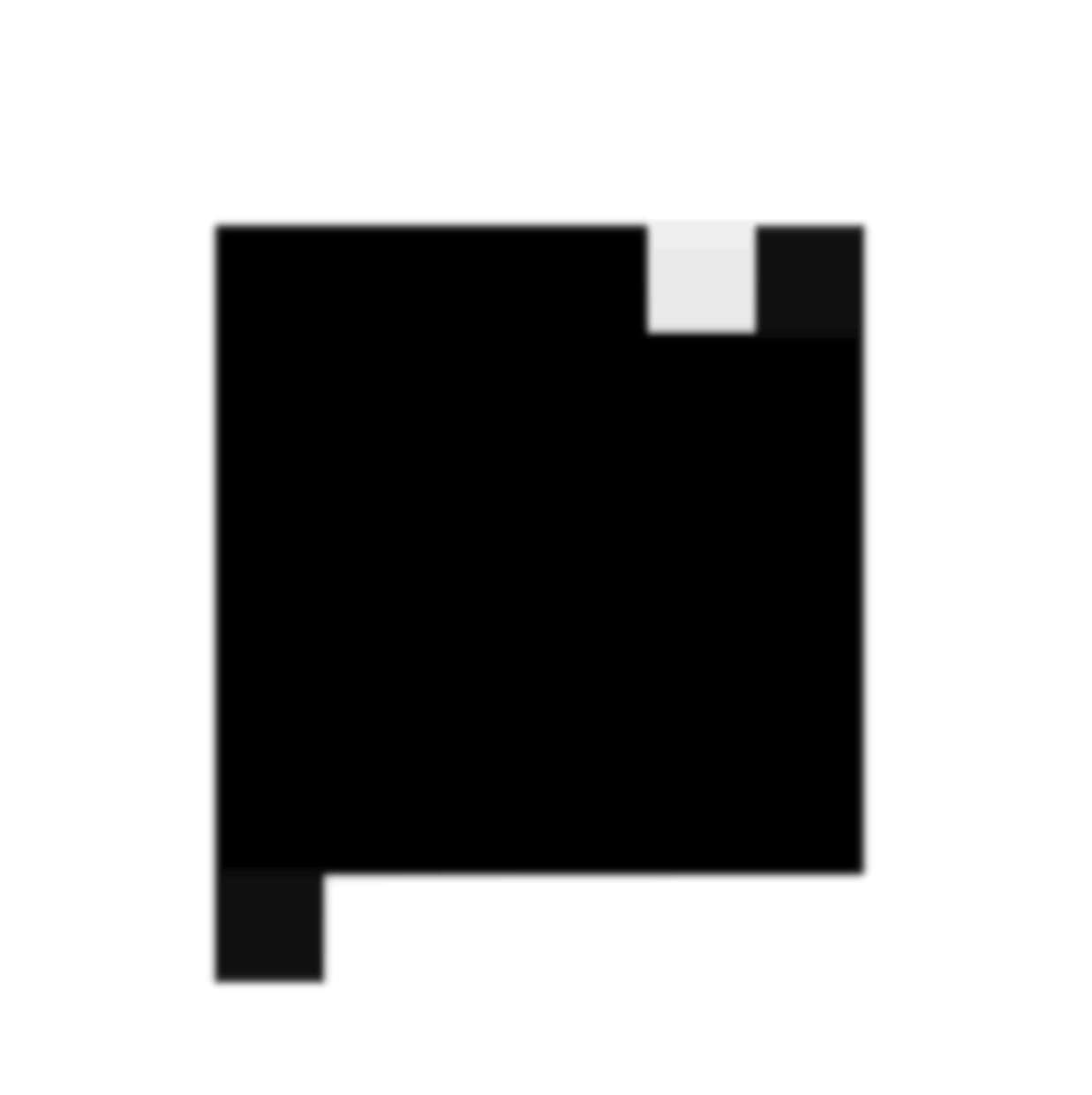}}\hspace{
0.2cm}
\fbox{\includegraphics[scale=0.28]
{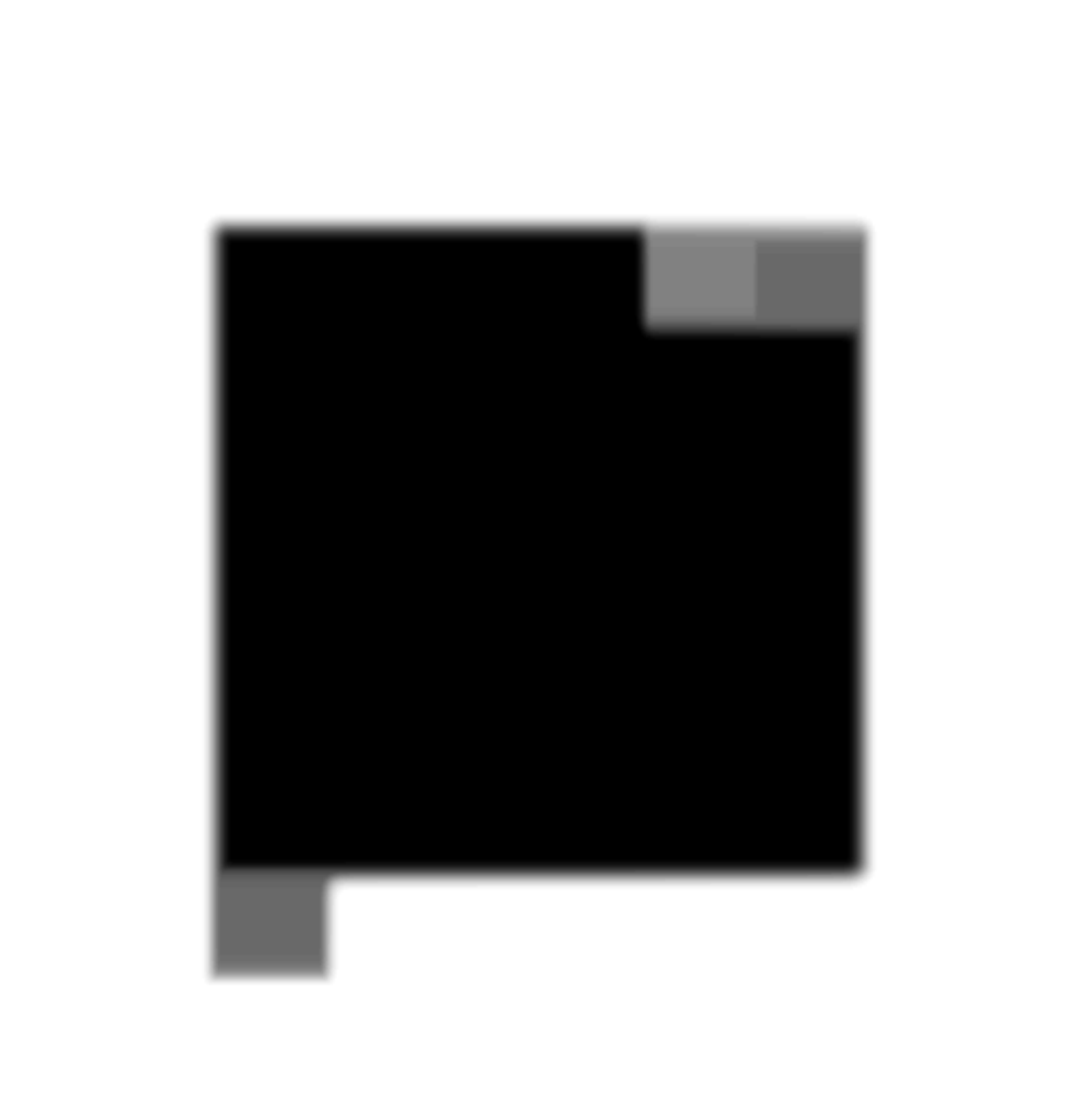}}\hspace{
0.2cm}
\fbox{\includegraphics[scale=0.28]
{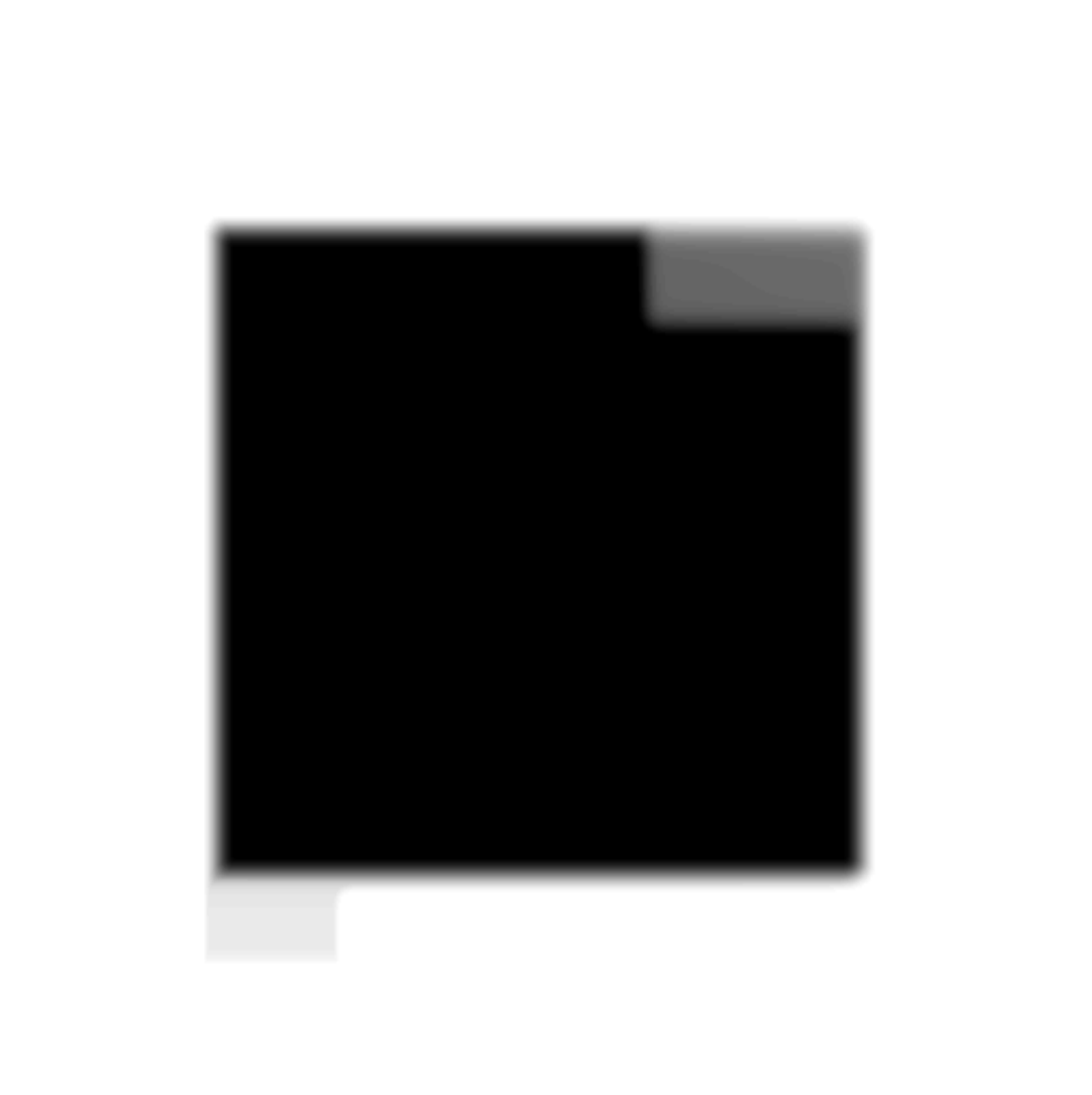}}\hspace{
0.2cm}
\fbox{\includegraphics[scale=0.28]
{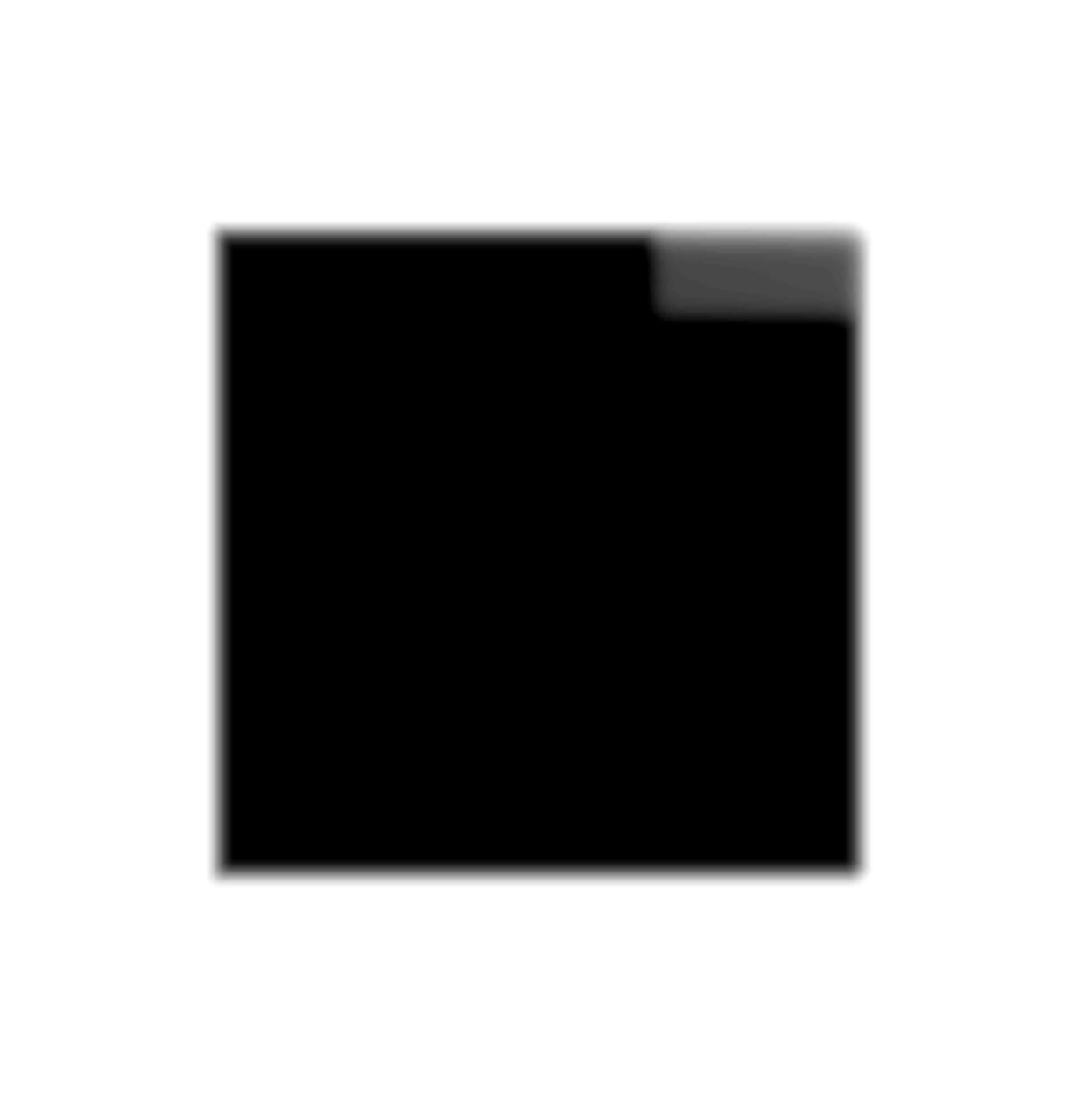}}\hspace{
0.2cm}\\[0.2cm]
\caption{Snapshots of evolution by (\ref{eq_anisoTV}) (first row) and
(\ref{eq_diffanisoTV}) (second row) of the same initial datum}\label{fig0}
\end{center}
\end{figure}

Applications to the phase transition theory are of particular interest, 
\cite{Spohn}, however, systems (\ref{eq_anisoTV})
and (\ref{eq_diffanisoTV}) are rather a simplification of more complex models. 
In the image processing, usefulness may be more straightforward to see.
The 
pictures presented in Section \ref{ronu} show a
possibility of reconstructing images.
The diffusion in the second system helps us to restore the picture. Positivity of $\gamma$ gives  averaging effects,
 but
strong anisotropic nonlinearity keeps edges in the chosen directions.

The upshot of these experiments is the following. The regularization of a very
singular system yields not only smooth solutions but also it preserves the main 
features of the original equation. We should keep in mind this important
observation in our future studies of systems with 
very singular nonlinear operators. 

At this point we note that the system with the added isotropic diffusion behaves
like 
phase field models with respect to free boundary problems including the mean
curvature flow. We mention just a 
few papers exploring the link, \cite{Danielle}, \cite{garcke}, \cite{goglione}.

Here, we do not plan to
present a consistent theory, but rather to 
pinpoint a few interesting results and phenomena to find motivation for
our future
deeper mathematical analysis. In fact, this note can be viewed as an attempt to
extend results for one-dimensional systems 
\cite{KMR,joss,new,M} on the two-dimensional case.  
To be more precise, we establish existence of solutions to both equations
(\ref{eq_anisoTV}) and (\ref{eq_diffanisoTV})
by using the theory of nonlinear semigroups. For this purpose we exploit the
gradient structure of (\ref{eq_anisoTV}) and (\ref{eq_diffanisoTV}).
Uniqueness is automatically guaranteed. This is  presented in the next section.
There we also present exact formulas for solutions. The advantage is that they
provide insight into the facet formation problem. Since the formulas do not
always fit the framework of semigroup solutions, we recall the notion of a weak
solution. 
The explicit solutions suggest 
that the flows of
(\ref{eq_anisoTV}) and (\ref{eq_diffanisoTV}) make ruled surfaces out of the
initial data, provided additional conditions are satisfied. This is rigorously
established in
Section \ref{secon}. This considerations require quite precise regularity
estimates established in Theorem in Section \ref{secreg}.

The paper is organized as follows. In Section 2 we state basic existence results for
systems (\ref{eq_anisoTV})
and (\ref{eq_diffanisoTV}) coming from the general theory. We point also to a
few interesting explicit solutions illustrating 
typical shapes. In Section  \ref{secreg} we show that the solutions are of
better
regularity, provided the initial data are smooth 
enough. Next we prove conditional results, which explain why flat regions and ruled surfaces are typical for graphs of
solutions. 
In Section \ref{ronu} we concentrate on numerical analysis and obtain a few
interesting numerical solutions. These results
show more direct phenomena which are able to be captured by the systems. In the appendix we present more complex example
of an explicit solutions to (\ref{eq_anisoTV}).

\section{Existence}\label{secexe}
We will use general tools exploiting the structure of the problem. 
In order to use the semigroup theory we notice that we  present
equations (\ref{eq_anisoTV}), (\ref{eq_diffanisoTV}) as gradient flows of
corresponding functionals on
$L^2(\Omega)$. We set
\begin{equation} \label{fi1}
\Phi_0(u)=\left\{
\begin{array}{ll}
\int_\Omega \beta(|u_{x_1}| + |u_{x_2}|) & \hbox{if } u\in
BV (\Omega),\\
+\infty & \hbox{if } u\in L^2(\Omega)\setminus BV (\Omega).
\end{array}\right.
\end{equation}
\begin{equation} \label{fi0}
\Phi_1(u)=\left\{
\begin{array}{ll}
\int_\Omega\frac\gamma2 |\nabla u|^2 + \beta(|u_{x_1}| + |u_{x_2}|) & \hbox{if }
u\in
H^1(\Omega),\\
+\infty & \hbox{if } u\in L^2(\Omega)\setminus H^1(\Omega);
\end{array}\right.
\end{equation}
Here, $\Omega$ is an open subset of $\bR^2$,
possibly unbounded, e.g. $\Omega=\bR^2$. We study the above equations on a
square with homogeneous Neumann boundary conditions
which is convenient from the numerical point of view, i.e. it is easier to
implement a numerical scheme on rectangular domain. We also consider periodic
boundary conditions. In general  $\beta,\gamma>0$, however we may scale the time
and from now on we put $\beta =1$ and admit $\gamma >0$ or $\gamma =0$ to have
a possibility to study both cases simultaneously.

It is obvious that $\Phi_1$ is well-defined and finite iff $u\in H^1$.
The correctness of the definition of $\Phi_0(u)$ is less obvious. In fact, this
is an
example of a more general situation studied in \cite[Section 2]{moll}.
Formula (\ref{fi1}) should be understood as follows,
\begin{equation}
 \int_\Omega  (|u_{x_1}| + |u_{x_2}|) := \sup\{\int_\Omega (z,Du)\,dx:\
z\in C^1_0(\Omega;\bR^2), \ |z|_\infty \le 1\},
\end{equation}
where $|(p_1,p_2)|_\infty := \max\{ |p_1|, |p_2|\}$.
It is now easy to check that these two functionals are convex, proper and lower
semicontinuous on $L^2(\Omega)$.

We notice that formally, the elliptic operator 
$\gamma\De u +
 \hbox{div}\,\left(\frac{u_x}{|{u_x}|},\frac{u_y}{|{u_y}|}\right) $ is the
first variation of functional $\Phi_1$ while
$ \hbox{div}\,\left(\frac{u_x}{|{u_x}|},\frac{u_y}{|{u_y}|}\right) $ is the
first variation of functional $\Phi_0$. Thus, equation (\ref{eq_diffanisoTV}) is
the gradient flow of 
$\Phi_1$ and equation (\ref{eq_anisoTV}) is the gradient flow of $\Phi_0$.
Keeping this in mind, we infer the following statement, where $A_i(u)=-\partial \Phi_i(u)$, $i=0,1$.
\begin{theorem}\label{tw1}
 Let us suppose that $u_0\in D(A_i)$, $i=0,1$. 
Then there
exists a unique function $u:[0,\infty)\to L^2(\Omega)$ such that:\\
(1) for all $t>0$ we have $u(t)\in D(A)$; \\
(2) $\frac{du}{dt} \in L^\infty(0,\infty, L^2(\Omega))$ and 
$\|\frac{du}{dt}\|_{L^\infty(0,\infty, L^2(\Omega))} \le \|A^o_i (u_0)\|_{L^2}$;\\
(3) $\frac{du}{dt} \in A_i(u(t))$ a.e. on $(0,\infty)$;\\
(4) $u(0) = u_0$.

In addition, $u$ has a right derivative at all $t\in[0,\infty)$ and
\begin{equation}\label{rnsem}
 \frac{d^+u}{dt} + A^o_i(u(t)) =0,
\end{equation}
where $A^o_i(u(t))$ is the minimal section of $A_i(u(t))$, (see \cite{brezis}).
\end{theorem}

Actually, since $A_i$ are subdifferential of convex functional, we  say more.
\begin{theorem}\label{tw2}
 Let us suppose that $u_0\in L^2$ and $A_i$ are as in Theorem \ref{tw1}, then there exists a unique solution to equation 
$$
\frac{d^+u}{dt} + A_i (u(t)) \ni 0,\quad u(\cdot,0) = u_0.
$$
Moreover, for all $t>0$ 
$u(t)$ belongs to $D(A)$ and (\ref{rnsem}) holds.
\end{theorem}

We notice that Theorem \ref{tw1} follows from \cite[Theorem 3.1]{brezis}, while
\cite[Theorem 3.2]{brezis} implies  our Theorem \ref{tw2}.
In  Theorem \ref{tw1} we refer to the domains  $D(\partial\Phi_0)$ and $D(\partial\Phi_1)$, however, we 
abstain from exact description of theses sets.
The semigroups obtained by these theorems are contraction semigroups, thus if
$u^n_0\to u_0$ in $L^2(\Omega)$, then for all fixed $t$ we have $u^n(t)\to
u(t)$. This observation will be used in the constructions of examples of
solutions based on explicit calculations.

\smallskip 
This general result gives us the justification for our exact formulas for
solutions.
They are particularly valuable when we strive to study motion of facets or
other special properties.
First, for special data we cook up explicit formula for a solution to
(\ref{eq_anisoTV}).
To keep the simplest setting we consider the equations in the whole plane.
We  construct $u$ (see formula (\ref{wzor2})) a
 solution  to a differential inclusion 
\begin{equation}\label{ink}
 u_t + A_0(u) \ni 0\qquad \hbox{in }\bR^2 \times \bR_+
\end{equation}
in place of (\ref{eq_anisoTV}), with the initial datum
\begin{equation}\label{r-dane-pp}
 u_0(x_1,x_2) = x_1^2+x_2^2 -2R^2.
\end{equation}
This initial condition does not belong to  $L^2(\bR^2)$, but 
$u_0 \in L^2_{loc}(\bR^2),$ $u_0, \nabla u_0\in BV_{loc}(\bR^2).$
The same property will be valid for $u(\cdot, t)$.
%
%
Hence the notion of solution introduced in Theorem \ref{tw1} by (\ref{rnsem}) is not quite appropriate. This is
 why we introduce in (\ref{rnsl}) the notion of a weak solution.

\begin{proposition}\label{pr0}
Formula (\ref{wzor2}) below yields a weak solution  to
(\ref{eq_anisoTV}) in $\bR^2$ with data (\ref{r-dane-pp}), understood as 
(\ref{rnsl}). 
Moreover, (\ref{ink}) is satisfied in
$\bR^2\times \bR_+$ in a
pointwise manner with the exception of a one dimensional set and the solution
is Lipschitz continuous, but not $C^1$. 
\end{proposition}
{\it Proof.} Let us define
\begin{equation*}
 \xi^+(t) \equiv \xi(t) = \left(\frac{3}{2}\right)^{\frac{1}{3}} t^{\frac{1}{3}}
\equiv -\xi^-(t)
\mbox{ \ and \ } 
 h(t) 
= \left(\frac{3}{2}\right)^{\frac{2}{3}} t^{\frac{2}{3}} +
x_1^2+ x_2^2 -2R^2.
\end{equation*}
We notice that for $t\ge0$ the quantities $\xi^\pm(t)$ are uniquely defined by the condition
$$
u_0(\xi^\pm(t), 0) = h(t) = u_0( 0,\xi^\pm(t)) \mbox{ \ \ with \ } 
\xi^\pm(t) = \pm\sqrt{h(t)}.
$$
The final observation is that these functions satisfy the equation
$$
\frac{dh}{dt} = \frac{2}{\xi^+(h) - \xi^-(h)},\quad h(0) = 0.
$$
Now, we  write the advertised formula for solutions to (\ref{eq_anisoTV}),
\begin{equation}\label{wzor2}
 u(x,t) =\left\{
\begin{array}{ll}
2h(t) 
& |x_1|, |x_2| \le \xi(t),\\
h(t) + x_2^2 -2R^2 & |x_1| \le \xi(t), |x_2| > \xi(t),\\
h(t) + x_1^2 -2R^2  & |x_2| \le \xi(t), |x_1| > \xi(t),\\
x_1^2 + x_2^2 - 2R^2 & |x_1|, |x_2| > \xi(t).
\end{array}
\right.
\end{equation}
This formula defines a Lipschitz continuous, but not a $C^1$ function.

We shall calculate $u_t$. We obviously obtain
\begin{equation*}
 u_t(x,t) =\left\{
\begin{array}{ll}
2 h'(t) & |x_1|, |x_2| < \xi(t),\\
h' (t)& |x_1| < \xi(t), |x_2| >\xi(t),\\
h'(t)& |x_2| < \xi(t), |x_1| >\xi(t),\\
0 & |x_1|, |x_2| > \xi(t).
\end{array}
\right.
\end{equation*}

The point is to calculate a
selection of $\cL(\nabla u):= (u_{x_1}/|u_{x_1}|, u_{x_2}/|u_{x_2}|) $, where at least one of the components of
$\nabla u$ vanishes.
For this purpose, we take advantage of the special structure of this operator,
permitting us
to use what we learned about the one dimensional case, see \cite{joss},
\cite{new}. This yields
\begin{equation}\label{wzor2-1}
 \cL(\nabla u)(x,t) =\left\{
\begin{array}{ll}
\frac{1}{\xi(t)}(x_1, x_2) & |x_1|, |x_2| \le \xi(t),\\
(\frac{x_1}{\xi(t)}, \sgn x_2)   & |x_1| \le \xi(t), |x_2| > \xi(t),\\
(\sgn x_1, \frac{x_2}{\xi(t)})& |x_2| \le \xi(t), |x_1| >
\xi(t),\\
(\sgn x_1, \sgn x_2) & |x_1|, |x_2| > \xi(t).
\end{array}
\right.
\end{equation}
This is a Lipschitz continuous vector field.  Let us check if $u$ is a weak
solution to (\ref{eq_anisoTV}). We
recall that $u$ is a weak solution iff
\begin{equation}\label{rnsl}
 (u_t,\phi)+(\sigma,\nabla \phi)=0 \mbox{ \ \ in } \mathcal{D}'([0,T))
\end{equation}
for each $\phi \in C^\infty_c(\R^2 \times [0,T))$ and $\sigma^i \in \sgn u_{x_i}$, $i=1,2$.

Here we put $\sigma=\cL(\nabla u)(x,t)$, where $\cL(\nabla u)(x,t)$ is given
by (\ref{wzor2-1}). Since $\sigma$ is Lipschitz continuous, we are allowed to
integrate by 
parts in the second term of the LHS in (\ref{rnsl}), getting $u_t =
\di\cL(\nabla u)$.
If we take into account the
explicit form of $h(t)$, it is easy to see that
the identity holds everywhere, except a two-dimensional subset 
$\{(x_1,x_2,t): \ |x_1|=\xi(t)\hbox{ or }|x_2|=\xi(t)\}$ of
$\bR^2 \times \bR_+$. \qed

\smallskip
This example was relatively easy to present, because the problem was consider on
the whole $\bR^2$. It is also interesting to see if a similar formula works  on
a
bounded domain with a boundary condition. In Proposition \ref{pr2} in the
Appendix, we present a similar, but more messy formula for a square with Neumann
boundary data.

\smallskip  
The same notion of  weak solutions like introduced in (\ref{rnsl}) may be used
also when Theorems \ref{tw1}
 and \ref{tw2} are applicable. However, it is easy to see that if $u$ satisfies
(\ref{rnsem}), then
 it is a weak solution in  the sense of  (\ref{rnsl}). 
In addition, if $u^1$ and  $u^2$ are weak solutions with the same initial data, then they must coincide. 

%
%

\smallskip 

Next, we study solutions to (\ref{eq_anisoTV}) with data just in $BV$ space.

\begin{proposition}\label{pr3}
 Let us suppose that $\Omega = (-L,L)^2$,  and  $\alpha\in(0,L)$, $M>0$. We set
$$
u_0(x_1,x_2) = - M \chi_{(-\alpha, \alpha)^2}.
$$
Then a unique solution to (\ref{eq_anisoTV}) with the above initial data is
given by formula (\ref{wzor3}) below.
\end{proposition}

{\it Proof.} Let us set 
\begin{equation}\label{wzor3}
 u(x,t) =\left\{
\begin{array}{ll}
 \frac{2t}{\alpha} - M  & |x_1|, | x_2| \le \alpha,\\
-\frac{2\alpha}{L^2-\alpha^2} t & \hbox{otherwise }
\end{array}\right.
\end{equation}
Checking correctness requires defining $\cL(\nabla u)$ in a proper way. We
define two auxiliary functions
\begin{equation*}
Z_1(x) = \left\{
\begin{array}{ll}
 \frac{-L-x}{L-\alpha} & x \in (-L,-\alpha),\\[3pt]
\frac{1}{\alpha}x & |x|\le \alpha,\\[3pt]
\frac{L-x}{L-\alpha} & x \in (\alpha,L),
\end{array} \right.
Z_2(x) = \left\{
\begin{array}{ll}
- \alpha \frac{x+L}{L^2-\alpha^2} & x \in (-L,-\alpha),\\[3pt]
-\frac{1}{L+\alpha}x & |x|\le \alpha,\\[3pt]
 -\alpha\frac{x-L}{L^2-\alpha^2}  & x \in (\alpha,L),
\end{array} \right.
\end{equation*}
We now define $\cL(\nabla u)$  by setting
$$
\cL(\nabla u) = (Z_1(x_1)\chi_{\{|x_2| \le \alpha\}} + Z_2(x_1)\chi_{\{|x_2| >\alpha\}},
Z_1(x_2)\chi_{\{|x_1| \le \alpha\}} + Z_2(x_2)\chi_{\{|x_1| >\alpha\}}).
$$
It is now easy to check that
$$
u_t = \di \cL (\nabla u).
$$
We use the same argumentation as in the proof of Proposition \ref{pr0},
however the difference is that here the defined above $\cL(\nabla u)$
is not so regular. In order to take the divergence we are required to
control only appropriate directional 
derivatives, so the  form of $\cL(\nabla u)$ and Lipschitz continuity of $Z_1$ and $Z_2$ allow us to obtain 
the desired identity. This equality holds pointwise in $\bR^2\times \bR_+$
except for a two-dimensional set.
This formula is valid until the time when two facets merge into a constant
stationary state at the extinction time $t=T_{ext}$,
$$
T_{ext} = M\left(\frac{2}{\alpha} + \frac{2\alpha}{L^2-\alpha^2}\right)^{-1}.
\eqno\Box
$$ 

\smallskip

Finally we point one special solution to the second system.
We  show existence of a moving front for (\ref{eq_diffanisoTV}),
but without any boundary conditions.

\begin{proposition}\label{pr4}
 Let us fix $\alpha>0$, then each of the functions given by the formula below is
traveling front solution to (\ref{eq_diffanisoTV}),
\begin{equation*}
  u^\alpha(x,t) =\left\{
\begin{array}{ll}
\frac{2t}{\alpha}  & |x_1|, |x_2| \le \alpha,\\[3pt]
\frac{2t}{\alpha}  +\frac{1}{\alpha} x_2^2  & |x_1| \le \alpha, |x_2| >
\alpha,\\[3pt]
\frac{2t}{\alpha}  +\frac{1}{\alpha}  x_1^2  & |x_2| \le \alpha,
|x_1| > \alpha,\\[3pt]
\frac{2t}{\alpha}  +\frac{1}{\alpha}  (x_1^2 + x_2^2) & |x_1|, |x_2| >
\alpha.
\end{array}
\right.
\end{equation*}
\end{proposition}
Checking the correctness of this formula is easier than in the previous case.
The above formula makes it clear that no  traveling front solution is possible for
(\ref{eq_anisoTV}). In the Appendix we point an extra explicit solution to (\ref{eq_anisoTV}).

\section{Extra regularity}\label{secreg}

In this part we show that solutions to (\ref{eq_anisoTV}) and (\ref{eq_diffanisoTV}) 
 obtained via Theorems \ref{tw1} and \ref{tw2} are of a better regularity. It
will be very important for deducing 
some qualitative features of solutions.

\begin{theorem}\label{tmxx}
 Let $u_0 \in H^1(K)$, then the solution to  (\ref{eq_diffanisoTV}) given by Theorem \ref{tw1} fulfills the
following estimate
\begin{equation}\label{e1}
 \|u_t\|_{L_2( 0,T\times K)} + \sup_{t\in [\delta,T]} \|u_t,\gamma \nabla^2
u\|_{L_2(K)}(t) \leq DATA(\delta).
\end{equation}
\end{theorem}

\smallskip 

{\bf Proof.} We consider both cases at ones: $\gamma =0$ and $\gamma >0$.
After mollifying the system we test it by $u_t$ getting
\begin{multline}\label{e2}
\qquad \int_0^T \int_K u_t^2 dxdt + \sup_{t \in [0,T]} \int_K [\frac\gamma2 |\nabla u|^2 + |u_{x_1}| + |u_{x_2}|] dx\\
\leq 2 \int_K [\frac\gamma2|\nabla u_0|^2 +|u_{0,x_1}|+|u_{0,x_2}|] dx.
\end{multline}
The structure of the equation allows us to differentiate the system with respect
to $t$.
\begin{equation}\label{e3}
 u_{tt} - (\d_{x_1} (\sgn u_{x_1})_t + \d_{x_2} (\sgn u_{x_2})_t + \gamma \Delta u_t)=0.
\end{equation}
Let $\eta$ be a time dependent function such that $\eta (0)=0$ and  for $\delta >0$ $\eta \equiv 1$, 
then we test (\ref{e3}) by $u_t \eta$ getting
\begin{multline*}
 \sup_{t\in [0,T]} \int_K u_t^2 \eta dx + \int_0^T \int_K \eta
 [\delta(u_{x_1}) u_{x_1t}^2 +\delta(u_{x_2}) u_{x_2t}^2 +\gamma |\nabla u_t|^2]dxdt\\
\leq 2 \int_0^T \int_K u_t^2 \eta' dxdt. \qquad
\end{multline*}
But the r.h.s. is bounded by (\ref{e2}), so we have 
$u_t \in B(\delta,T;L_2(K)).$
Taking into account the above information, we consider (\ref{eq_diffanisoTV}) in
the following modification
\begin{equation}\label{e5}
 -[ \d_{x_1} (\delta(u_{x_1}) u_{x_1x_1}) + \d_{x_2} (\delta(u_{x_2}) u_{x_1x_2}) + \gamma \Delta u_{x_1}]=-u_{tx_1}
\end{equation}
here time is a fixed parameter. Testing (\ref{e5}) by $u_{x_1}$, we get 
\begin{equation*}
 \int_{K} ( \delta(u_{x_1}) u_{x_1x_1}^2 + \delta(x_2) u_{x_1x_2}^2  + \gamma |\nabla u_{x_1}|^2 ) dx \leq \int_K |u_t u_{x_1x_1}|dx
\end{equation*}
which gives the estimates on $\gamma \int_K|\nabla u_{x_1}|^2 dx$.
The same we have for $x_2$. The estimate (\ref{e1}) is proved. \qed

\smallskip 

If we use $t^2$ as a test function $\eta$ above, then we obtain
information on the blow up of $\|u_t\|$. Namely, it is easy to see that
\begin{corollary}
 Under the assumptions of Theorem \ref{tmxx}, we have
$\|u_t\|_{L^2} \le C  t^{-1/2}.$
\end{corollary}

We shall emphasize that the terms $\delta(u_{x_1}) u_{x_1x_1}^2$ are considered
just formal, to be precise
we shall treat them as limits coming from analysis done on the level of approximation. 

\smallskip 

%

\section{Ruled surface and convexity}\label{secon}

The first phenomenon, which is very expected for this type of systems, are features of 
minimizers and maximizers of the solution. We ask about a possible structure of
sets where 
the function $u$, for fixed time $t$, admits extrema. Since the issue of regularity is not well
studied, here, we prove only the following result.

\begin{proposition}
 Let $u$ be a solution to system (\ref{eq_anisoTV}) or (\ref{eq_diffanisoTV}). Let $t>0$ and 
for $x_0$ in the domain $u(\cdot,t)$ has a minimum at $x_0$ and in addition, $u(\cdot,t)$
is a convex function different from a constant in a neighborhood $N$ of set
$u(\cdot,t)=u(x_0,t)$, then the set 
\begin{equation*}
M= \{x: u(x,t)=u(x_0,t)\}\cap N
\end{equation*}
is a closed set with nonempty interior.
\end{proposition}

{\bf Proof.} 
We deduce that there is a sequence $m_n$ converging to $m:= u|_{M}$ from above
and such that each
level set $\{u(\cdot,t)=m_n\}$ is a convex closed curve.
Moreover, the sets $M_n=\{u(\cdot,t) \le m_n\}$ are
convex too. We integrate equation (\ref{eq_diffanisoTV}) over this set
\begin{equation*}
 \int_{M_n} u_t- \gamma \Delta u - \div (\sgn u_{x_1},\sgn u_{x_2})
\,dx_1dx_2 =0.
\end{equation*}
Integration by parts
leads us to the following conclusion,
\begin{equation*}
 \int_{\{u=m_n\}} (\gamma \frac{\d u}{\d n} + n_1 \sgn u_{x_1} + n_2 \sgn
n_{x_2} ) d\cH^1 =  \int_{M_n} u_t \,dx_1dx_2.
\end{equation*}
But convexity implies that
$\frac{\d u}{\d n} \geq 0$ at $\d M_n$.
At the same time for almost all $y$ functions
$x_1 \mapsto u(x_1,y,t)\quad\hbox{and}\quad x_2 \mapsto u(y,x_2,t)$
are monotone, hence in a neighborhood of $M_n$
\begin{equation*}
 n_1\sgn u_{x_1} + n_2 \sgn u_{x_2} = |n_1| + |n_2| \geq |n|=1.
\end{equation*}
We conclude that we obtain 
\begin{equation*}
 \int_{ \{u=m_n\} } d\cH^1\leq |M_n|^{1/2} (\int_{M_n} u_t^2 dx_1dx_2)^{1/2}.
\end{equation*}
Moreover, since $u$ is not constant, then the sets $M_n$ must have a positive
two-dimensional Lebesgue measure. However, due to the isoperimetric inequality
we have
$$
\cH^1(\partial M_n) \ge \frac{1}{2\sqrt{\pi}}|M_n|^{1/2},
$$
the identity holds for the disc. Hence 
\begin{equation}\label{r46}
 C\leq (\int_{M_n} u_t^2 dx_1dx_2)^{1/2}.
\end{equation}
However, due to Theorem \ref{tmxx}, $u_t$ is square integrable, so the RHS of
(\ref{r46}) above cannot go to zero when $n\to \infty$. Thus,    $M$ is a
convex set of positive two-dimensional measure, hence it must have nonempty
interior.
\qed

%
%

\smallskip 

The next feature concerns the shape of the graph of solutions.
The example presented in the earlier section suggests that the graph of the
solution develops parts which are ruled surfaces. To be more precise, we will
show that if the level sets of a convex solution $u(\cdot, t)$ at $t>0$ are 
regular, then the graph contains
ruled surfaces which are of positive two-dimensional measure. The tangent is
orthogonal to  vector $(0,1)$ or $(1,0)$. The
precise phenomenon is prescribed by the lemma below.

\smallskip 

\begin{lemma}\label{lemxx}
 Let $u$  be a sufficiently regular solution to (\ref{eq_anisoTV}) or
(\ref{eq_diffanisoTV}), (in other words $\gamma$ 
is equal to $0$ or $1$). That means, for a fixed $t$ the restriction of
$u(\cdot,t)$ to an open set $U$ is convex. 
Furthermore, we assume that for given $c\in \R$, the  level  set
\begin{equation*}
 S(c)=\{x\in K: u(t,x)=c\}
\end{equation*}
is regular, i.e. $\nabla u|_{S(c)}$ exists $\cH^1$--a.e. on $S(c)$ and $\nabla
u|_{S(c)}\neq 0$ $\cH^1$ a.e. Then sets 
\begin{equation*}
 M^+_1=\{ x: x=(m^+_1,x_2) \in S(c)\}, \mbox{ where } m^+_1=\max\{x_1: (x_1,x_2)
\in S(c)\};
\end{equation*}
\begin{equation*}
 M^-_1=\{ x: x=(m^-_1,x_2) \in S(c)\}, \mbox{ where } m^-_1=\min\{x_1: (x_1,x_2)
\in S(c)\};
\end{equation*}
\begin{equation*}
 M^+_2=\{ x: x=(x_1,m^+_2) \in S(c)\}, \mbox{ where } m^+_2=\max\{x_2: (x_1,x_2)
\in S(c)\};
\end{equation*}
\begin{equation*}
 M^-_2=\{ x: x=(x_1,m^-_2) \in S(c)\}, \mbox{ where } m^-_2=\min\{x_2: (x_1,x_2)
\in S(c)\}
\end{equation*}
do not contain isolated points.
\end{lemma}

{\bf Proof.} It suffices to consider just one of these sets, e.g. $M^+_1$. Let
us
suppose that our claim fails and $M^+_1 = \{ p \}$, in
other words, function
$x_2\mapsto u(m^+, x_2)$ has a strict minimum at $x_2=p$ in the interval
$[p-\ell, p+\ell]$. Due to the continuity of $u$ we notice that if
$u(x_1,\cdot)$
restricted to $[p-\ell, p+\ell]$ attains its minimum on $[p^-(x_1), p^+(x_1)]$,
then $p^-(x_1)$, $p^+(x_1)$ converge to $p$ as $x_1$ goes to $m^+$. In
particular, $u(x_1,  p\pm\ell)> u(m^+,p)$ for $x_1$ close to $m^+$. The last
observation combined with monotonicity of $u^-_{x_2}(x_1,\cdot)$,
$u^+_{x_2}(x_1,\cdot)$ 
implies that
$$
u^\pm_{x_2}(x_1, p+\ell)>0,\qquad u^\pm_{x_2}(x_1, p-\ell)<0,
$$
for all $x_1$ close to $m^+$. Thus, we can consistently define
\begin{equation*}
 \sgn u_{x_2} =
\left\{
\begin{array}{ll}
 1&\hbox{ on }\{(x_1, p+\ell): \ x_1\in (m^+ - \delta, m^+
+\delta\},\\
-1&\hbox{ on }\{(x_1, p-\ell): \ x_1\in (m^+ - \delta, m^+
+\delta\}
\end{array}
\right.
\end{equation*}
Let us take rectangles, 
$R_k =[m^+,m^+ + \delta_k]\times [p-\ell,p+\ell]$, where $\delta_k \le \delta$
and $\delta$ is so small that the above considerations are valid. We integrate 
$(\sgn u_{x_2})_{x_2}$ over $R_k$. 
We obtain,
\begin{equation*}
 \int_{R_k} (\sgn u_{x_2})_{x_2} dx_1dx_2 = \int_{\partial R_k} \sgn u_{x_2}
n_2 
= 2\cdot 2\ell.
\end{equation*}
  We may assume that function $x_1\mapsto u(x_1,p)$ is increasing on
$[m^+,m^+ +\delta]$, otherwise we could consider $u(-x_1,x_2)$, in place of
$u(x_1,x_2)$. 

Since $x_2\mapsto u(m,x_2)$ is convex, with minimum at $x_2=p$, then it must
be increasing on $[p, p+\delta]$ and due to our assumption $u(m^+,x_2)>u(m^+,p)$
for $x_2\neq p$.
Moreover, all lines $l_a=\{ (x,a): x\in \bR\}$ intersect $S(c)$ for $a$ close to
$p$, i.e. $|p-a|<\delta$, otherwise $S(c)$ would be a point, i.e. a singular
level set.
Let us suppose 
\begin{equation}\label{pola}
 (\tilde x_1, \tilde x_2) \in l_a\cap S(c),
\end{equation}
with $x_2$ close to $p$. Then,
$$
u^\pm_{x_1} (m,\tilde x_2)>0.
$$
Equality above is excluded because it contradicts (\ref{pola}) and
the monotonicity
of $u_{x_1}(\cdot, \tilde x_2)$. By the monotonicity of the derivative of a
convex
function we also obtain 
$u^\pm_{x_1} (m,\tilde x_2) < u^\pm_{x_1} (m,\tilde x_2+\delta_k)$. Thus, we may
consistently define $\sgn u_{x_1} = 1$ on the sides of $R_k$ parallel to the
vertical axis. 

Let us now integrate our equation over $R_k$,
\begin{equation*}
 \int_{R_k} \d_{x_1}(\sgn u_{x_1}) + \d_{x_2}(\sgn u_{x_2}) dx_1 dx_2 =
\int_{R_k} (u_t-\gamma\Delta u) dx_1 dx_2.
\end{equation*}
performing integration by parts on the LHS and taking into account observations
collected above, we conclude that
\begin{equation*}
 \int_{\d R_k} \sgn u_{x_2} n_2 d\sigma = \int_{R_k} (u_t -\gamma \Delta
u)\,dx_1 dx_2.
\end{equation*}
We continue the calculations. Using the square integrability
of
$u_t -\gamma \Delta u$ established in Theorem \ref{tmxx} we obtain that
\begin{equation*}
 4 \ell \leq |\int_{R_k} (u_t -\gamma \Delta u)\, dx_1 dx_2|\leq  
(2\ell\delta_k)^{1/2} 
\left( \int_{R_k} | u_t -\gamma \Delta u|^2\,dx_1 \right)^{1/2}
\end{equation*}
i.e.
\begin{equation*}
 4 \ell^{1/2}  
\leq  (2\delta_k)^{1/2} 
\left( \int_{R_k} | u_t -\gamma \Delta u|^2\,dx_1 \right)^{1/2}.
\end{equation*}
If  $\delta_k$ goes to zero, then we reach a contradiction. Thus, $M_1^+$
may not be a point. \qed

\begin{theorem}
 Assume that for $t>0$ and a region $A$ the solution $u(\cdot, t)$,  restricted
to $A$, is convex and the level sets of $u(\cdot, t)$ satisfy
the regularity assumption of Lemma \ref{lemxx}, then
sets 
\begin{equation}
 S_1=\{(x,u(x,t)): x \in A, u_{x_1}(x)=0\} \mbox{ and }  S_2=\{(x,u(x,t)): x \in A, u_{x_2}(x)=0\}
\end{equation}
are ruled surfaces, provided that $u_t,\gamma \nabla^2 u$ is bounded pointwisely, and
$\nabla u$ is continuous for  $\gamma=0$.
\end{theorem}

The proof of the above lemma follows immediately from Lemma \ref{lemxx}.

\section{Numerical experiments} \label{ronu}

The algorithm used to perform numerical experiments is based on the duality
approach considered by Chambolle \cite{Cham04}. He computed a minimizer of the
total variation model for the image denoising proposed by
Rudin~et~al.~\cite{ROF}. 
In order to adapt this approach to solve the equation (\ref{eq_diffanisoTV}), we
note  first that the semi-discretization of (\ref{eq_diffanisoTV}) yields the
following iterative scheme
\begin{equation}
\label{semidiscret} 
\frac{u^{m}- u^{m-1}}{\delta t}=\gamma \Delta u^m + \beta\,
\nabla\cdot\left(\frac{u^{m}_{x_1}}{|u^{m}_{x_1}|}, \frac{u^{m}_{x_2}}{|u^{m}_{x_2}|}\right)\,,
\end{equation}
where $u^m(x):=u(x,t_m)$ for $m=1,2,...$ and $x\in\mathbb{R}^2$, the initial data $u^0(x):= f(x)$ for $x\in\mathbb{R}^2$, where $f\in L^{\infty}(\Omega)$ is a given function, and $0<\delta t = t_m-t_{m-1}$ for $m=1,2,...$ denotes the time discretization step. For the convenience of notation, assume that $\delta t=1$ and consider the case $m=1$. 
Then, we note that the equation (\ref{semidiscret}) can be seen as the optimality condition for the minimization problem
\begin{equation}
\label{primal} 
\min_{u\in H^1(\Omega)}\left(\dfrac{1}{2}\int_{\Omega}(u-f)^2 +\gamma |\nabla u|^2\,dx + \beta\int_{\Omega} 
\abs{u_{x_1}} +  \abs{u_{x_2}}\,dx\right)\,.
\end{equation} 

Let us introduce the differential operator $A_{\gamma}:H^1(\Omega)\rightarrow H^{-1}(\Omega)$ defined by $A_{\gamma} u := u-\gamma\Delta u$. Using standard results of convex analysis (see, e.g., Ekeland~and~T\'emam~\cite{Ekel99}), we can show that the dual problem to (\ref{primal}) is
\begin{equation}
\label{dual} 
\begin{split}
&\min_{g \in C^1_c(\Omega; \mathbb{R}^2)} \left(
 \dfrac{1}{2}\int_{\Omega} A_{\gamma}^{-1}(f-\beta \nabla\cdot g)\,   (f-\beta \nabla\cdot g) \,dx \right) \\[0.2cm]
&\qquad\qquad\qquad\qquad\text{subject to}\ \ |g|_{\infty}\leq 1,
\end{split}
\end{equation}
where $g = (g_1, g_2)$ is a vector function  and $|g|_{\infty}:= \max\{\abs{g_1},\abs{g_2}\}$.

From the Karush–-Kuhn–-Tucker conditions (see, e.g., Ciarlet \cite[Theorem 9.2-4]{Ciar89}), 
we get that there exist constants $\mu_1$, $\mu_2\geq 0$, such that
$$\left( A_{\gamma}^{-1}(f-\beta\nabla\cdot g)\right)_{x_k}-\mu_k g_k =0\,, \quad k=1,2\,,$$
with either $\mu_k>0$ and $\abs{g_k} =1$ or $\mu_k=0$ and $|g_k|<1$ for $k=1,2$.
In any case, we have that $\mu_1 = \abs{u_{x_1}}$ and $\mu_2 = \abs{u_{x_2}}$,
and therefore, we conclude that the solution $u$ to problem (\ref{primal})
can be found by solving the system of equations
\begin{equation}
\label{sys_eq_cont}
\left\{\begin{array}{rr} 
 A_{\gamma} u= f-\beta\nabla\cdot g \,,\\[0.1cm]
 -u_{x_k}+\abs{u_{x_k}} g_k=0\,,  & k=1,2\,.
\end{array}\right.
\end{equation}

In order to introduce the algorithm to solve (\ref{sys_eq_cont}), we need to
turn into the discrete setting. From now on let $\Omega
=(-L,L)^2\subset\mathbb{R}^2$ and values of the initial data $f$ be given in the
discrete set of $N^2$ uniformly distributed points in $\Omega$. 
To simplify notation, we can fix the number $N$ and take $L$ such that $N = 2L+1$. Now let $\bar{f}$ be a vector in the Euclidean space $X = \mathbb{R}^{N^2}$, defined by $\bar{f}(\abs{x_2-L}+1 + \abs{x_1+L} N) := f(x_1,x_2)$, for $x_1,\, x_2=-L,-L+1,...,L-1,L$, and let us define vectors $\bar{g_1}$, $\bar{g_2}$ and $\bar{u}$ in $X$ in a similar way. Using this notation, we can introduce the discrete version of the system (\ref{sys_eq_cont}), given by
\begin{equation}
\label{sys_eq_discr}
\left\{\begin{array}{rr} 
 \bar{A}_{\gamma} \bar{u}= \bar{f}-\beta\sum_{k=1}^2 D_k \bar{g}_k
 \,,\\[0.1cm]
- D_k \bar{u} + \abs{D_k \bar{u}} \bar{g}_k = 0\,, &  k=1,2\,,
\end{array}\right.
\end{equation}
where $\bar{A}_{\gamma}\in Y$ with $Y=\mathbb{R}^{N^2 \times N^2}$ is a discrete version of the operator $A_{\gamma}$ derived by the standard finite difference scheme taking into account the Neumann boundary conditions and $(D_1, D_2)\in Y\times Y$ corresponds to the discrete version of the gradient operator. To solve the last equations in (\ref{sys_eq_discr}), we follow Chambolle \cite{Cham04} and propose the fixed point iteration
$$\bar{g}_k^{n}=\bar{g}_k^{n-1}+\tau \left(D_k \bar{u}^{n}-\abs{D_k \bar{u}^{n}}\bar{g}_k^{n}\right)\,, \quad k = 1,2,$$
for $n = 1,2,...$.
Finally, the algorithm to solve (\ref{sys_eq_cont}) is given by
\begin{equation}
\label{num_scheme}
\left\{\begin{array}{rr} 
 \bar{A}_{\gamma} \bar{u}^{n}= \bar{f}-\beta \sum_{k=1}^2 D_k \bar{g}_k^{n-1}
 \,,\\[0.1cm]
 \bar{g}_k^{n}=\dfrac{\bar{g}_k^{n-1}+\tau D_k \bar{u}^{n}}{1+\tau \abs{D_k \bar{u}^{n}}}\,, & k=1,2\,,
 \end{array}\right.
\end{equation}
for $n=1,2,...$.

\begin{theorem}
Let $\tau<(8\lambda_1)^{-1}$, where $\lambda_1$ is the smallest eigenvalue of the operator $\bar{A}_{\gamma}$. Then, the sequence $(\bar{u}_n,\bar{g}_n)$ defined  by the scheme (\ref{num_scheme}) converges to the solution $(\bar{u},\bar{g})$ of the equations (\ref{sys_eq_discr}) as $n\rightarrow \infty$.
\end{theorem}

{\bf Proof.} The proof can be carried out in a similar way as in Chambolle~\cite[Theorem 3.1]{Cham04} using the fact that $\bar{A}_{\gamma}$ is a symmetric positive define matrix, what implies that $\lambda_1>0$ and $\langle \bar{A}_{\gamma}^{-1} w,v \rangle = \langle  w, \bar{A}_{\gamma}^{-1} v \rangle$, for all $w,v\in X$. 

\begin{figure}[h!]  
   \begin{center}
    \setlength{\fboxsep}{0pt}
    \setlength{\fboxrule}{0.5pt}
    \fbox{\includegraphics[scale=0.28]{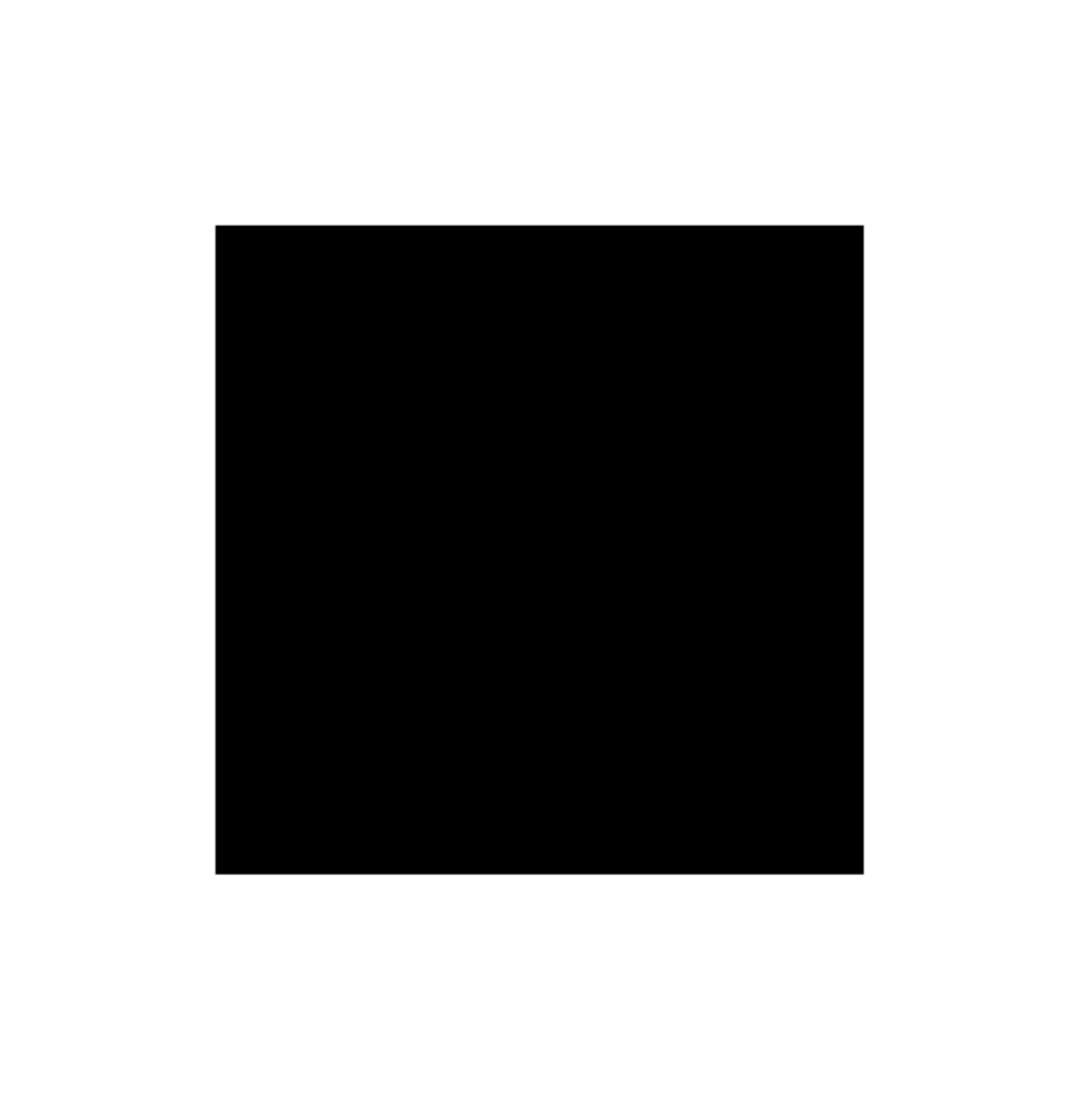}}\hspace{0.2cm}
    \fbox{\includegraphics[scale=0.28]{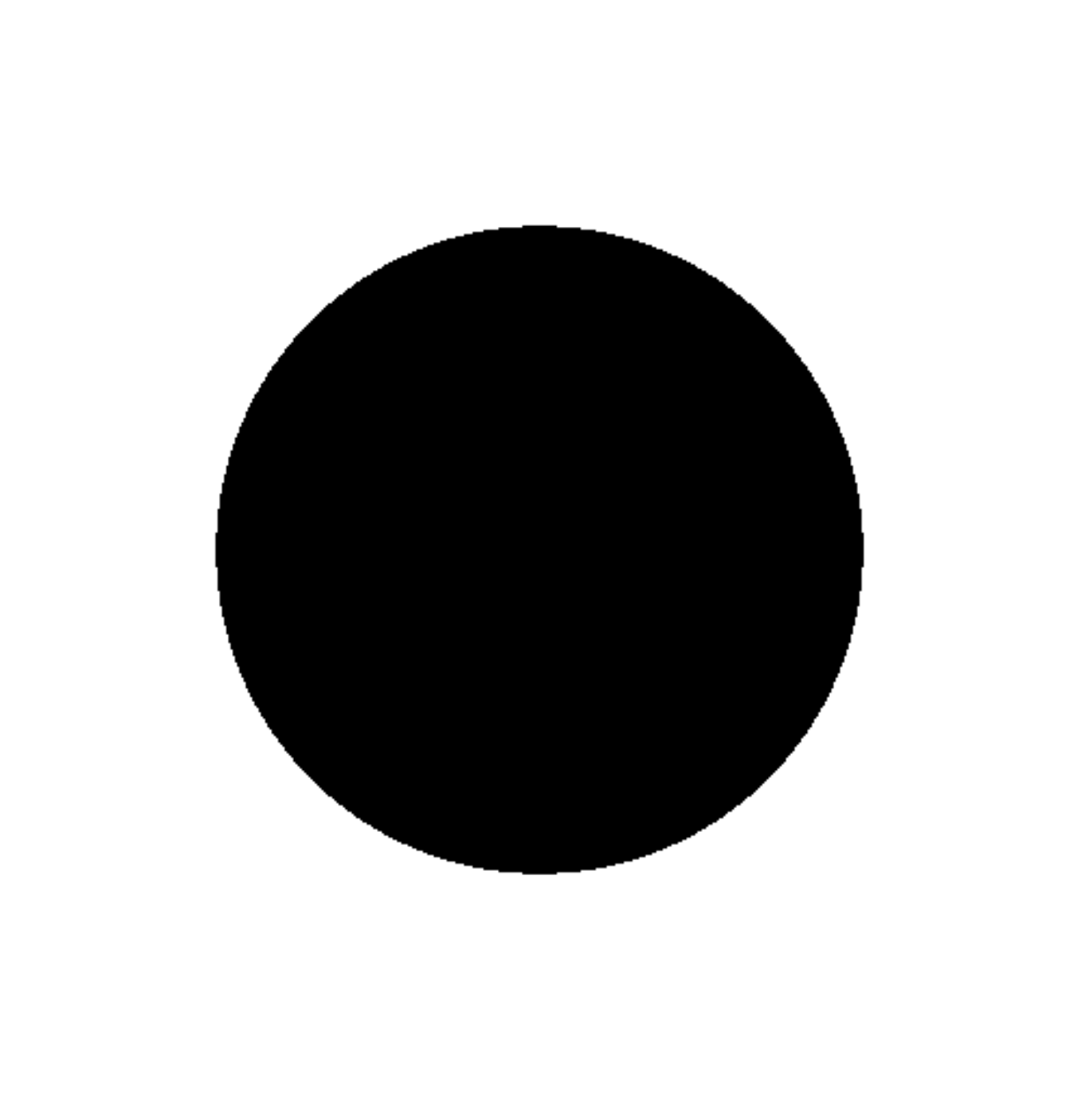}}\hspace{0.2cm}
    \fbox{\includegraphics[scale=0.28]{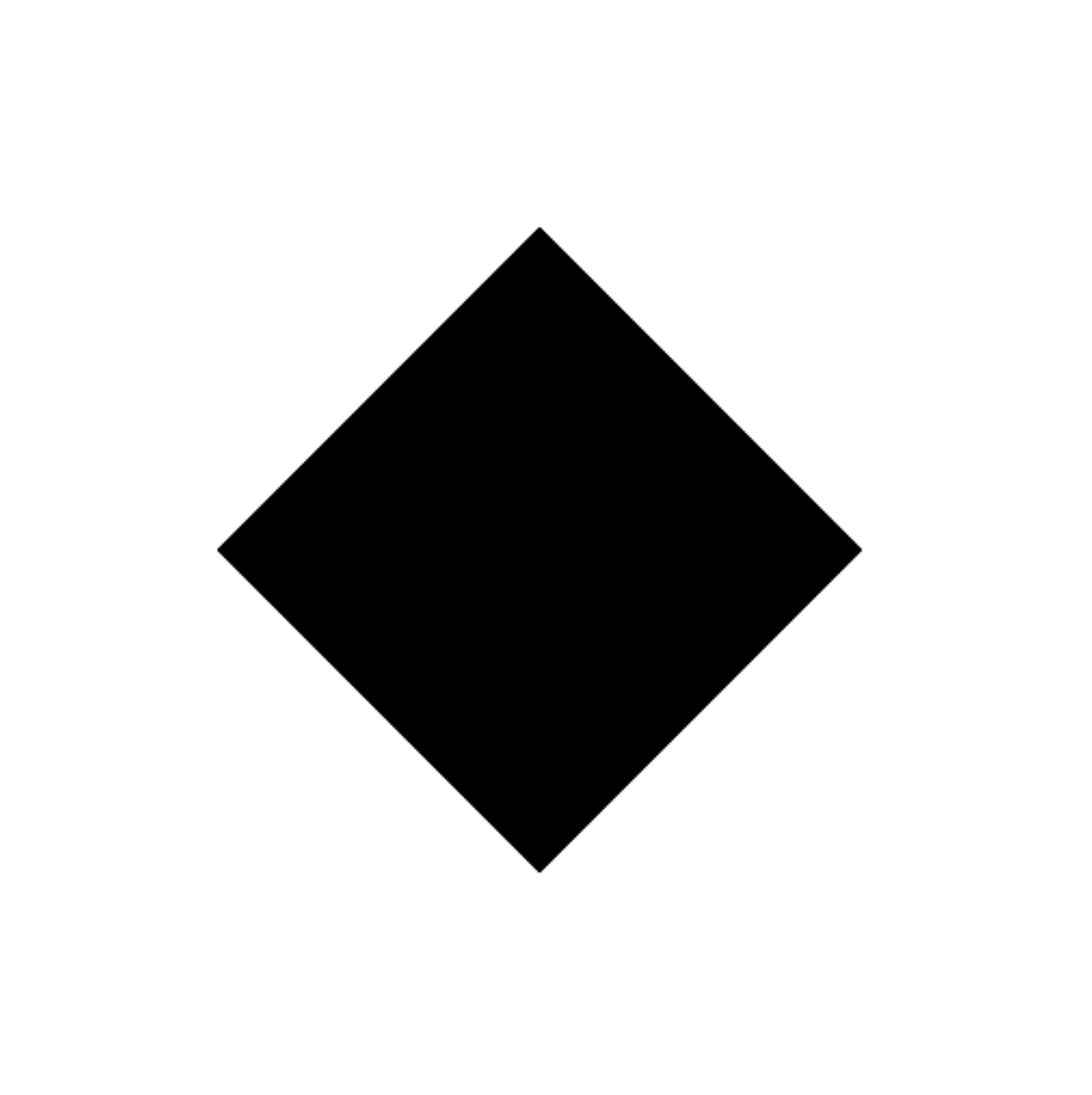}}\hspace{0.2cm}
    \fbox{\includegraphics[scale=0.28]{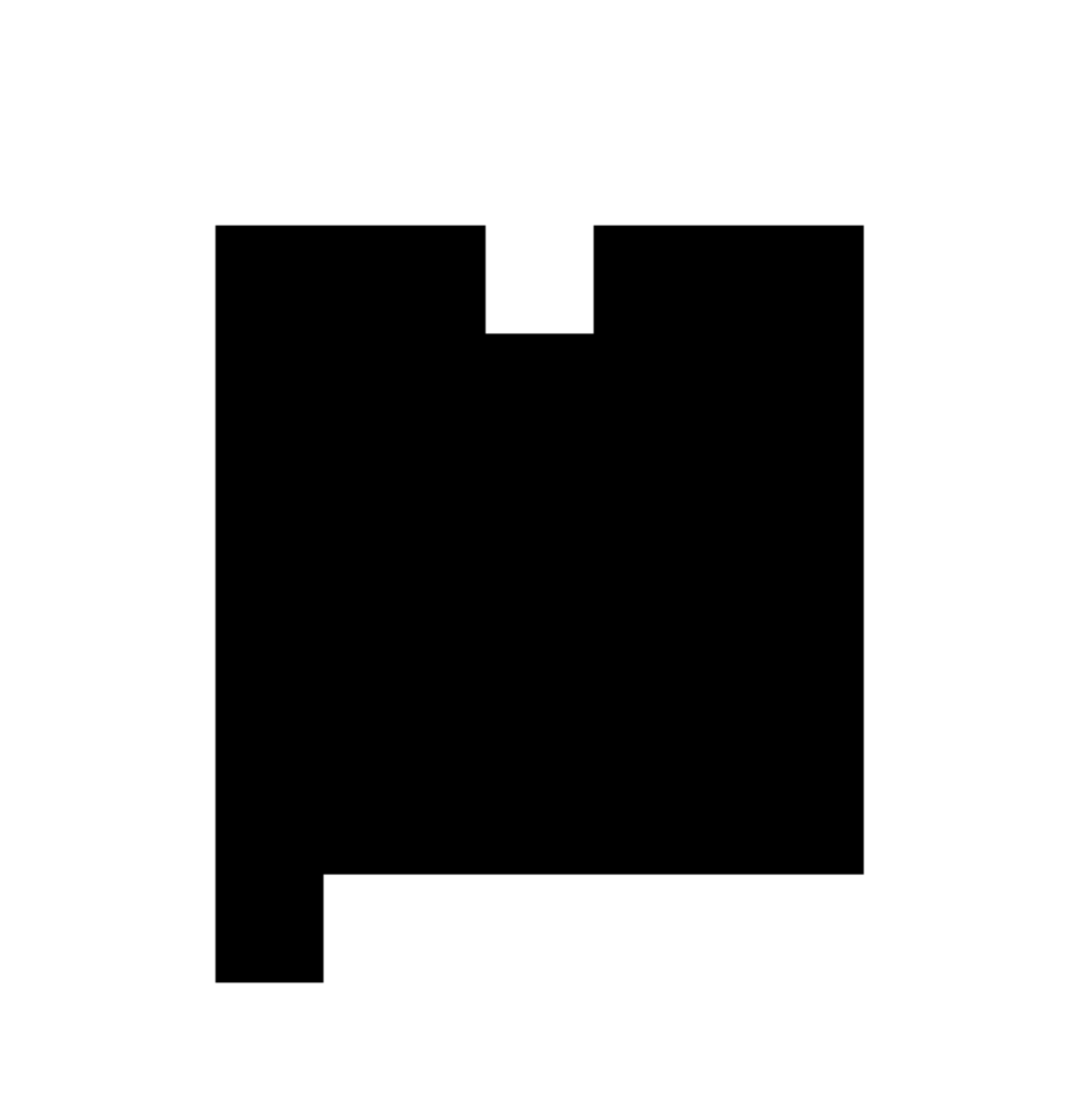}}
    \caption{Images $f_{S_1}$,  $f_{S_2}$, $f_{S_3}$ and $f_{S_4}$.}
    \label{fig1}
    \end{center}
\end{figure}

In the further part of this section, we present numerical solutions to the equations (\ref{eq_anisoTV}) and (\ref{eq_diffanisoTV}) with the Neumann boundary conditions and the initial data $f_{S} = -M \chi_{S}$, where $\chi_{S}:(-L,L)^2\rightarrow\{0,1\}$ is a characteristic function of the set $S \subset (-L,L)^2$. For experiments, we have taken $L=250$, $M=50$ and considered the following four sets:
\begin{equation}
\begin{array}{l} 
  S_1 = \{ x\in \mathbb{R}^2 : \norm{x}_1\leq 150\}\,, 
  S_2 = \{ x\in \mathbb{R}^2 : \norm{x}_2\leq 150\}\,,
  S_3 = \{ x\in \mathbb{R}^2 : \norm{x}_{\infty}\leq 150\}\,,\\[0.1cm]
  S_4 = (S_1\cup \{ x\in \mathbb{R}^2 : \norm{x-(125,175)}_1\leq 25\})\setminus \{ x\in \mathbb{R}^2 : \norm{x+(0,125)}_1\leq 25\}\,.
\end{array}
\nonumber
\end{equation}
Images $f_{S_1}$,  $f_{S_2}$, $f_{S_3}$ and $f_{S_4}$ are presented in Figure \ref{fig1}. 

All experiments were performed with the same values for parameters involved in
the algorithm, i.e., $\gamma=5^{-1}$, $\beta=10$, $\delta t = 1$, $\tau=8^{-1}$.
As the stopping criterion for the iterative scheme (\ref{num_scheme}), we have
used
$\norm{\bar{u}^{n-1}-\bar{u}^{n}}_2\norm{\bar{u}^{n}}_2^{-1}<tol$,
with the tolerance $tol = 10^{-5}$.

\begin{figure}[h!]  
   \begin{center}
    \setlength{\fboxsep}{0pt}
    \setlength{\fboxrule}{0.5pt}
    \fbox{\includegraphics[scale=0.28]{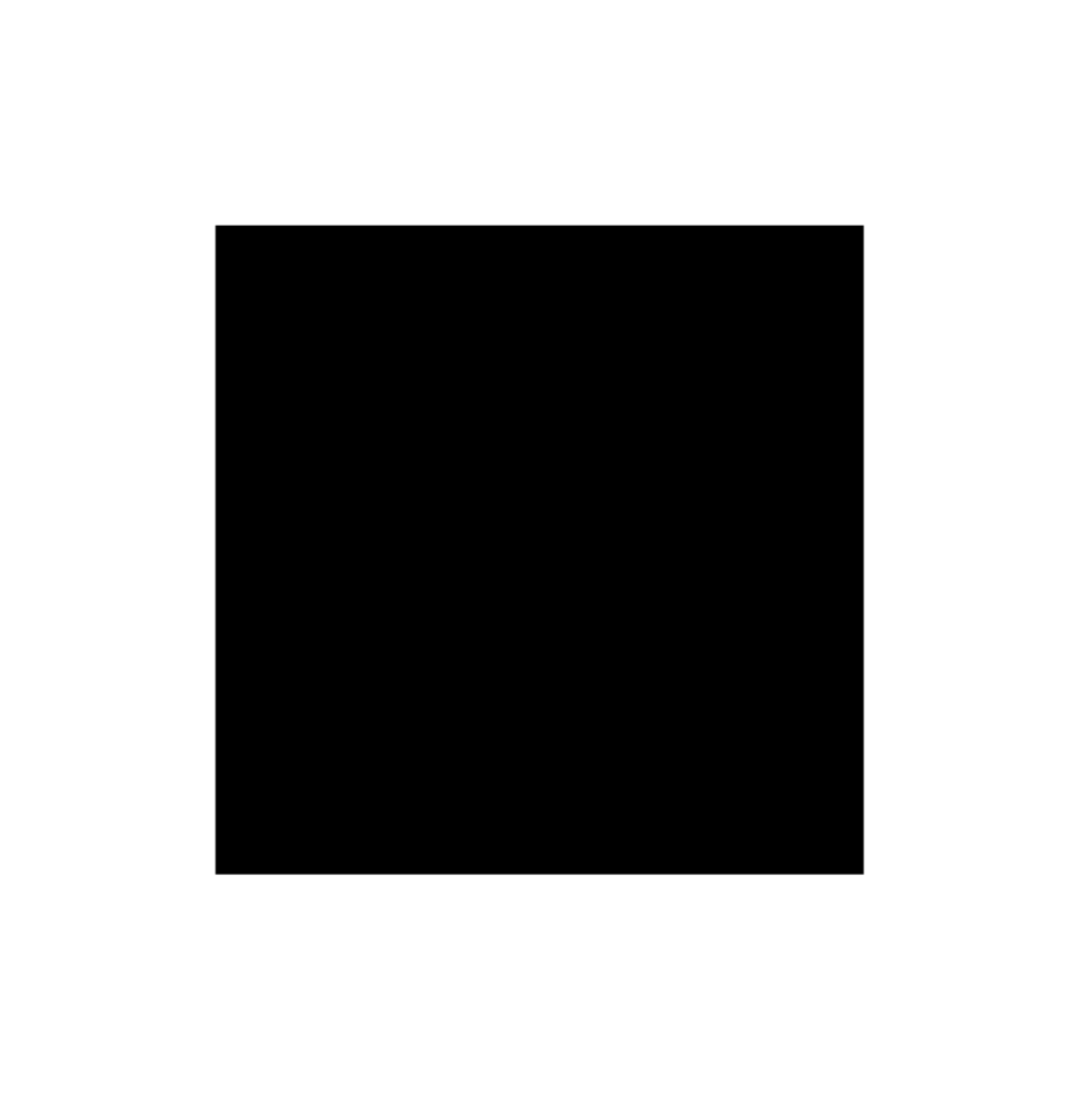}}\hspace{0.2cm}
    \fbox{\includegraphics[scale=0.28]{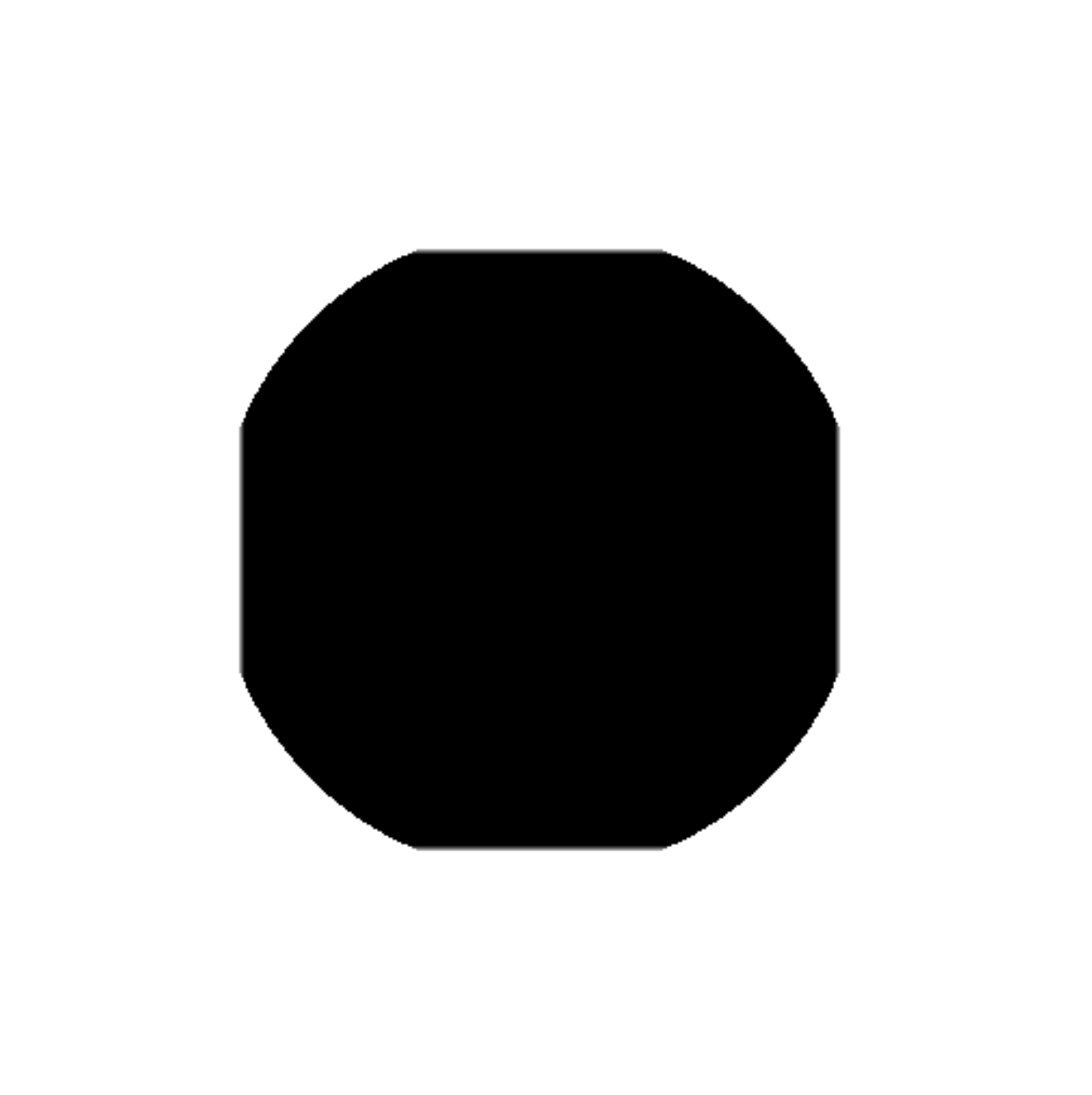}}\hspace{0.2cm}
    \fbox{\includegraphics[scale=0.28]{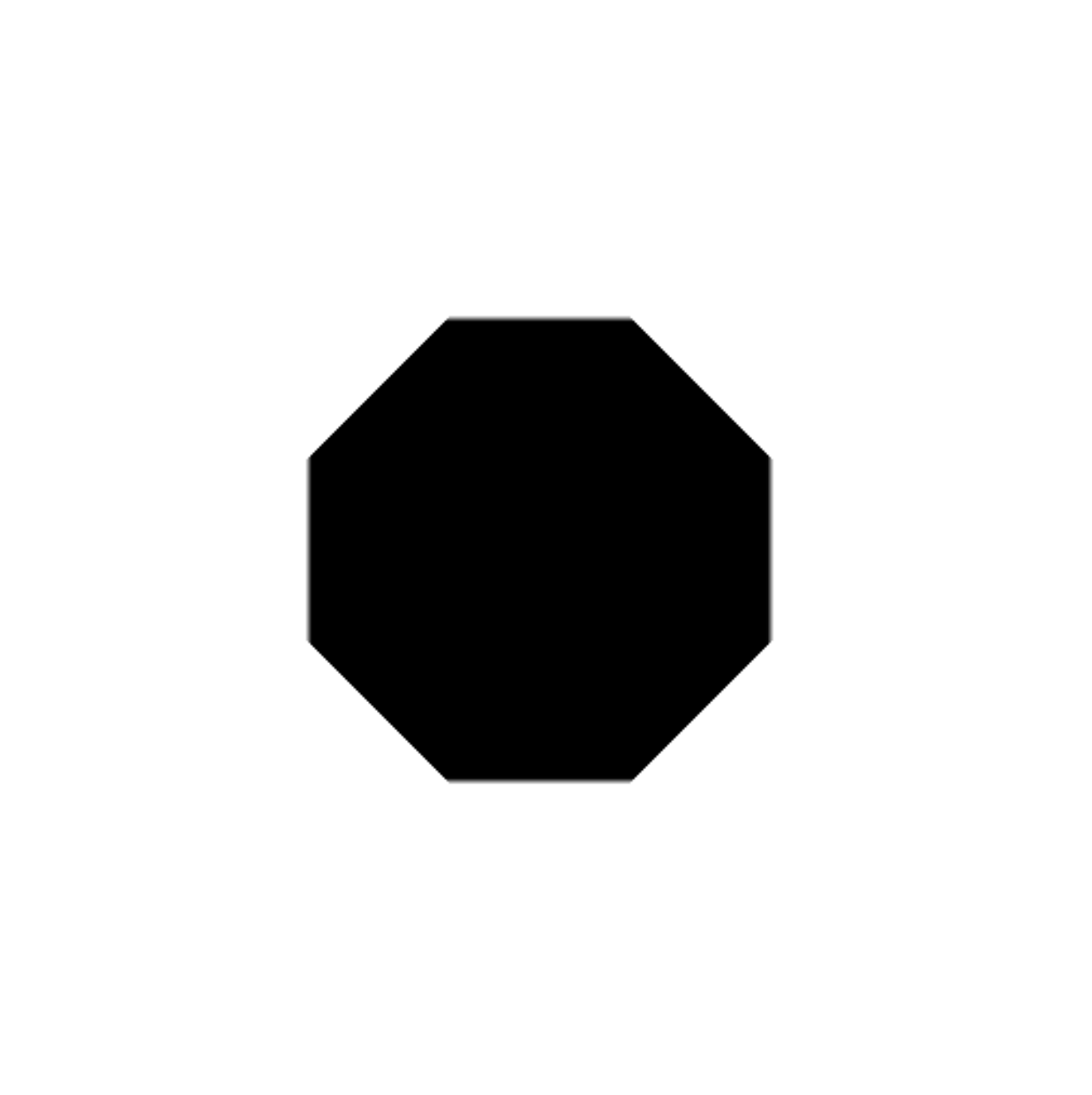}}\hspace{0.2cm}
        \fbox{\includegraphics[scale=0.28]{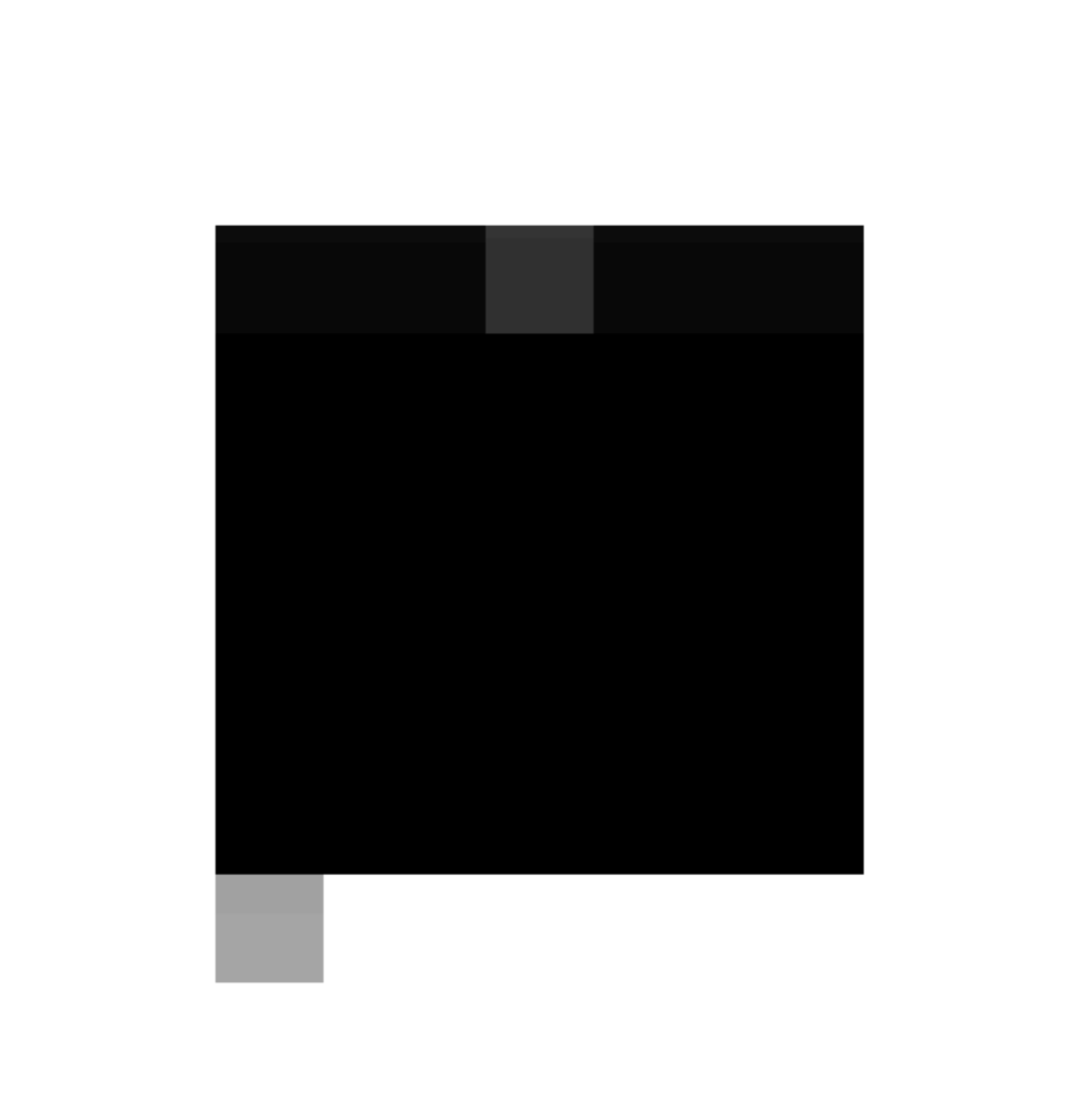}}\\[0.2cm]
    \fbox{\includegraphics[scale=0.28]{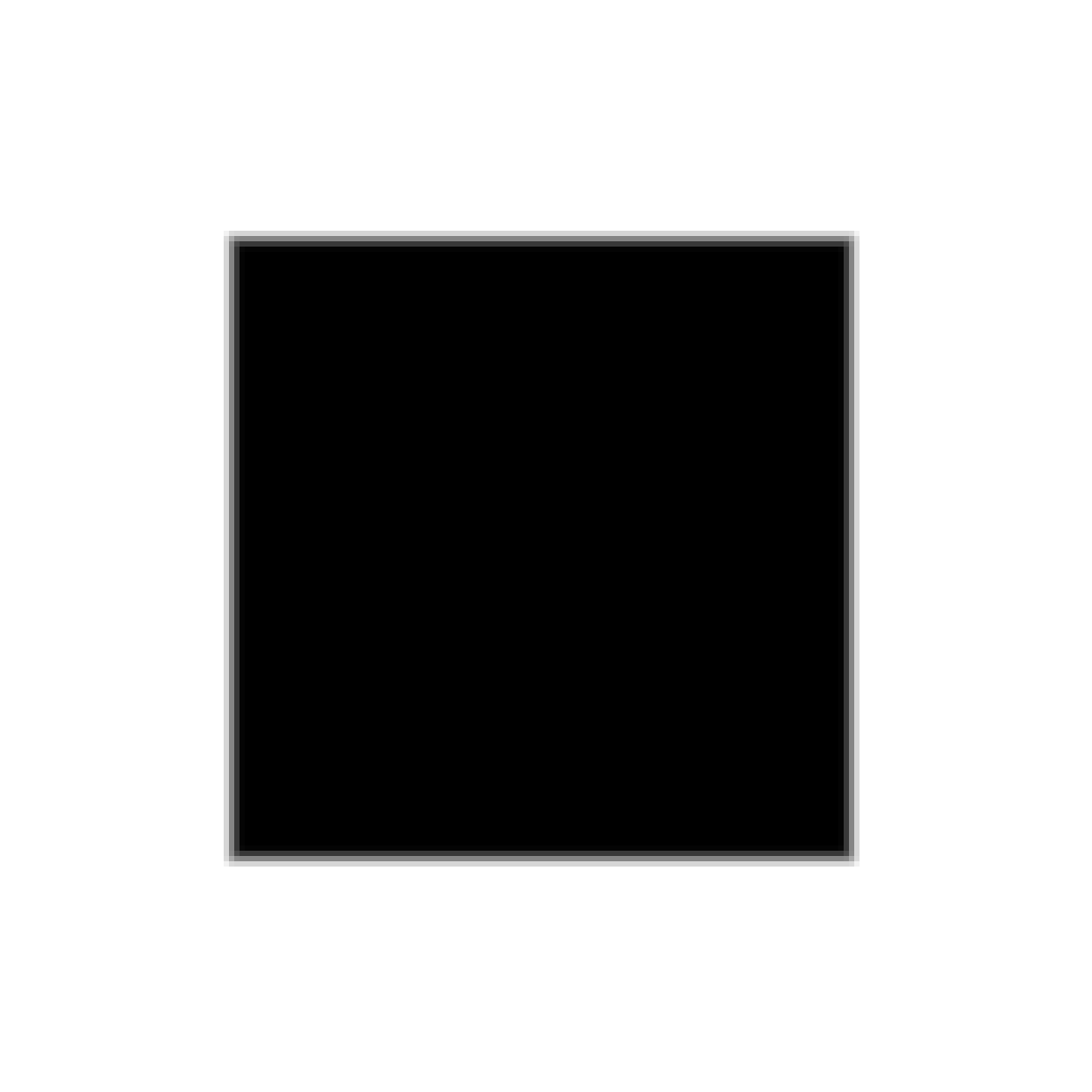}}\hspace{0.2cm}
    \fbox{\includegraphics[scale=0.28]{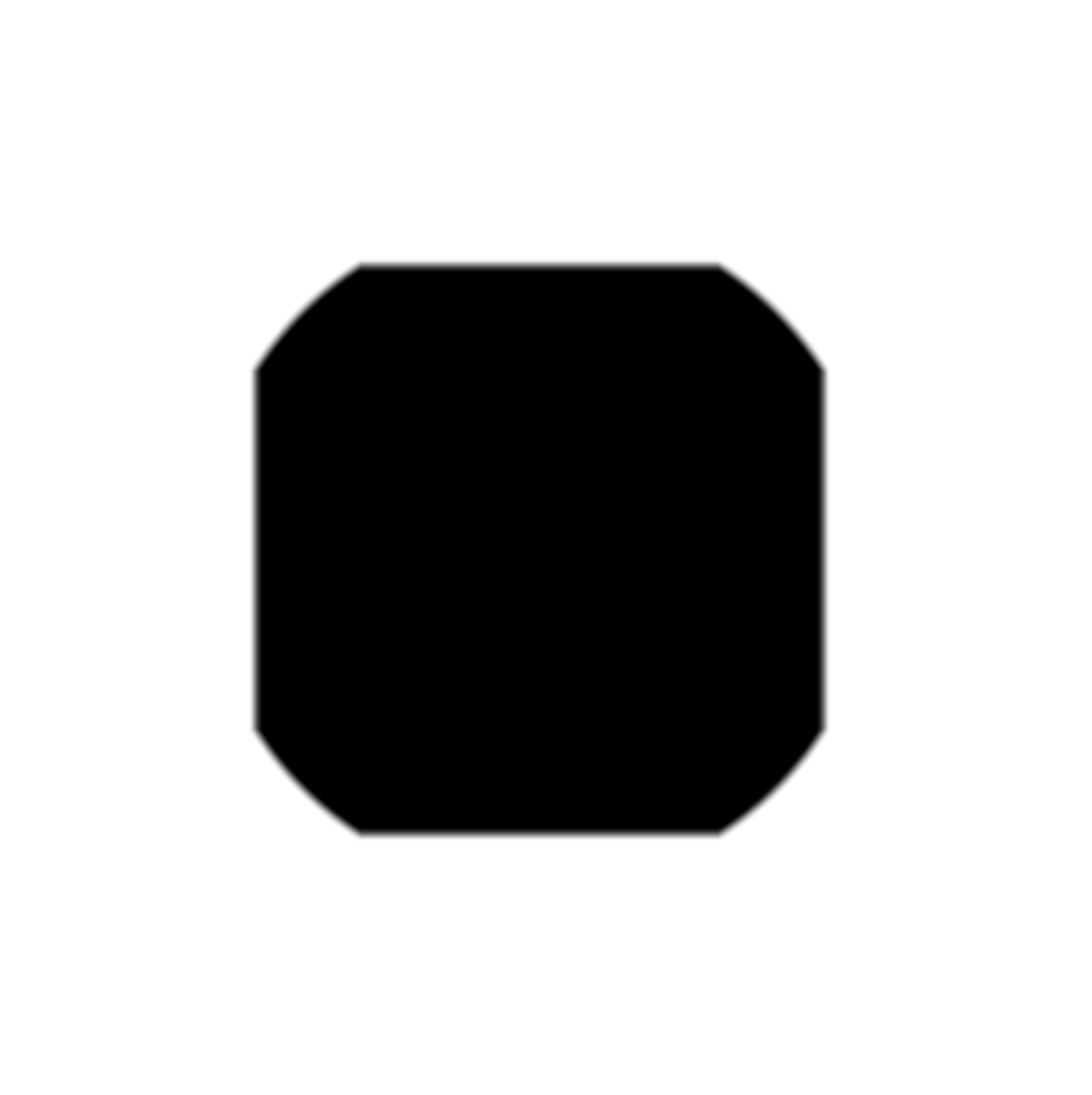}}\hspace{0.2cm}
    \fbox{\includegraphics[scale=0.28]{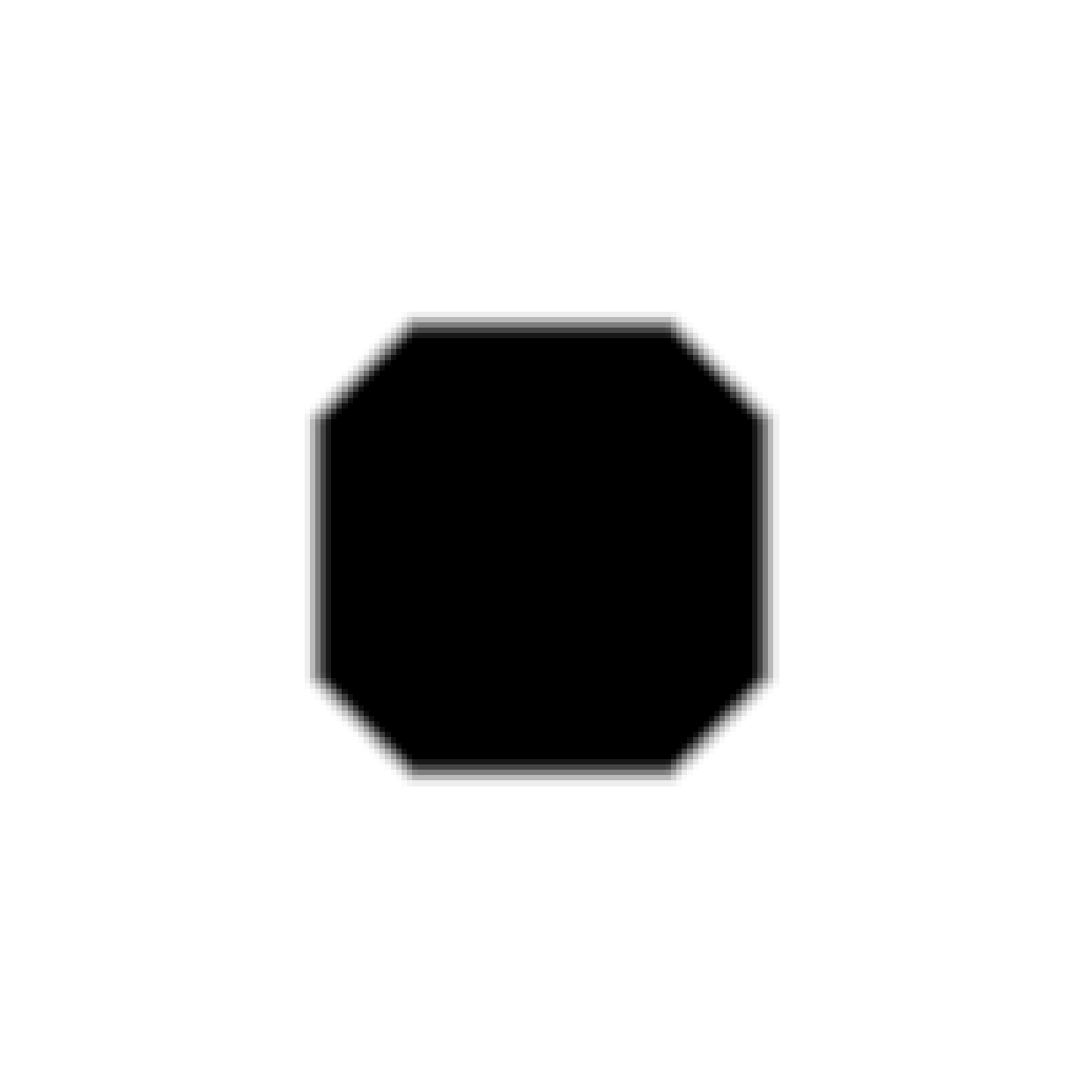}}\hspace{0.2cm}
    \fbox{\includegraphics[scale=0.28]{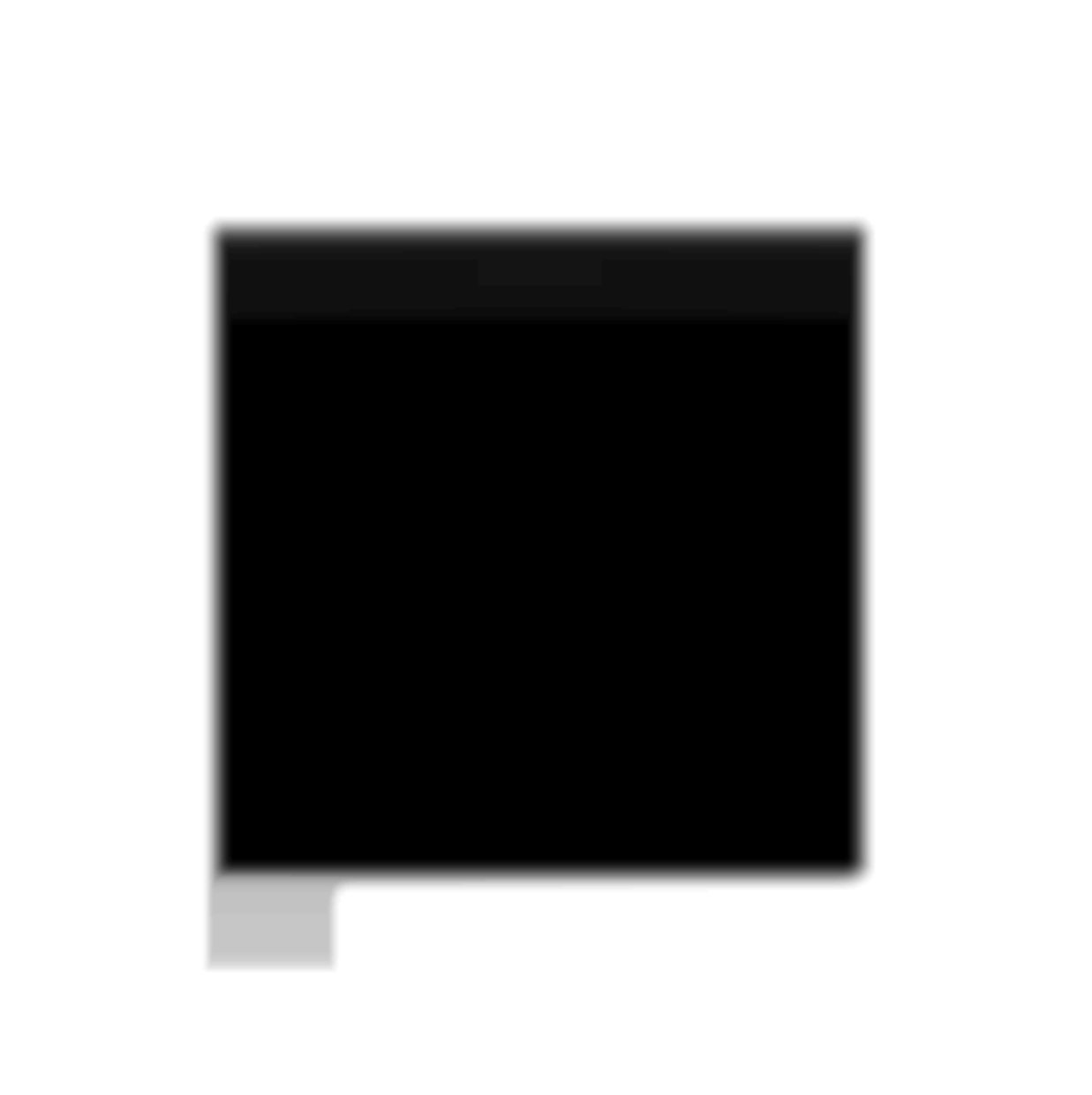}}
    \caption{Numerical solutions to the equations (\ref{eq_anisoTV}) (upper row) and (\ref{eq_diffanisoTV}) (lower row) with the initial data $f_{S_1}$,  $f_{S_2}$, $f_{S_3}$ and $f_{S_4}$, respectively.}
    \label{fig2}
    \end{center}
\end{figure}

Numerical solutions to the equation (\ref{eq_anisoTV}) with the initial data
accordingly equal to $f_{S_1}$,  $f_{S_2}$, $f_{S_3}$ and $f_{S_4}$ are
presented in the upper row of Figure \ref{fig2}. The first two results have been
obtained for $m=200$ , whereas the next two results, for $m=170$ and $m=90$,
respectively. We recall that $m$ denotes the number of iteration of the scheme
(\ref{semidiscret}). Numerical solutions to the equation (\ref{eq_diffanisoTV}),
with the same initial data and for the same numbers of iterations as before are
presented in the lower row of Figure \ref{fig2}.

\begin{figure}[h!]  
   \begin{center}
    \setlength{\fboxsep}{3pt}
    \setlength{\fboxrule}{0pt}
    \fbox{\includegraphics[scale=0.28]{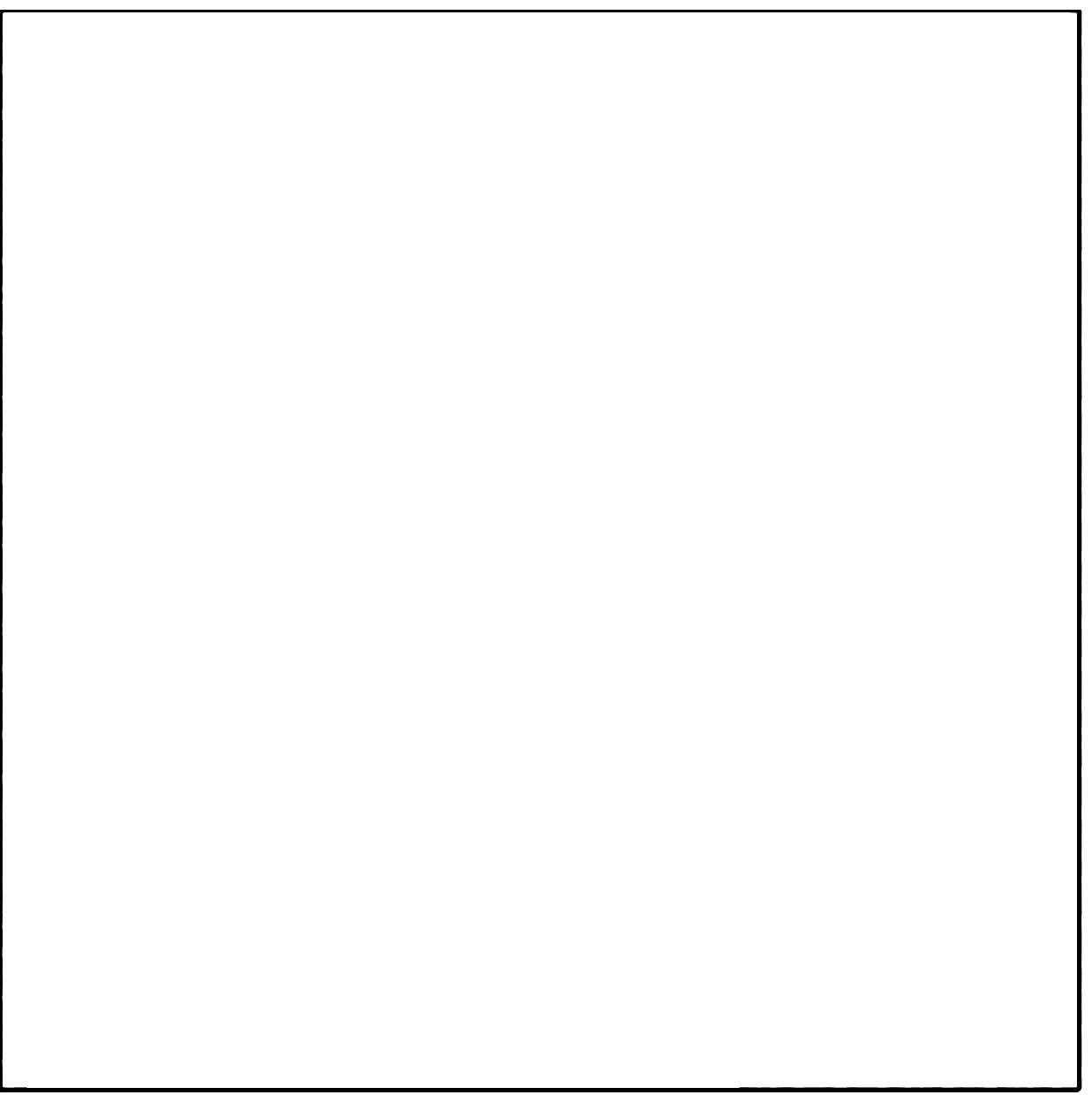}}\hspace{0.2cm}
    \fbox{\includegraphics[scale=0.28]{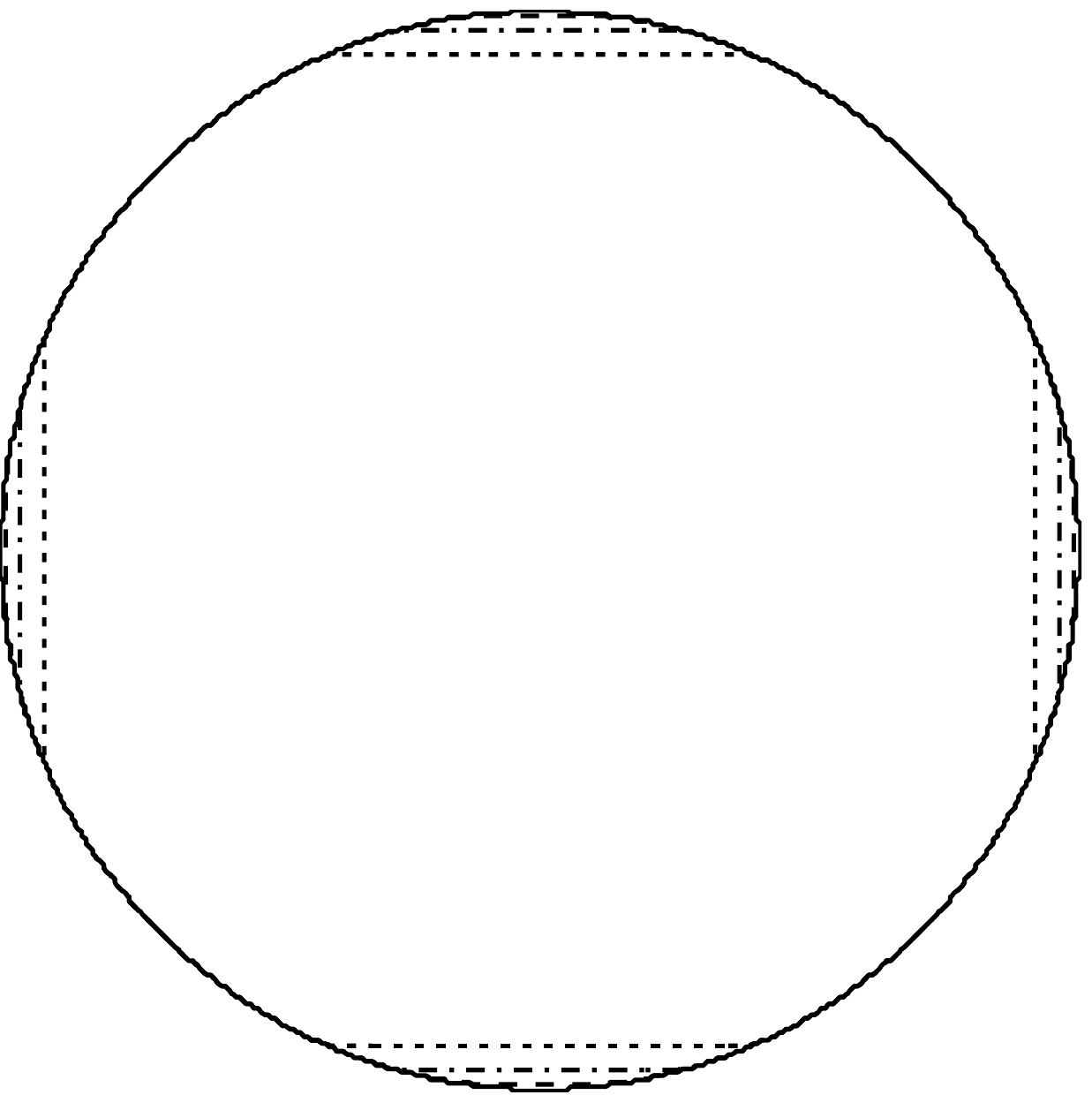}}\hspace{0.2cm}
    \fbox{\includegraphics[scale=0.28]{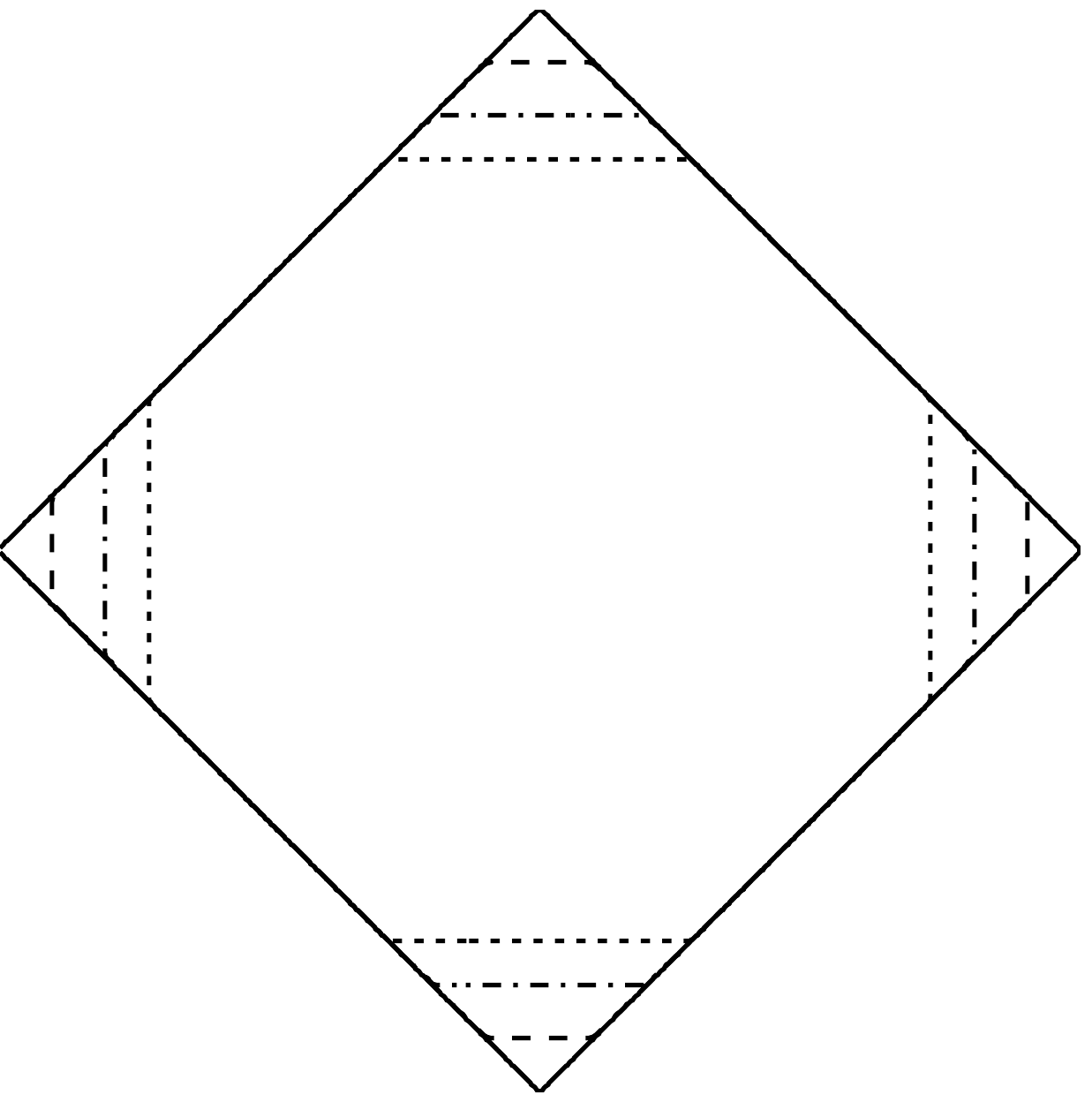}}\hspace{0.2cm}
  \fbox{\raisebox{\depth}{\scalebox{1}[-1]{\includegraphics[scale=0.28]{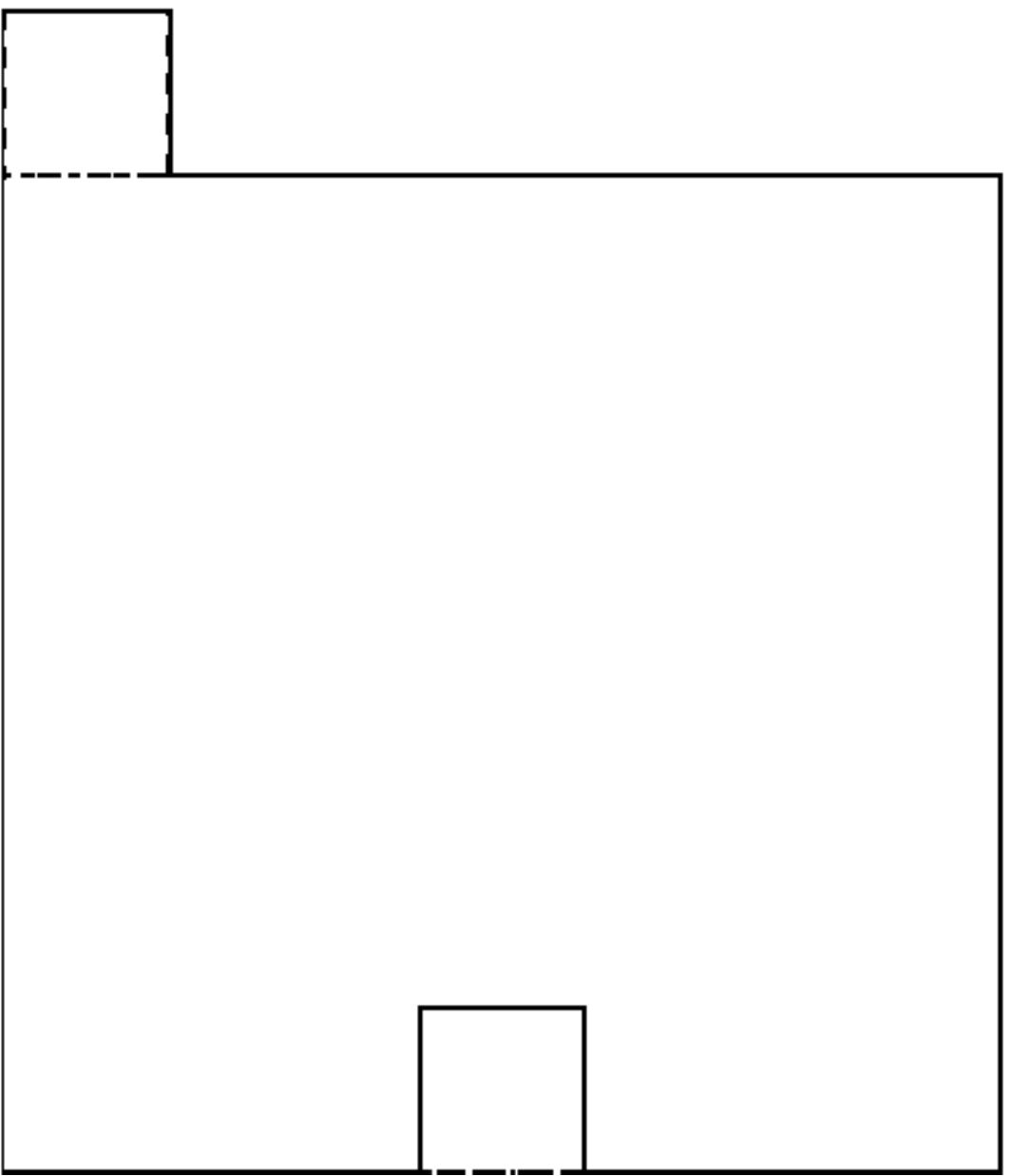}}}}\\[0.2cm]
    \fbox{\includegraphics[scale=0.28]{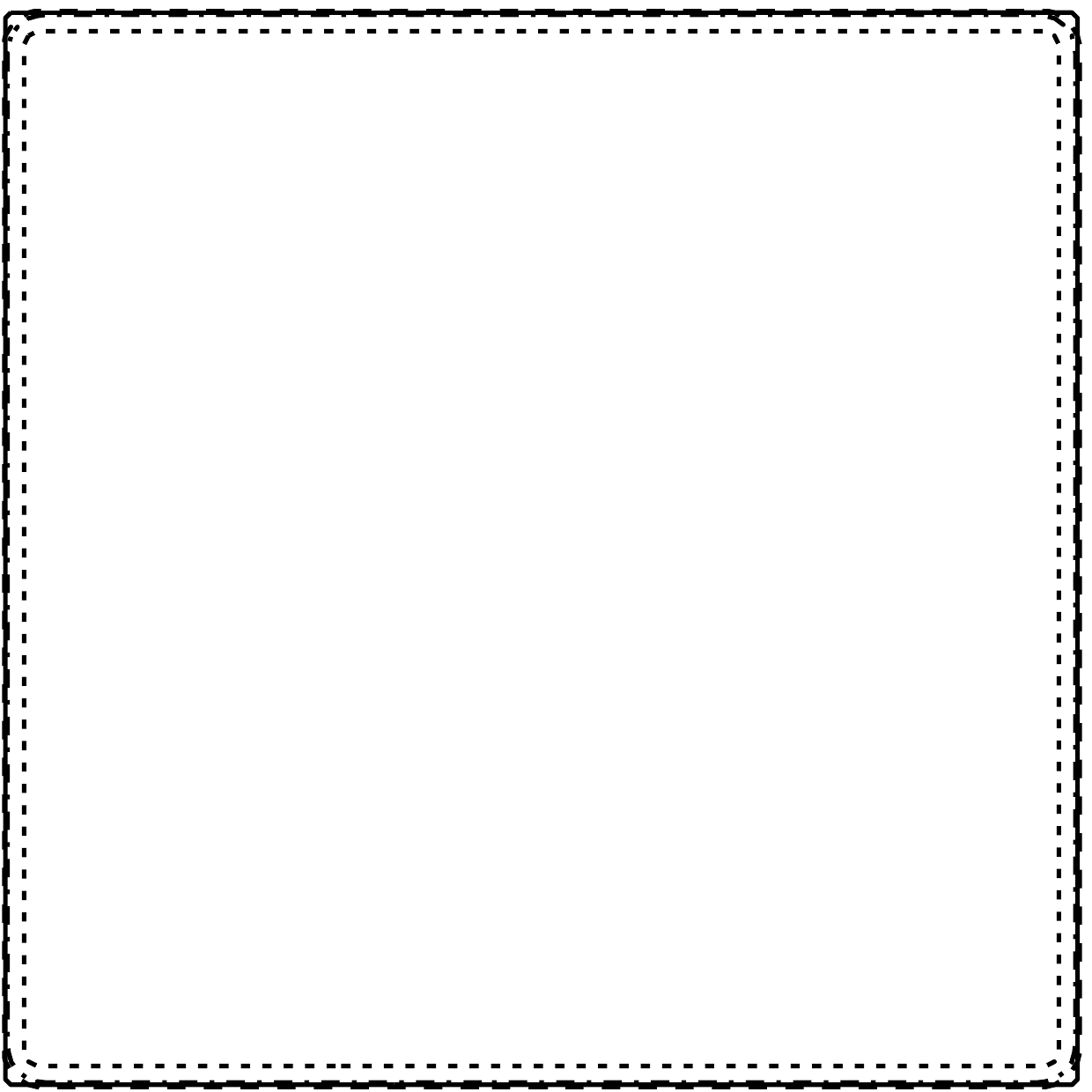}}\hspace{0.2cm}
    \fbox{\includegraphics[scale=0.28]{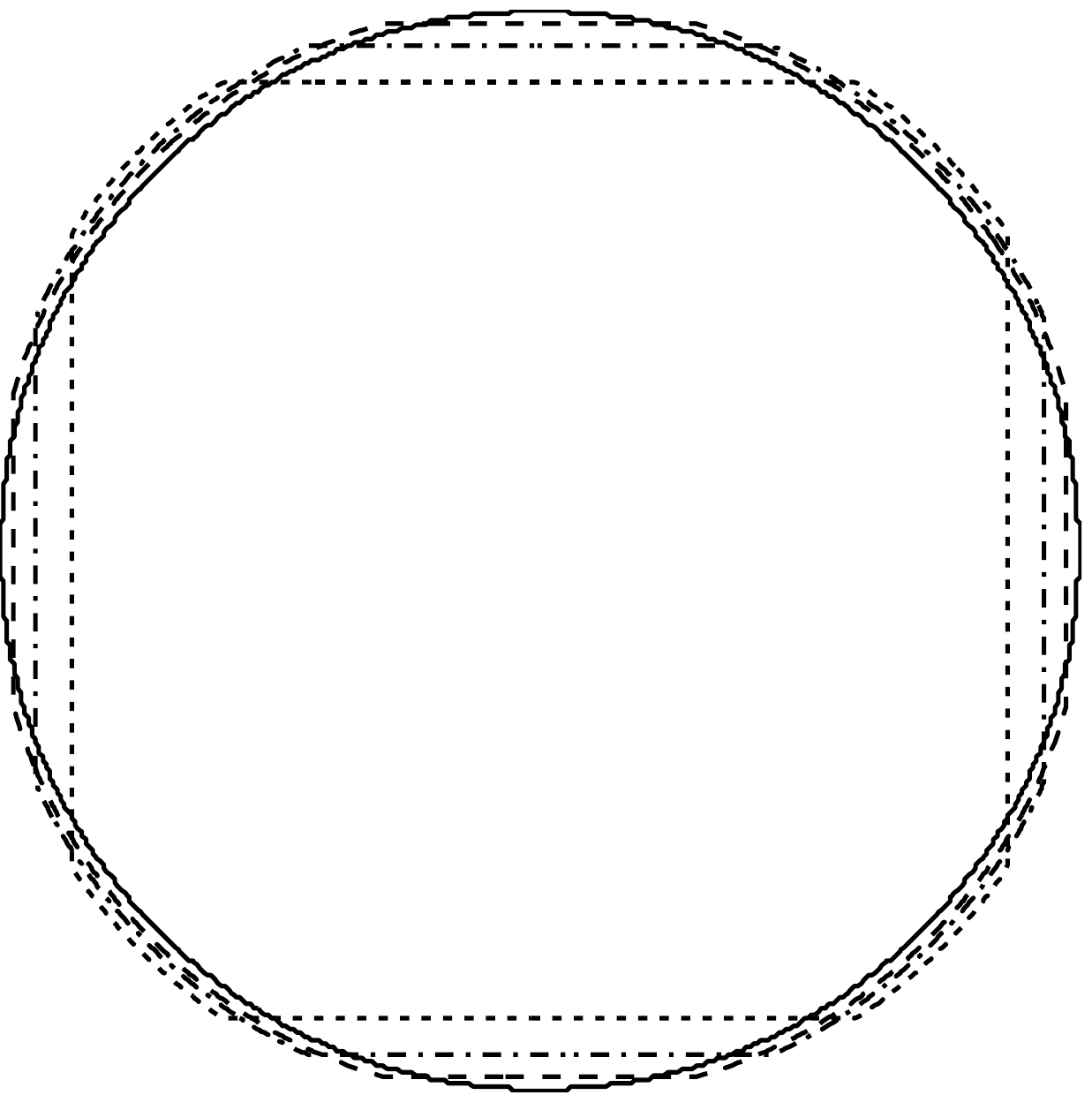}}\hspace{0.2cm}
    \fbox{\includegraphics[scale=0.28]{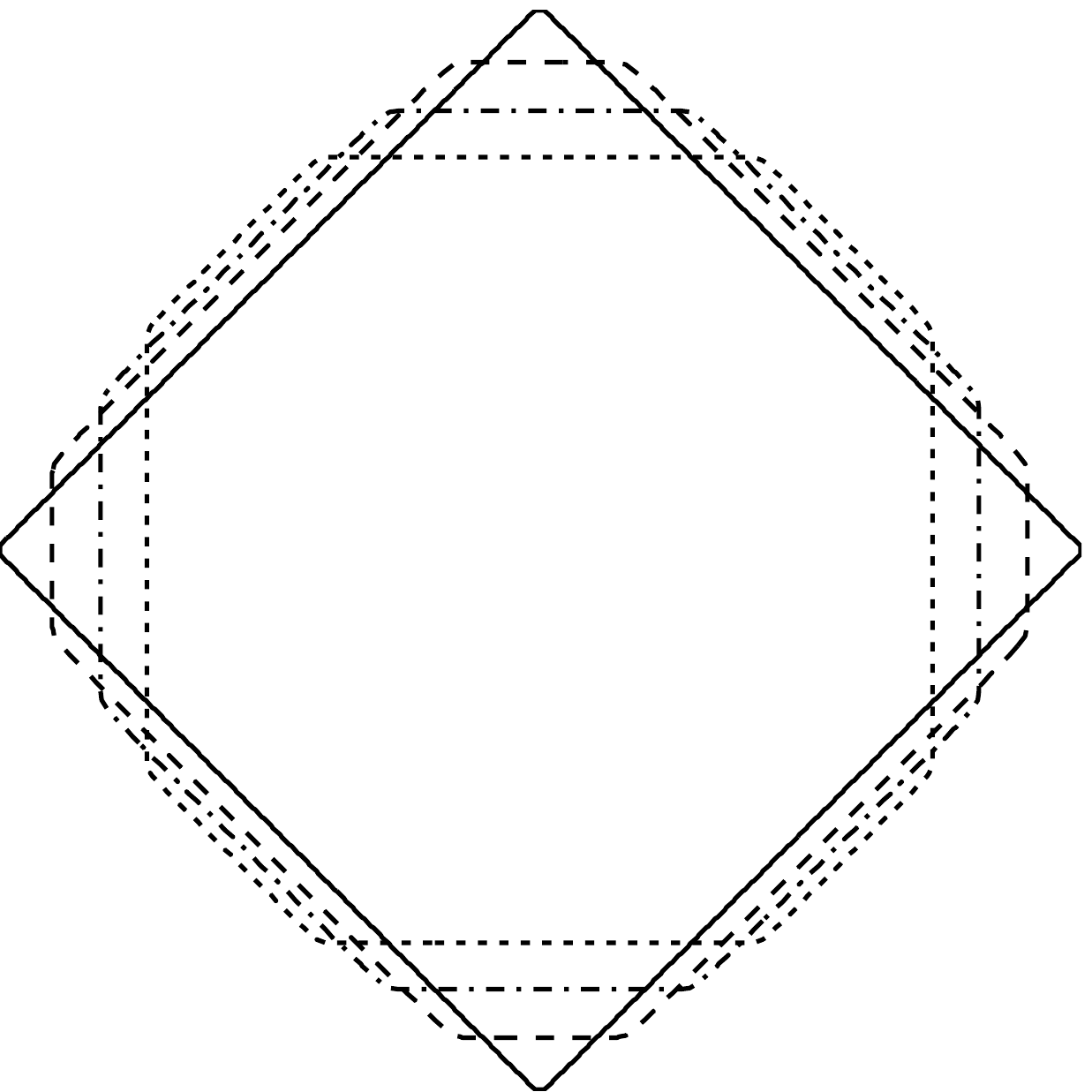}}\hspace{0.2cm}
\fbox{\raisebox{\depth}{\scalebox{1}[-1]{\includegraphics[scale=0.28]{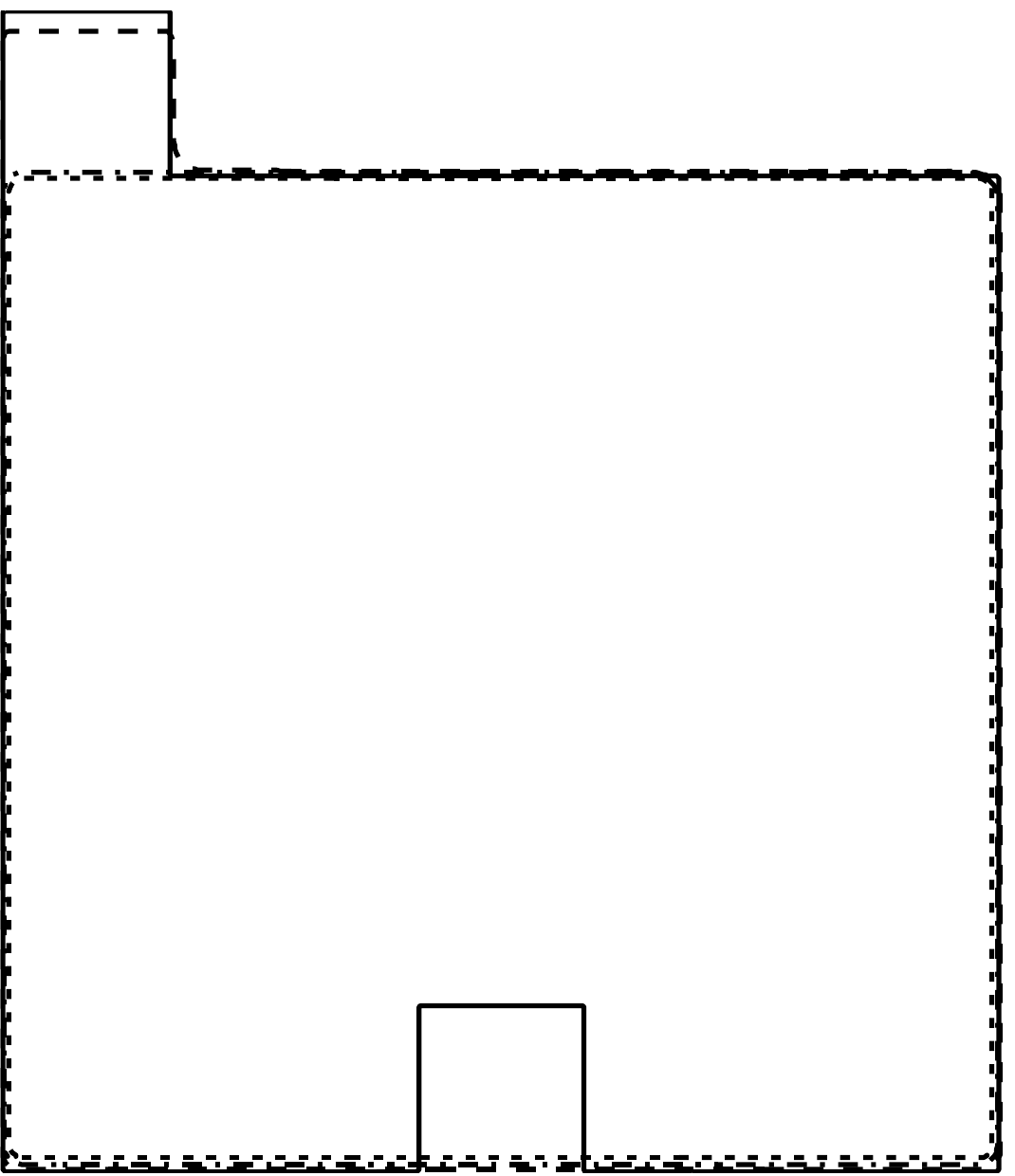}}}}
    \caption{Evolution of contours corresponding to solutions to the equations (\ref{eq_anisoTV}) (upper row) and (\ref{eq_diffanisoTV}) (lower row) with the initial data $f_{S_1}$,  $f_{S_2}$, $f_{S_3}$ and $f_{S_4}$, respectively.}
    \label{fig3}
    \end{center}
\end{figure}

The first two graphs in the upper row of Figure \ref{fig3} present evolution of
contour lines of solutions to the equation (\ref{eq_anisoTV}) with the initial
data $f_{S_1}$ and $f_{S_2}$, respectively. In each graph, contours are plotted
for the level equal to the average value of a given initial data and correspond
to solutions of the equation (\ref{eq_anisoTV}) for $m = 0$, $70$, $140$ and
$210$. The contour lines of solutions to the same equation but with the initial
data $f_{S_3}$ and $f_{S_4}$ and for $m=0$, $60$, $120$ and $170$ are presented
in the next two graphs in the same row. The lower row of Figure \ref{fig3}
presents the evolution of contour lines corresponding to solutions of the
equation (\ref{eq_diffanisoTV}) with the same initial data and for the same
numbers of iterations as before.

\begin{figure}[h!]  
   \begin{center}
   \includegraphics[scale=0.31]{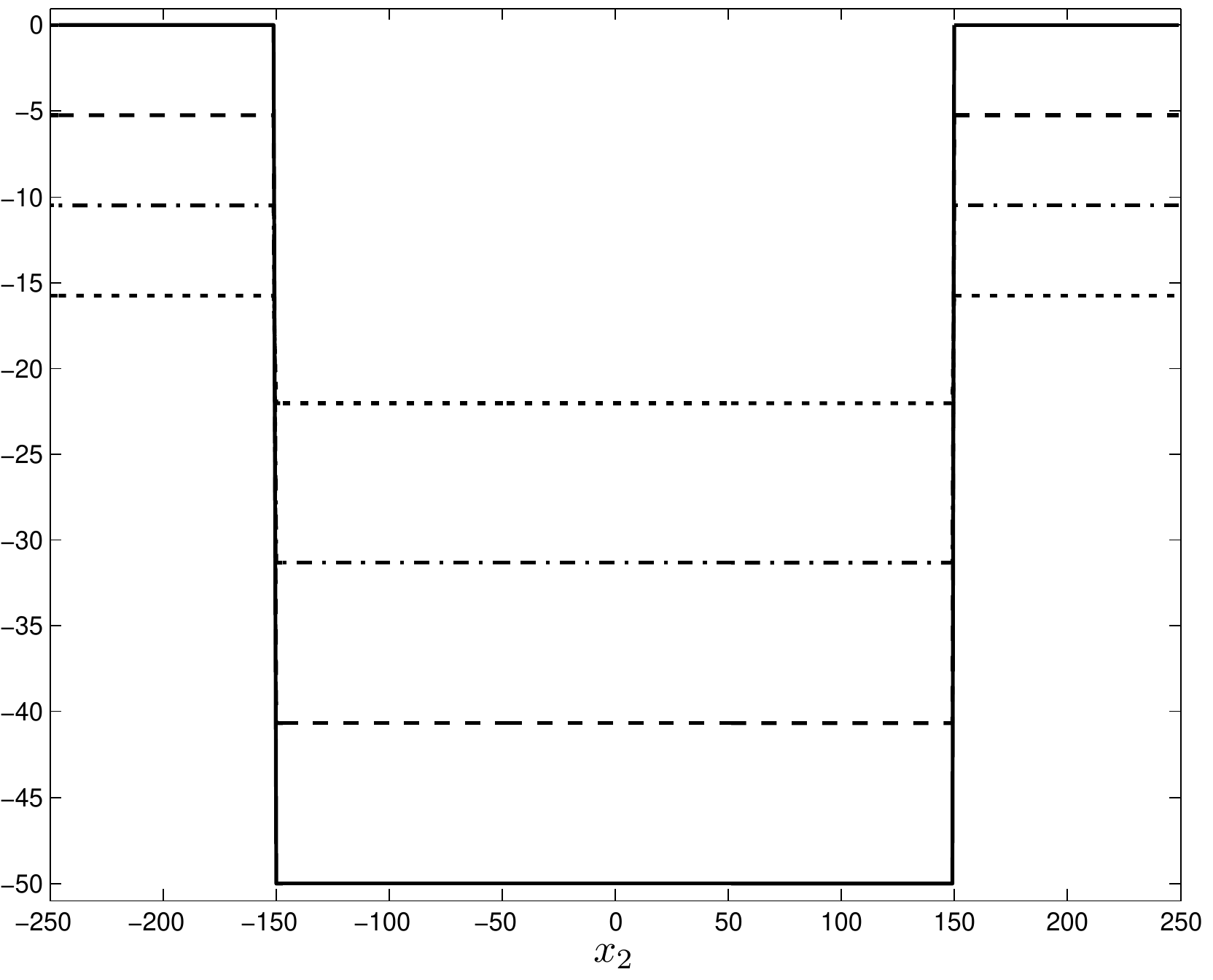}
   \includegraphics[scale=0.31]{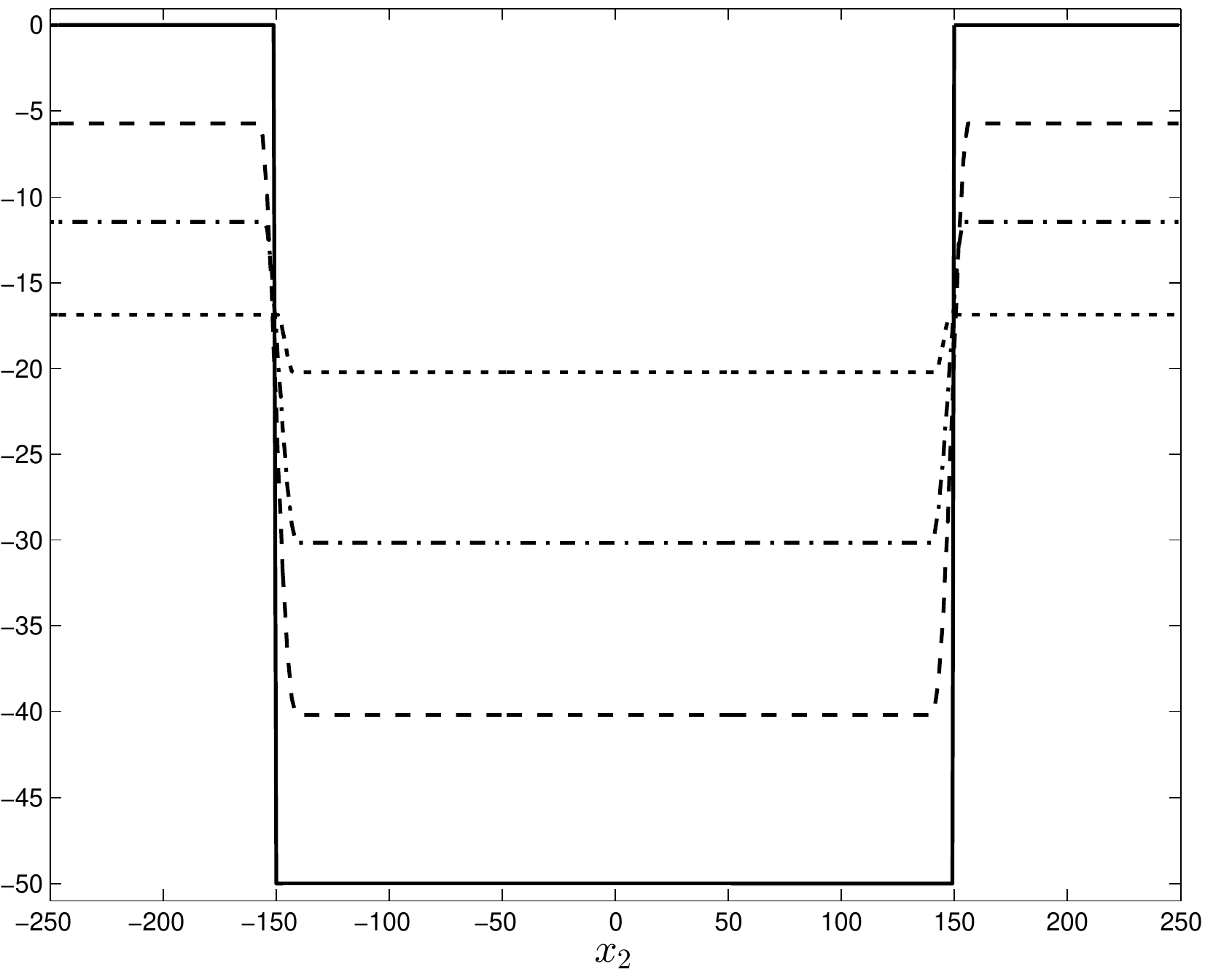}
    \includegraphics[scale=0.31]{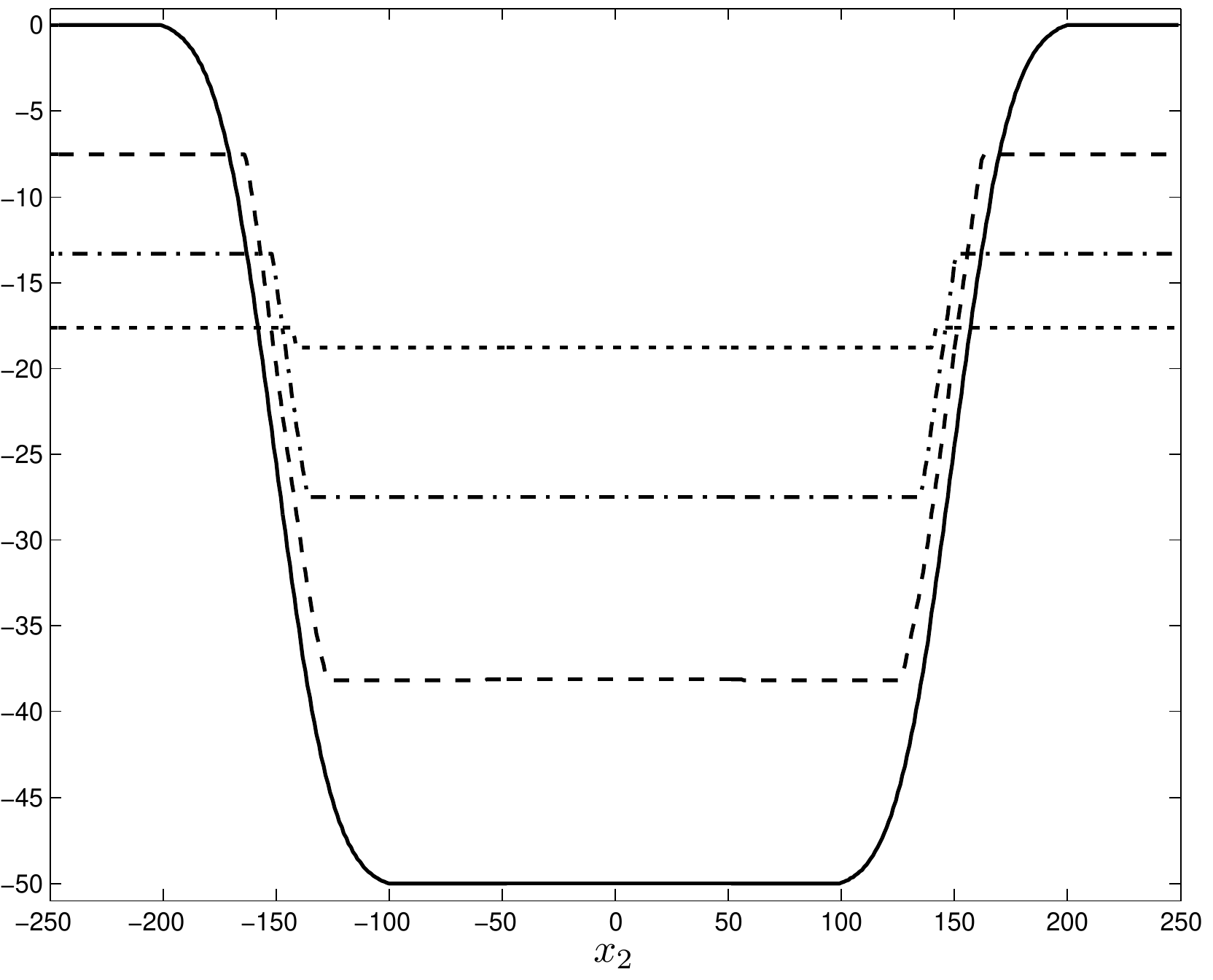}
    \caption{Evolution of numerical solutions to the equations (\ref{eq_anisoTV}) and (\ref{eq_diffanisoTV}) with the initial data $f_{S_1}$ and to the equation (\ref{eq_anisoTV}) with the initial data $G_\sigma\ast f_{S_1}$.}
    \label{fig4}
    \end{center}
\end{figure}

The first two plots in Figure \ref{fig4} show evolution of numerical solutions
to the equations (\ref{eq_anisoTV}) and (\ref{eq_diffanisoTV}), respectively, 
along cross-sectional line $x_1=0$ passing through the middle of the square
$S_1$ for $m=0$, $70$, $140$ and $210$. In the case of the solution $u$ to the
equation (\ref{eq_anisoTV}), obtained values were equal to:
$\{-5.25, -40.67\}$  for $m=70$,  
$\{-10.5, -31.33\}$ for $m=140$,
$\{-15.75, -22\}$ for $m=210$, where the first numbers in brackets correspond to values of $u$ in $\Omega\setminus S_1$, and the second ones, in $S_1$. 
We note that these results coincide with the exact values given by the formula
(\ref{wzor3}) for $t=\beta\,\delta t\,m$. The third plot in Figure \ref{fig4}
presents the evolution of the numerical solution to the equation
(\ref{eq_anisoTV}) for smooth initial data obtained by convolution of the image
$f_{S_1}$ and the Gaussian kernel $G_{\sigma}$ with the standard deviation
$\sigma=10$. Similarly as in the one dimensional version of the equation
(\ref{eq_anisoTV}) studied in \cite{joss}. Here, we may observe
propagation of facets.


%

\begin{figure}[h!]  
   \begin{center}
    \setlength{\fboxsep}{0pt}
    \setlength{\fboxrule}{0.5pt}
    \fbox{\includegraphics[scale=0.28]{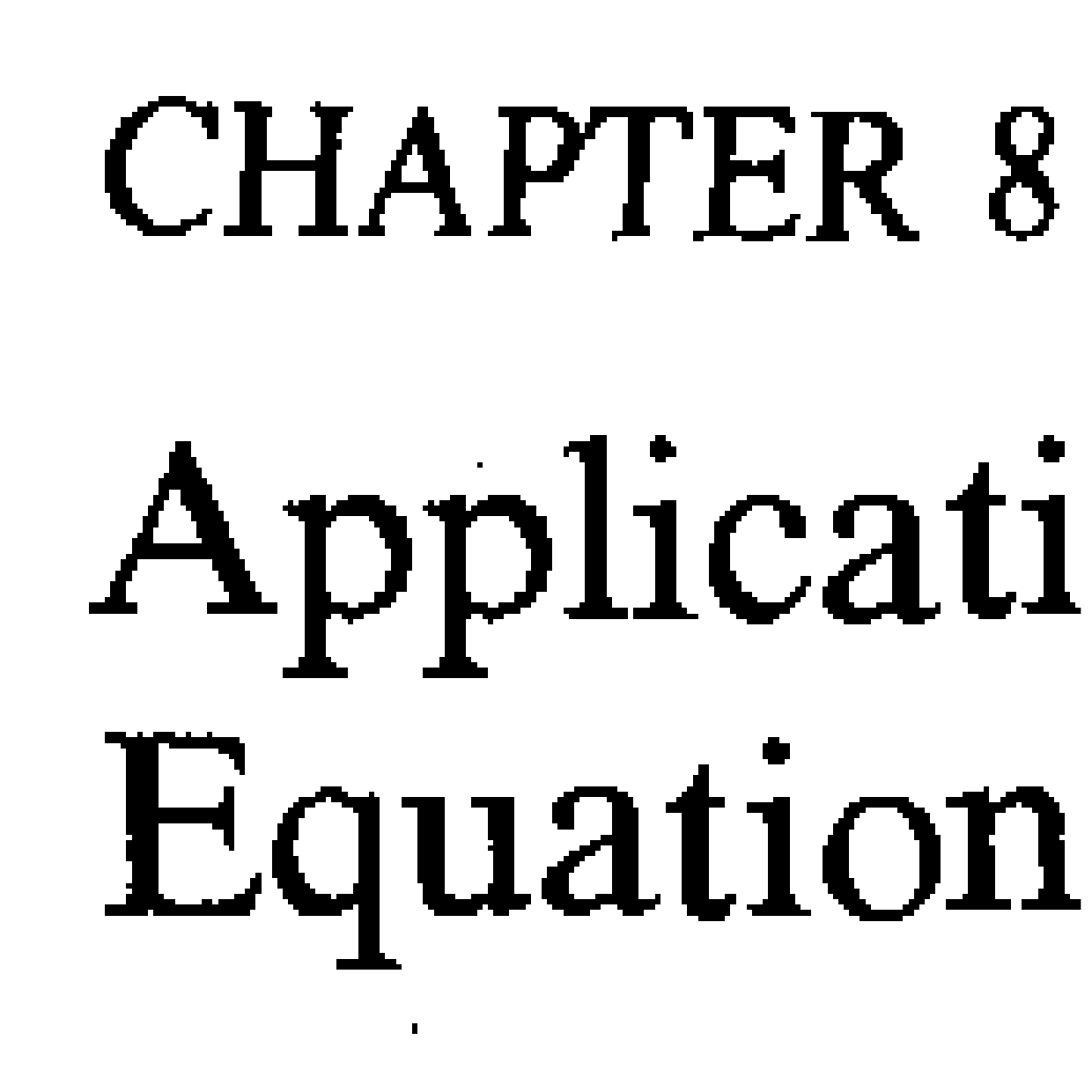}}\hspace{0.2cm}
    \fbox{\includegraphics[scale=0.28]{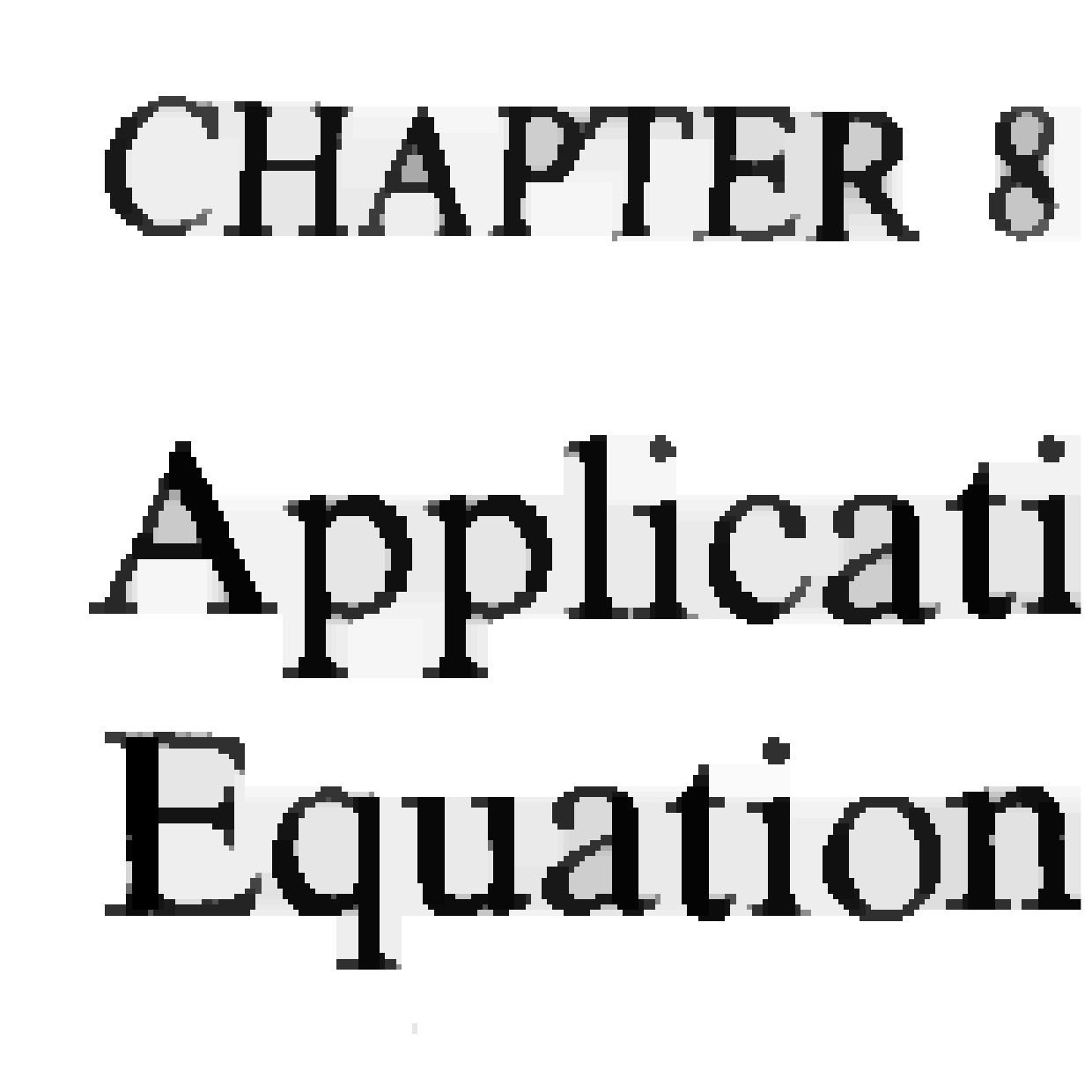}}\hspace{0.2cm}
    \fbox{\includegraphics[scale=0.28]{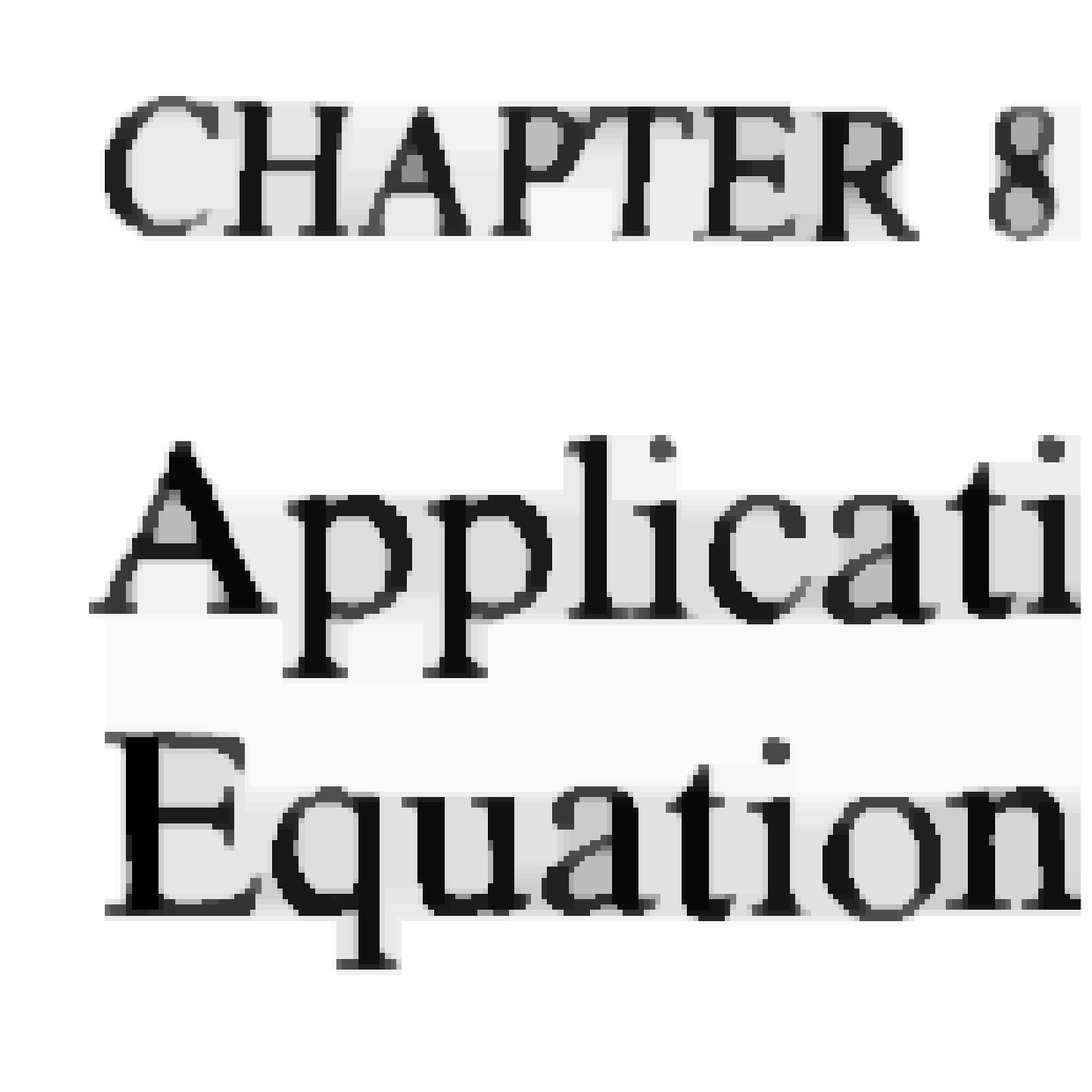}}\hspace{0.2cm}
    \fbox{\includegraphics[scale=0.28]{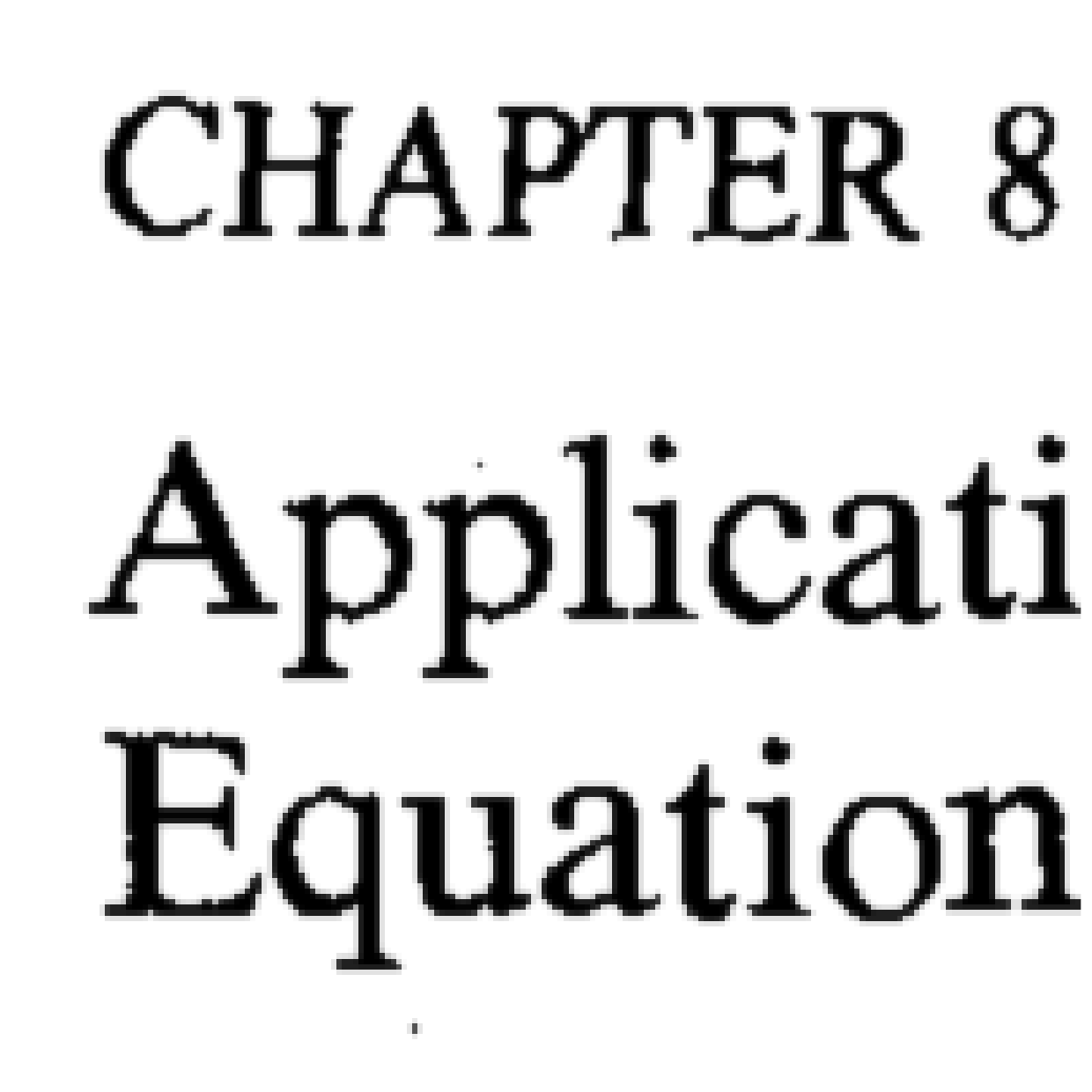}}    
    \\[0.2cm]
    \fbox{\includegraphics[scale=0.28]{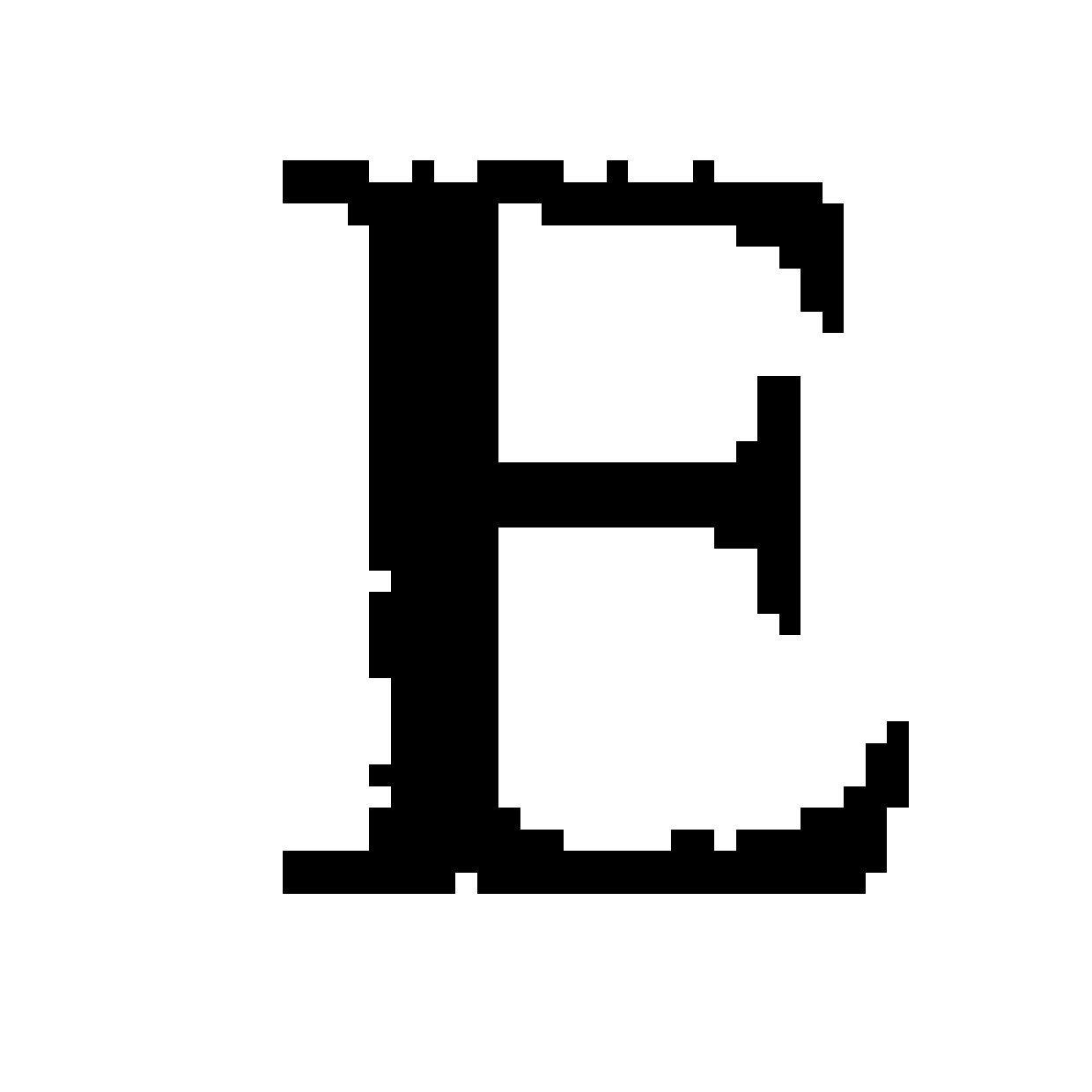}}\hspace{0.2cm}
    \fbox{\includegraphics[scale=0.28]{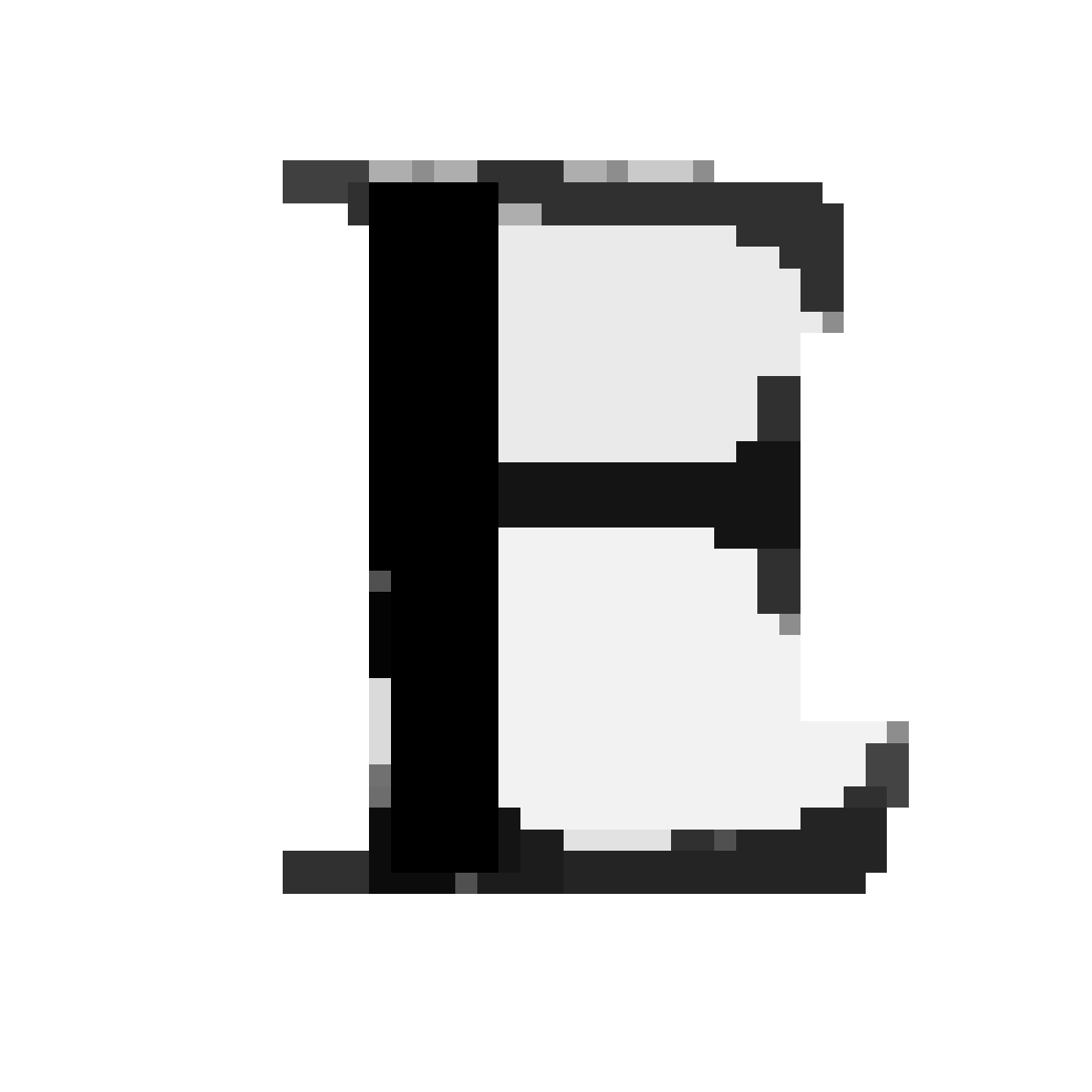}}\hspace{0.2cm}
    \fbox{\includegraphics[scale=0.28]{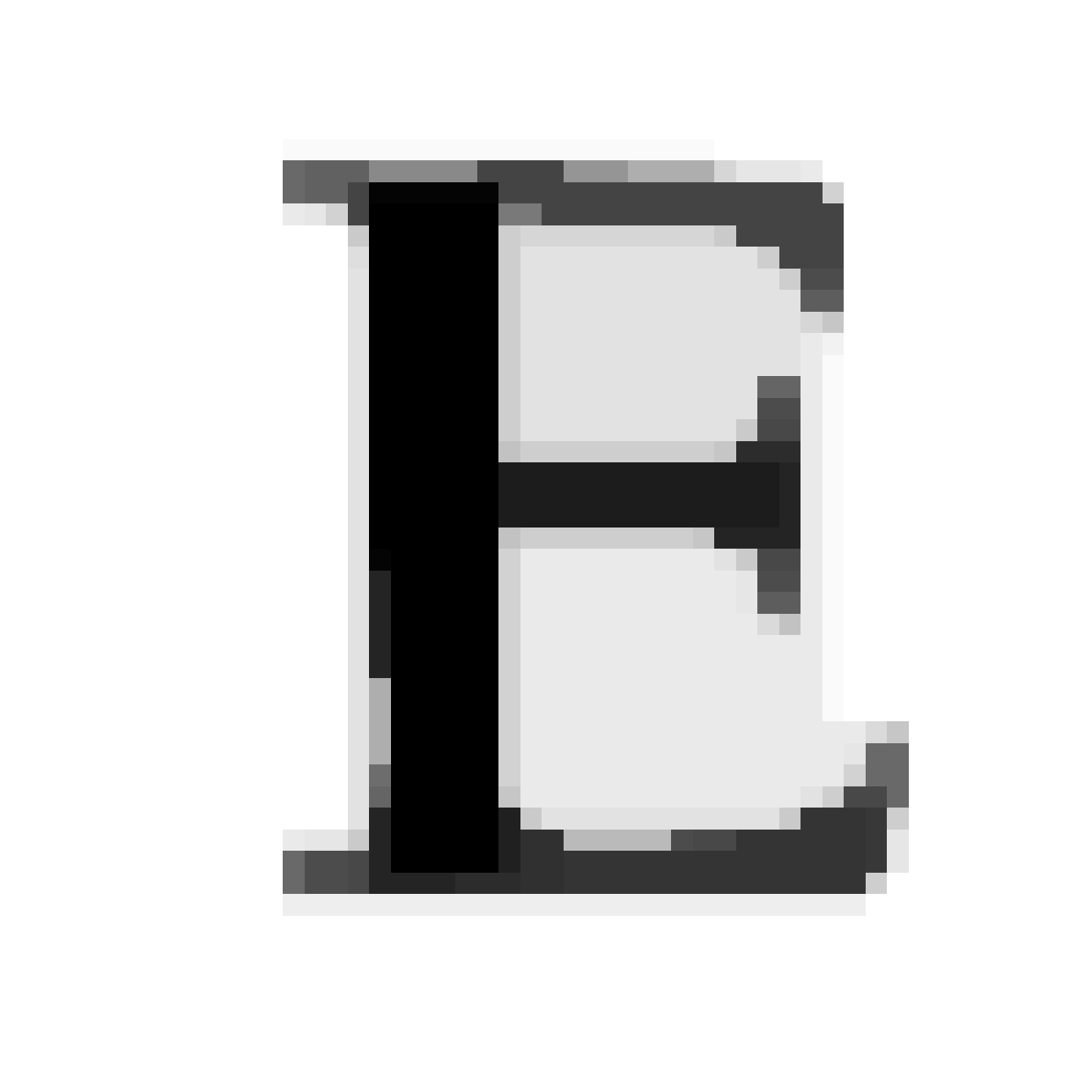}}\hspace{0.2cm}
    \fbox{\includegraphics[scale=0.28]{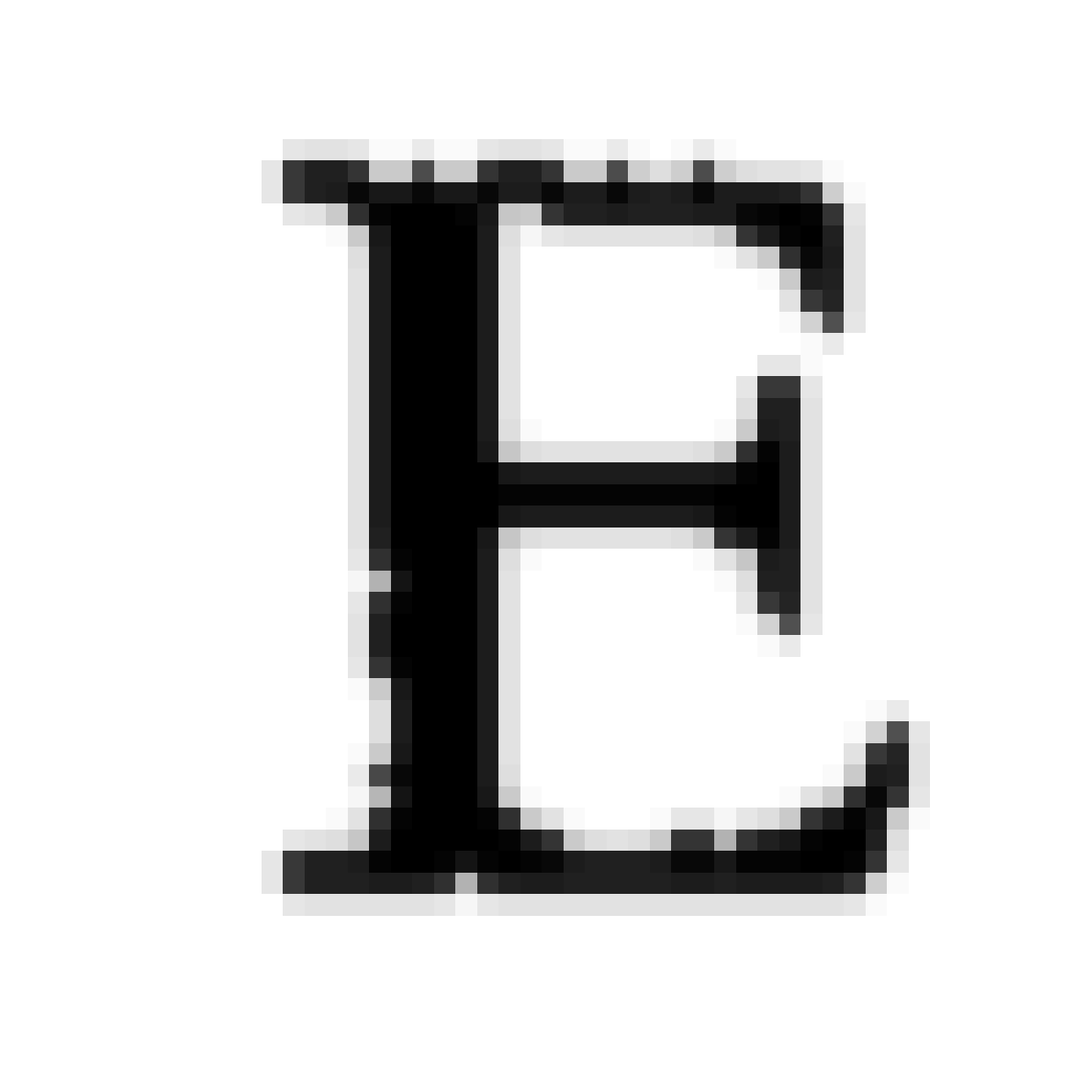}}       
    \caption{In columns: (a) initial data -- two binary images of the scanned text, (b) solutions to the equation (\ref{eq_anisoTV}), (c) solutions to the equation (\ref{eq_diffanisoTV}), (d) solutions to the linear diffusion equation.}
    \label{fig5}
    \end{center}
\end{figure}

\begin{figure}[h!]  
   \begin{center}
    \setlength{\fboxsep}{0pt}
    \setlength{\fboxrule}{0.5pt}
    \fbox{\includegraphics[scale=0.28]{fig2_letters-eps-converted-to.pdf}}\hspace{0.2cm}
    \fbox{\includegraphics[scale=0.28]{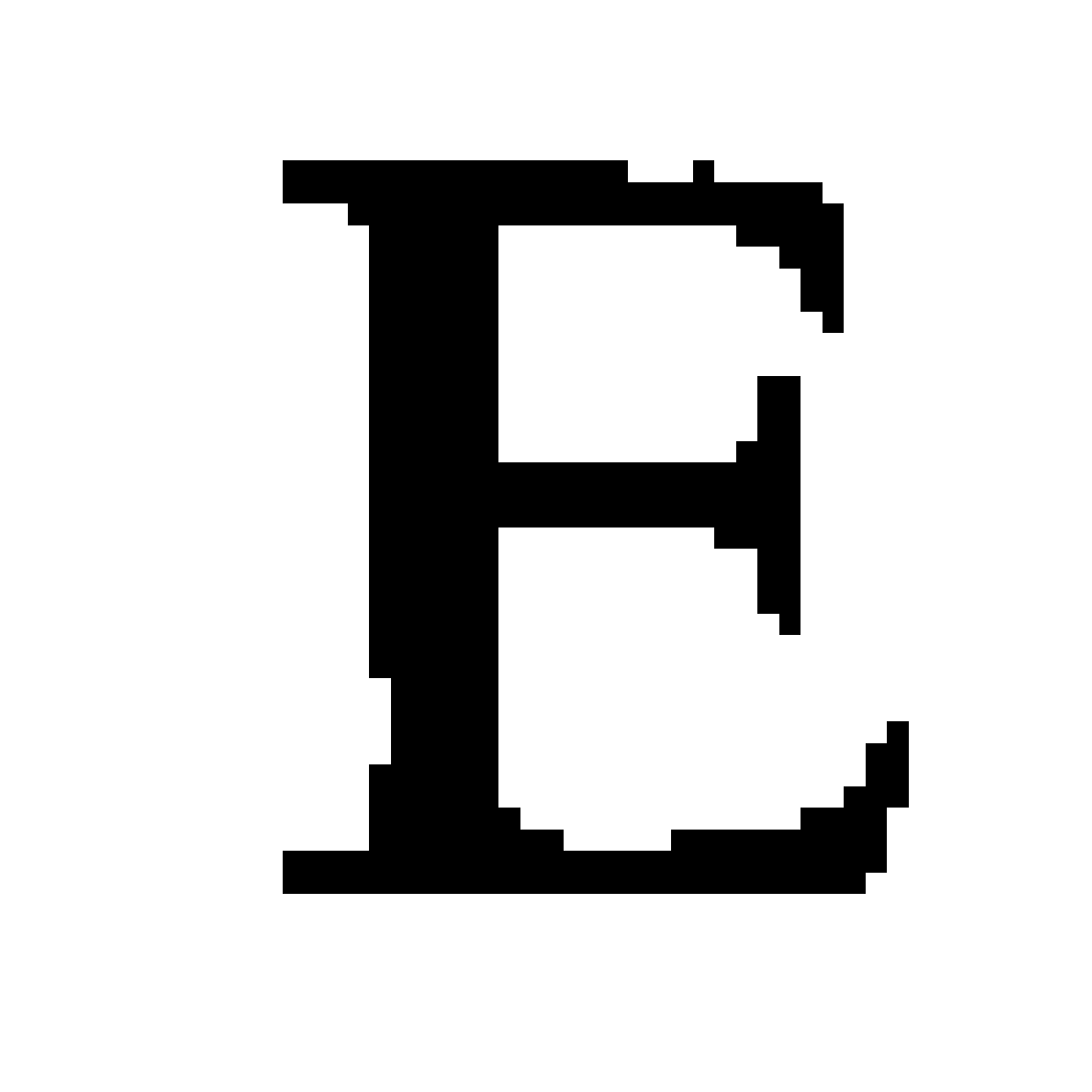}}\hspace{0.2cm}
    \fbox{\includegraphics[scale=0.28]{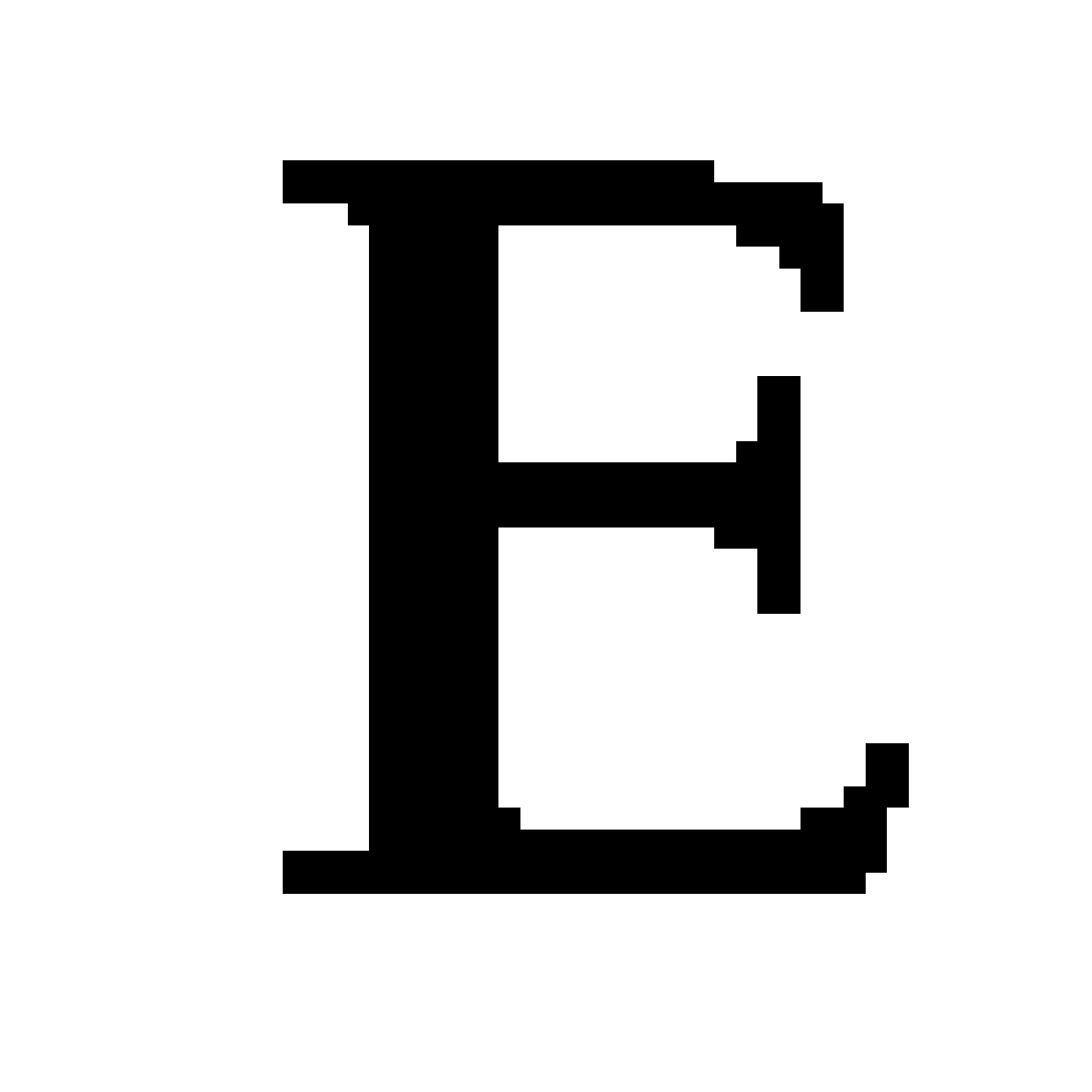}}\hspace{0.2cm}
    \fbox{\includegraphics[scale=0.28]{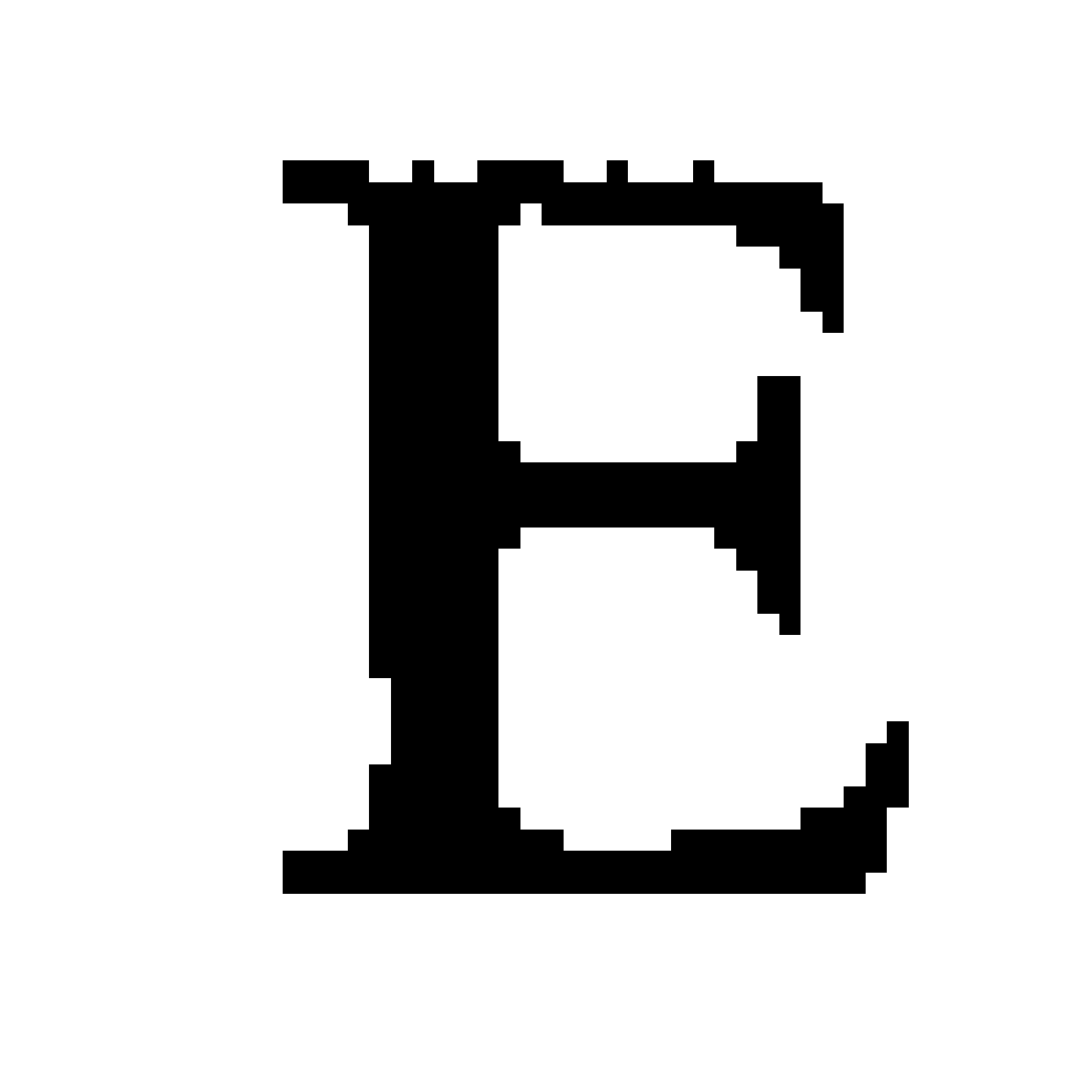}} \\[0.2cm]
    \fbox{\includegraphics[scale=0.28]{fig2_letters-eps-converted-to.pdf}}\hspace{0.2cm}
    \fbox{\includegraphics[scale=0.28]{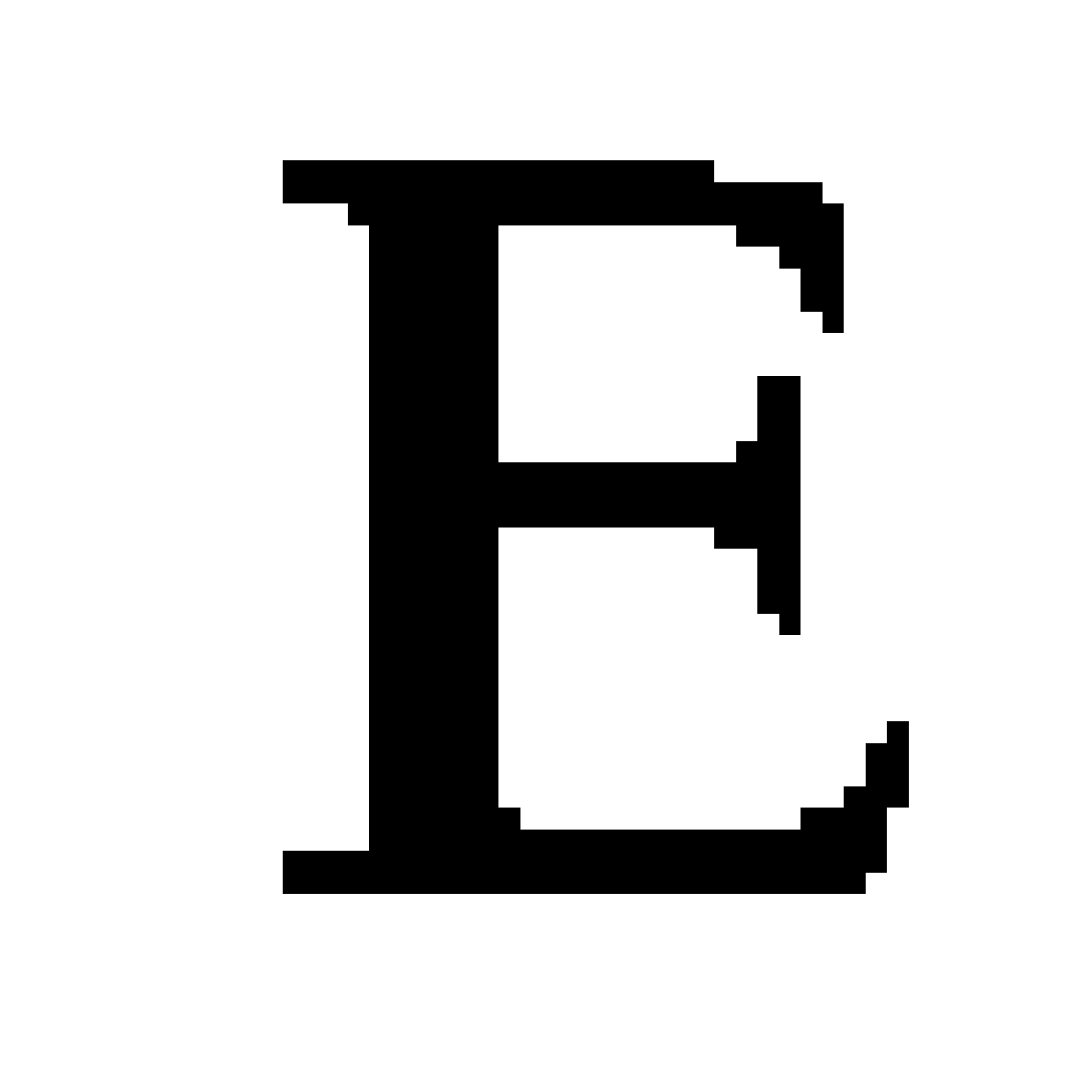}}\hspace{0.2cm}
    \fbox{\includegraphics[scale=0.28]{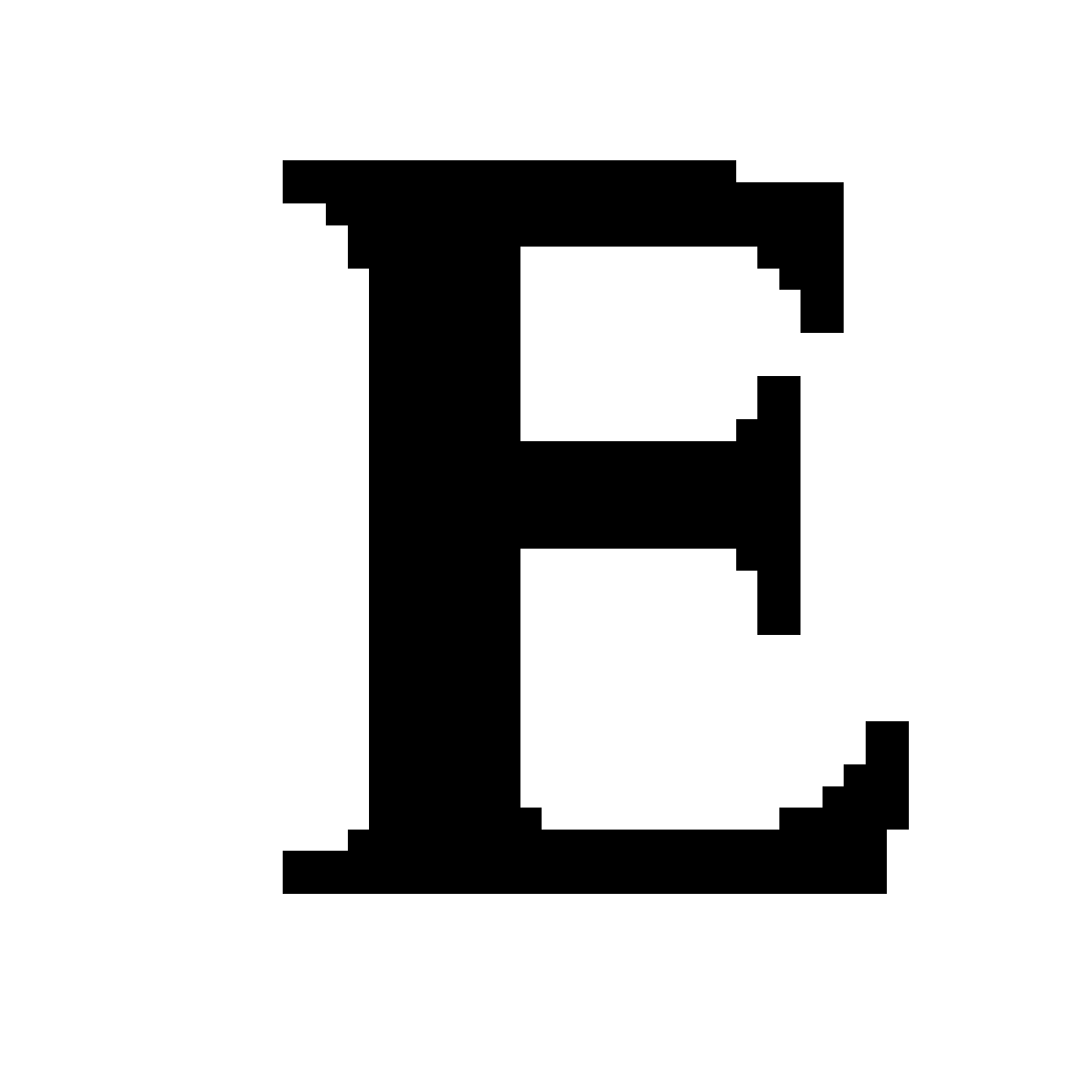}}\hspace{0.2cm}
    \fbox{\includegraphics[scale=0.28]{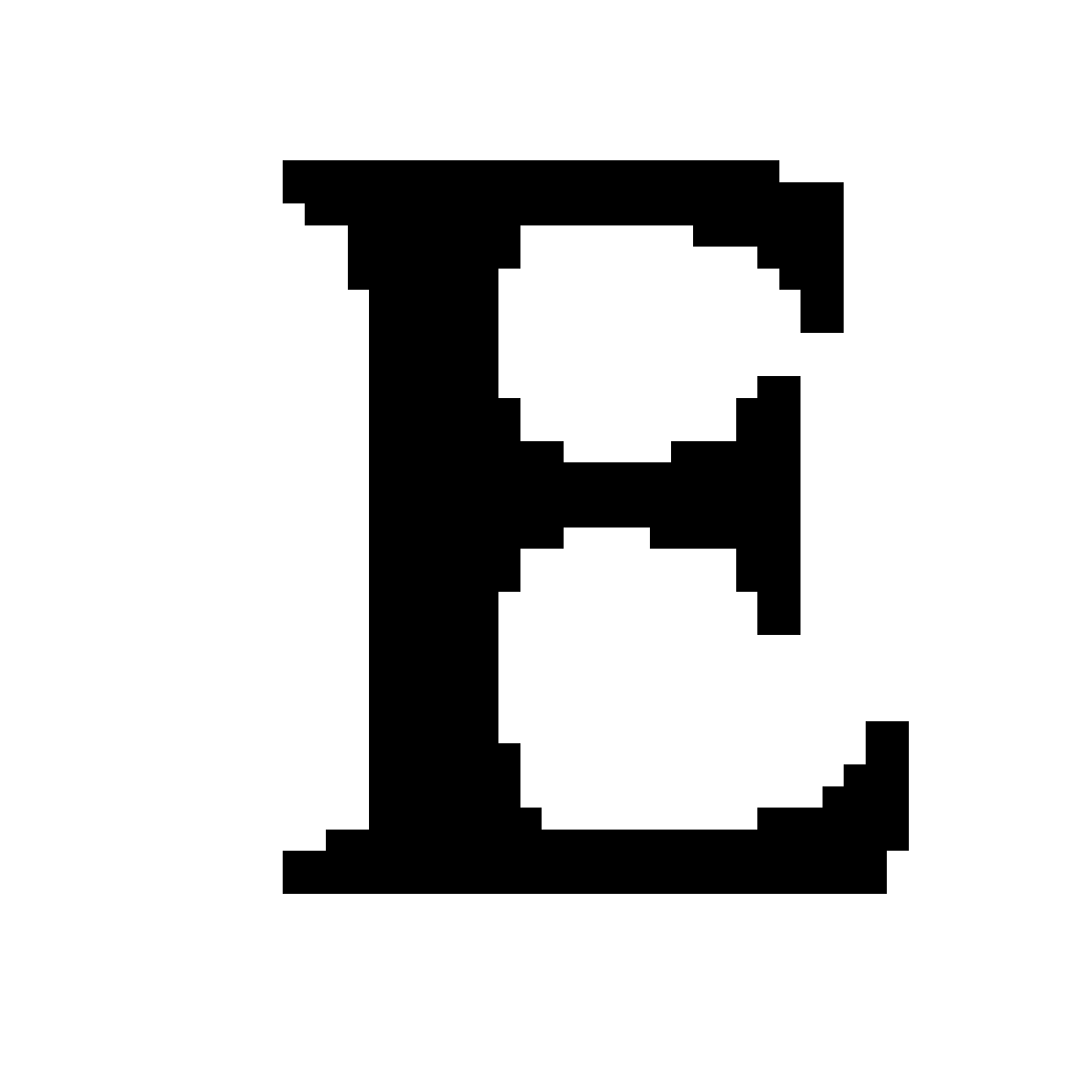}}
    \caption{Results obtained by thresholding of images presented in the lower row of Figure~\ref{fig5}
    on~the~level~$-10$ (upper row) and  on~the~level~$-5$ (lower row).}
    \label{fig6}
    \end{center}
\end{figure}

In the last experiment, we were testing a possible application of the
anisotropic total variation flow equations (\ref{eq_anisoTV}) and
(\ref{eq_diffanisoTV}) to solve the real problem of improving the quality of the
scanned text. In this experiment, we were considering two binary images
presented in the first column of Figure~\ref{fig5} with values scaled to
$\{-50,0\}$. Images in the second and third column of Figure~\ref{fig5}
correspond to numerical solutions to the equations (\ref{eq_anisoTV}) and
(\ref{eq_diffanisoTV}), respectively, for parameters $\gamma = 10^{-2}$, $\beta
= 1$, $\delta t = 1$, $\tau = 8^{-1}$, $m=15$, and with images in the first
column of Figure~\ref{fig5} as initial data.
For comparison, in the last column of Figure \ref{fig5}, we present numerical
solutions to the linear diffusion flow (the equation (\ref{eq_diffanisoTV}) for
$\gamma = 10^{-2}$, $\beta=0$, $\delta t = 1$, $m=15$) with the same initial
data. We observe that in fact equation (\ref{eq_diffanisoTV}) represents the 
interplay between an anisotropic total variation flow and the linear diffusion.
It allows to fill corrupted parts of letters and at the same time slightly blur
their boundaries. We notice that these properties are also visible in the
results of
the experiment with the image $f_{S_4}$, presented in the last columns of
Figures \ref{fig2} and \ref{fig3}.

In Figure~\ref{fig6}, we present results obtained by thresholding  images in
the lower row of Figure~\ref{fig5} on the level $-10$ and $-5$, respectively. We
see that application of  equation (\ref{eq_anisoTV}) gives basically better
results, however in the case when larger parts of the letters are corrupted, the
properties of  equation (\ref{eq_diffanisoTV}) may be useful. In general, we 
infer from the experiments we performed that both total variation flow models
analysed in this paper provide better results when applied to a class of
real problem, than
the standard linear diffusion equation.

\section{Appendix}

Formula (\ref{wzor2}) must be modified in order to accommodate the boundary
conditions. This is done below.

\begin{proposition}\label{pr2}
Formula (\ref{wzor4}) below yields a weak solution to (\ref{eq_anisoTV}) in
$\bR^2$
with the data 
\begin{equation*}
 u_0(x_1,x_2) = (x_1^2+x_2^2 -2R^2)\chi_{B(0,R)}(x_1,x_2)\in L^2(\bR^2).
 \end{equation*}
 in the sense specified in Theorem \ref{tw1}. Moreover, the equation is
satisfied in $\bR^2$ in a
pointwise manner with the exception of a one dimensional set and the solution
is Lipschitz continuous, but not $C^1$. 
\end{proposition}

{\it Proof.} Formula (\ref{wzor2}) shows the creation of a square facet and
ruled
surfaces over strips $|x_1| \le \xi(t)$ and   $|x_2| \le \xi(t)$. Now, we have
to take into account their interaction with the boundary of the ball
$x_1^2+x_2^2 \le R^2$. The result is region $\Omega(t)$, where $u$ is different from zero. This set is defined as follows,
$\Omega(t) = B(0,R) \cap (-L(t), L(t))$,
where $L(t) = \sqrt{R^2- \xi^2(t)}$.

We shall see that the solution gets extinct, when the square
facet hits the plane $u=0$ at $t= t_1$. This is why for
$t\in[0, t_1)$, we set,
\begin{equation}\label{wzor4}
 u(x,t) =\left\{
\begin{array}{ll}
2h(t) 
& |x_1|, |x_2| \le \xi(t), (x_1,x_2)\in \Om(t),\\
h(t) + x_2^2 -2R^2 & |x_1| \le \xi(t), \xi(t)< |x_2| \le \sqrt{R^2-\xi^2(t)}, (x_1,x_2)\in \Om(t)\\
0& |x_1| \le \xi(t),  |x_2| > \sqrt{R^2-\xi^2(t)}, (x_1,x_2)\in \Om(t),\\
h(t) + x_1^2 -2R^2  & |x_2| \le \xi(t), \xi(t)< |x_1| \le
\sqrt{R^2-\xi^2(t)}, (x_1,x_2)\in \Om(t),\\
0& |x_2| \le \xi(t),  |x_1| >\sqrt{R^2-\xi^2(t)},\\
x_1^2 + x_2^2 - 2R^2 & |x_1|, |x_2| > \xi(t), (x_1,x_2)\in \Om(t),\\
0&  (x_1,x_2)\not\in \Om(t).
\end{array}
\right.
\end{equation}
This formula is valid up to $2\xi^2(t_1)= R^2$, i.e.
$t_1 = \frac{\sqrt 2}{6} R^3$.

Calculating $\nabla u$ is easy, but we have to modify $\cL(\nabla u)$. Namely,
we set,
\begin{equation*}
 \cL(\nabla u)_1 =\left\{
\begin{array}{ll}
\frac{x_1}{\xi(t)} & |x_1|\le \xi(t), |x_2| \le \sqrt{R^2-\xi^2(t)} ,\\
\sgn x_1 & |x_1| > \xi(t),\\
0 & |x_2| >\sqrt{R^2-\xi^2(t)} ;
\end{array}
\right.
\end{equation*}
\begin{equation*}
 \cL(\nabla u)_2 =\left\{
\begin{array}{ll}
\frac{x_2}{\xi(t)} & |x_2|\le \xi(t), |x_1| \le \sqrt{R^2-\xi^2(t)} ,\\
\sgn x_2 & |x_2| > \xi(t),\\
0 & |x_1| >\sqrt{R^2-\xi^2(t)} .
\end{array}
\right.
\end{equation*}
We notice that vector field $\cL(\nabla u)$ has jump discontinuities,
nonetheless its distributional divergence is in $L^2_{loc}$ and has the desired
properties. It is now easy to check that $u$ satisfies (\ref{eq_anisoTV})
pointwise except a two-dimensional set in $\bR^2\times \bR_+$.
We note the discontinuity of $u_t$ is responsible for the creation of the
two dimensional facet and its growth. \qed

\subsection*{Acknowledgement} The work has been supported by the MN grant IdP2011 000661.
A part of the research for this paper was performed while PR was visiting IMA,
University of Minnesota, whose hospitality is acknowledged.

%

\end{document}